\documentclass[11pt, a4paper]{amsart}
\usepackage{epsfig}
\usepackage{amsmath}
\usepackage{stackrel, amssymb}
\usepackage{amscd}
\usepackage{latexsym}
\usepackage{tabularx}
\usepackage{a4wide}
\usepackage[usenames]{color}
\usepackage{subfigure}
\usepackage{url}
\usepackage{verbatim}
\usepackage{mathtools}
\usepackage{circuitikz}
\usepackage{mathabx}
\usepackage{tikz-cd}

\usepackage{enumitem} 
\setlist[enumerate]{label=\textnormal{(\arabic*)}}

\newcommand{\image}{\mathrm{im}}
\ctikzset{bipoles/resistor/height=0.15}
\ctikzset{bipoles/resistor/width=0.4}

\usepackage[
colorlinks,  citecolor=blue,
pdfauthor={Omid Amini, Eduardo Esteves}, 
pdftitle={tropicalization-linear-series-polymatroids},
pdfstartview ={FitV},
]{hyperref}
\usepackage{tikz, float} 
\usetikzlibrary {positioning}
\usetikzlibrary{calc,decorations.markings}
\usetikzlibrary{shapes,snakes}
\numberwithin{equation}{section}
\tikzstyle{Cwhite}=[scale = .8,circle, fill = white, minimum size=3mm] 
\tikzstyle{Cgray}=[scale = .4,circle, fill = gray, minimum size=3mm] 
\tikzstyle{Cblack2}=[scale = .4,circle, fill = black, minimum size=5mm] 
\tikzstyle{Cblack}=[scale = .7,circle, fill = black, minimum size=3mm]
\tikzstyle{C0}=[scale = .9,circle, fill = black!0, inner sep = 0pt, minimum size=3mm]
\tikzstyle{C1}=[scale = .7,circle, fill = black!0, inner sep = 0pt, minimum size=3mm]
\tikzstyle{Cred}=[scale = .4,circle, fill = red, minimum size=3mm] 
\usepackage[matrix,arrow,tips,curve]{xy}
\usepackage{pb-diagram,pb-xy}
\usepackage{verbatim}
\usepackage{mathrsfs}
\usepackage{color}
\usepackage[leqno]{amsmath}
\usepackage{scalerel}
\usepackage{marginnote}
\RequirePackage{thmtools}   
\RequirePackage{varioref}

\RequirePackage[capitalize, noabbrev]{cleveref}
\newtheorem{thm}{Theorem}[section]
\newtheorem{lemma}[thm]{Lemma}
\newtheorem{prop}[thm]{Proposition}

\theoremstyle{definition}
\newtheorem{defii}[thm]{Definition}
\newenvironment{defi}
  {\pushQED{\qed}\defii}
  {\popQED\enddefii}
\newtheorem{remm}[thm]{Remark}
\newenvironment{remark}
  {\pushQED{\qed}\remm}
  {\popQED\endremm}

  \newtheorem{conv}[thm]{Convention}
  \newenvironment{convention}
    {\pushQED{\qed}\conv}
    {\popQED\endconv}  
\newtheorem{question}[thm]{Question}

\numberwithin{equation}{section}

\newcommand{\K}{\ssub{\mathbf K}!}
\renewcommand{\k}{\ssub{\mathbf k}!}
\newcommand{\R}{\ssub{\mathbb R}!}
\renewcommand{\P}{\ssub[-2pt]{\mathbf P}!}
\newcommand{\Q}{\ssub{\mathbf Q}!}
\newcommand{\B}{\ssub{\mathbf B}!}

\newcommand{\F}{\ssub{\mathbf F}!}

\newcommand{\M}{\widehat{\mathcal M}}
\newcommand{\subM}{\widecheck{\mathcal M}}
\newcommand{\Z}{\mathbb Z}
\renewcommand{\:}{\colon}
\newcommand{\C}{{\mathscr C}}

\newcommand{\g}{\mathfrak{g}}
\newcommand{\red}{\ssub{\mathrm{red}}!}

\newcommand{\val}{\ssub{\mathfrak v}!}
\newcommand{\trop}{\mathrm{trop}}
\newcommand{\an}{{^\mathrm{an}}}

\newcommand{\innone}[2]{\langle#2\rangle_{_{\hspace{-.02cm}#1}}}

\NewDocumentCommand{\ssub}{O{0pt} O{.85} m t! e{_^}}{
  #3%
  \IfValueT{#5}{
    \IfBooleanTF{#4}{\sb{\hspace{#1}\scaleobj{#2}{#5}}}{\sb{#5}}
  }
  \IfValueT{#6}{
  \IfBooleanTF{#4}{\sp{\hspace{#1}\scaleobj{#2}{#6}}}{\sp{#6}}
}
}

\ExplSyntaxOn
\NewDocumentCommand{\tossub}{o o m}{
  \expandafter\let\csname old\cs_to_str:N #3\endcsname#3
  \renewcommand#3%
  {\ssub[#1][#2]{\csname old\cs_to_str:N #3\endcsname}}
}
\ExplSyntaxOff

\newcommand{\stable}{\mathscr A}

\newcommand{\rquot}[2]{#1\big/#2}
\newcommand{\adfun}{\ssub{\mathscr F}!}

\newcommand{\cd}{\mathrm{cd}}
\newcommand{\dd}{\mathit{d}}

\newcommand{\VS}{\ssub{\mathbf U}!}

\tossub\mu
\tossub\nu
\tossub\delta
\tossub\Gamma
\tossub\eta
\tossub\stable
\tossub\C
\tossub\iota
\tossub\rho
\tossub\g
\tossub\red
\tossub\cd
\tossub\dd
\tossub\varpi
\tossub\lambda
\tossub\epsilon
\tossub\theta
\tossub\zeta
\newcommand{\sschi}{\ssub{\chi}!}

\newcommand{\ssGamma}{\ssub[-2pt]{\Gamma}!}
\newcommand{\norm}[2]{\ssub{|#1|}!_{#2}}

\newcommand{\ssI}{\ssub{I}!}
\newcommand{\ssJ}{\ssub{J}!}

\newcommand{\sspi}{\ssub{\pi}!}
\newcommand{\ssG}{\ssub{G}!}
\newcommand{\ssD}{\ssub{D}!}

\newcommand{\ssW}{\ssub[-.5pt]{\mathrm W}!}
\newcommand{\ssn}{\ssub{n}!}

\newcommand{\ssq}{\ssub{q}!}
\newcommand{\ssh}{\ssub{h}!}

\newcommand{\ssfhat}{\ssub{\hat f}!}
\newcommand{\ssf}
{\ssub{f}!}
\newcommand{\ssE}{\ssub{E}!}

\newcommand{\sse}{\ssub{e}!}

\newcommand{\ssa}{\ssub{a}!}
\newcommand{\sst}{\ssub{t}!}
\newcommand{\ssc}{\ssub{c}!}
\newcommand{\ssX}{\ssub{X}!}

\newcommand{\ml}{\ssub{\mathrm{ml}}!}   
\newcommand{\ts}{\ssub{\mathrm{ts}}!}   
\newcommand{\ssF}{\ssub{F}!}
\newcommand{\N}{\mathbb N}
\newcommand{\ev}{\ssub{\mathrm{ev}}!}

\newcommand{\filt}{\ssub{\mathrm{F}}!}
\newcommand{\filtind}{\mathrm{F}}

\renewcommand{\setminus}{\smallsetminus}
\renewcommand{\emptyset}{\varnothing}

\newcommand{\Div}{\operatorname{Div}}
\renewcommand{\div}{\mathrm{div}}
\newcommand{\E}{\mathbb E}
\newcommand{\ord}{\mathrm{ord}}
\newcommand{\he}{\mathrm{h}}  
\newcommand{\te}{\mathrm{t}}  

\newcommand{\Rat}{\mathrm{Rat}}

\newcommand{\twist}{\tau}

\newcommand{\twistd}[1]{\ssub{\twist}!^{#1}}

\newcommand{\slztwist}[2]{\partial_{_{ \hspace{-.07cm}#1}}^{^{\hspace{-.02cm}#2}}\hspace{-.1cm}}

\newcommand{\rest}[1]{\ssub{\raisebox{-1pt}{$\vert$}}!_{#1}}

\newcommand{\proj}[1]{\theta_{_{ \hspace{-.07cm}#1}}}

\newcommand{\cl}{\ssub{\mathscr D}!}  
\newcommand{\Prin}{\operatorname{Prin}}
\newcommand{\Pic}{\operatorname{Pic}}

\newcommand{\ol}{\overline}

\newcommand{\st}{\bigm|} 

\newcommand{\Fun}{\mathscr E} 
\tossub\Fun 

\newcommand{\id}{\mathrm{Id}}
\tossub\id
\newcommand{\zero}{0}
\newcommand{\one}{\mathbf 1}
\tossub\zero
\tossub\one

\newcommand{\ssrho}{\ssub{\rho}!}
\newcommand{\ssphi}{\ssub{\phi}!}
\newcommand{\sspsi}{\ssub{\psi}!}

\tossub \chi
\newcommand{\Vor}{\ssub{\mathrm{Vor}}!}
\tossub{\varepsilon}
\newcommand{\varX}{\ssub{\mathbf X}!}
\newcommand{\varC}{\mathbf C}

\newcommand{\varD}{\mathbf D}
\newcommand{\varH}{\mathbf H}
\newcommand{\varM}{\ssub{\mathbf M}!}
\newcommand{\varN}{\mathbf N}
\newcommand{\varm}{\mathbf m}
\newcommand{\varR}{\mathbf R}
\newcommand{\varL}{\mathbf L}

\newcommand{\varf}{\mathbf f}
\newcommand{\varg}{\mathbf g}
\newcommand{\varG}{\ssub{\mathbf G}!}
\newcommand{\varGr}{\mathbf{Gr}}
\newcommand{\varY}{\mathbf Y}
\newcommand{\varZ}{\mathbf Z}
\newcommand{\varT}{\mathbf T}

\newcommand{\valgroup}{\Lambda}

\newcommand{\frakX}{\ssub{\mathfrak X}!}

\newcommand{\FC}{\ssub[-.5pt]{\mathcal FC}!}

\definecolor{cadmiumgreen}{rgb}{0.0, 0.42, 0.24}
\definecolor{darkred}{rgb}{.6,0,0}
\definecolor{byzant}{rgb}{0.74, 0.2, 0.64}
 \definecolor{pblue}{rgb}{0.11, 0.22, 0.73}
\definecolor{pgreen}{rgb}{0.0, 0.65, 0.58}
\definecolor{aqua}{rgb}{0.0, 1.0, 1.0}
\definecolor{lblue}{rgb}{0.0, 0.55, 1.0}

\makeatletter
\@namedef{subjclassname@2020}{\textup{2020} Mathematics Subject Classification}
\makeatother

\newcommand{\msc}[1]{\href{http://www.ams.org/msc/msc2020.html?t=&s=#1}{#1}}

\title{Tropicalization of linear series and tilings by polymatroids}
 \author{Omid Amini}
 \author{Eduardo Esteves}
 \date{\today}
  \address{CNRS - CMLS, \'Ecole Polytechnique, Palaiseau, France}
\email{\href{omid.amini@polytechnique.edu}{omid.amini@polytechnique.edu}}

\address{Instituto Nacional de Matem\'atica Pura e Aplicada, Rio de Janeiro, Brazil}
\email{\href{esteves@impa.br}{esteves@impa.br}}

\begin{document}

\begin{abstract}
 We show that tropicalization of linear series on curves gives rise to two-parameter families of tilings by polymatroids, with one parameter arising from the theory of divisors on tropical curves and the other from the reduction of linear series of rational functions in non-Archimedean geometry. In order to do this, we introduce a general framework that produces tilings of vector spaces and their subsets by polymatroids. We furthermore show that these tilings are regular and relate them to work by Kapranov and Lafforgue on Chow quotients of Grassmannians. 
\end{abstract}

\maketitle
\setcounter{tocdepth}{1}

\tableofcontents


\section{Introduction}

Understanding degenerations of linear series, that is, line bundles and their spaces of sections, has been a long-standing open problem in the theory of algebraic curves.  The results of this paper form the combinatorial and polyhedral foundation of a program aimed at addressing this problem.

To formulate our results, we will use the framework of tropical geometry, that provides a setting to study degenerations of algebraic varieties and geometric structures over them by enriching the classical techniques in algebraic geometry with polyhedral geometry. In the case of curves and line bundles, tropicalization gives rise to a metric graph, also known as tropical curve, and a linear equivalence class of divisors on it. Over the past fifteen years, tropical methods have been quite successful in the study of the geometry of curves and their moduli spaces. We refer to the survey papers~\cite{BJ} and~\cite{JP} and the references therein for an introduction to the topic and a sample of results.

In a degenerating family of linear series on curves that tropicalize to a given metric graph, we encode the combinatorics of the degeneration by a collection of \emph{admissible divisors} in the equivalence class of divisors given by the tropicalization of the underlying line bundles. Furthermore, we associate to each of the admissible divisors a \emph{limit space of sections} on the limit stable curve of the family in the Deligne--Mumford compactification of the moduli space of smooth curves. In this way, we obtain a collection of hybrid pairs consisting each of an admissible divisor on the tropical curve and a space of sections on the stable curve. 
 
Our main result shows that on the one side, the collection of admissible divisors on the metric graph gives rise to a tiling of an Euclidean space by polytopes, and on the other side, the collection of spaces of sections associated to the admissible divisors leads to a tiling of a simplex by polytopes. Quite remarkably, the two tilings are compatible with each other in the sense that, given any two positive real numbers as coefficients, the collection of polytopes obtained by taking the linear combination of the two tiles associated to each admissible divisor again gives rise to a tiling of an Euclidean space.  Furthermore, the polytopes appearing in this two-parameter family of tilings are all polymatroids. Figures~\ref{fig:tiling3} and \ref{fig:mixed-tiling} give two examples of such tilings.

We identify the polytope associated with an admissible divisor on a metric graph, when that divisor arises from a degeneration of line bundles, with the semistability polytope associated with the limit line bundle on a semistable model of the curve (more precisely, to the corresponding torsion-free, rank-one sheaf on the stable curve itself). The semistability polytope is the union of the chambers parameterizing the semistability conditions that the limit sheaf satisfies. Even though the various semistability conditions and moduli spaces associated to them have been considered by several authors \cite{OS79, Cap, MRV, KP}, including the second named author of the present work \cite{Esteves01}, to our knowledge the semistability polytopes have not been studied before, let alone tilings by the semistable polytopes of admissible divisors, corresponding to different limits of line bundles.

In equal characteristic, the tilings of simplices by polymatroids we consider are related to the pioneering  work by Kapranov~\cite{Kap-chow} on Chow quotients and by Lafforgue~\cite{Laf99} on compactifications of $\rquot{\mathrm{PGL}_r^{n+1}}{\mathrm{PGL}_r}$. Making this link precise requires the results we prove on reduction of spaces of sections in Sections~\ref{sec:tropicalization} and~\ref{sec:tilings-simplex-mixed}. This will be discussed in Section~\ref{sec:chow-quotients}, using more  recent work by Giansiracusa and Wu \cite{GW22} that extends both the work by Kapranov and that by Lafforgue. In our setting, each tile in the simplex is associated to an admissible divisor,  tying up the tiling by semistability polytopes to the tiling by the degenerations of the space of sections, leading to a hybrid setting that captures both. This way, the tiling of the simplex acquires geometric meaning for the nodal curve itself.

In order to prove our results, we introduce a combinatorial framework that produces tilings of real vector spaces and their subsets by polymatroids. We focus here on results that are directly related to the degeneration problem for linear series. The setup is however quite general and can be used in the study of other geometric problems, as we show in a sequel. 

In our use of this combinatorial framework, we prove results in tropical geometry, both in the theory of divisors on metric graphs and in the study of reduction of spaces of sections of line bundles in non-Archimedean geometry, that we hope to be of independent interest in further applications of tropical techniques in the study of the geometry of curves and their moduli spaces.

  In a sequel~\cite{AEG23},  using the approach introduced here, we describe limits of canonical series on smooth proper curves degenerating under general directions to a nodal curve that is general for its topology, and construct a parameter space for them, generalizing to any topology pioneering work by Eisenbud--Harris \cite{EH87}, Esteves--Medeiros \cite{EM} and Esteves--Salehyan \cite{ES07}.

\smallskip

In the remaining of this introduction, we give an overview of our results and expand on related work.

\subsection{Modular pairs and base polytopes} Let $V$ be a finite nonempty set. A function $\mu \colon \ssub{2}!^V \to \R$ is called \emph{supermodular} if $\mu(\emptyset) =0$ and for each pair of subsets $\ssI_1,\ssI_2\subseteq V$, we have the inequality 
\[
\mu(\ssI_1)+\mu(\ssI_2)\leq \mu(\ssI_1\cup \ssI_2)+ \mu(\ssI_1\cap \ssI_2).
\]
The \emph{adjoint} $\mu!^*$ of each function $\mu$ on $\ssub{2}!^V $ is defined by setting $\mu!^*(I) = \mu(V) - \mu(\ssI^c)$, with $\ssI^c=V\setminus I$. If $\mu$ is supermodular, the adjoint $\mu!^*$ is \emph{submodular}, meaning that $\mu!^*(\emptyset)=0$, and it verifies the set of inequalities 
\[
\mu!^*(\ssI_1)+\mu!^*(\ssI_2)\geq \mu!^*(\ssI_1\cup \ssI_2)+ \mu!^*(\ssI_1\cap \ssI_2)
\]
for each $\ssI_1, \ssI_2 \subseteq V$. We refer to $(\mu, \mu!^*)$ as a \emph{modular pair}.
 
Denote by $H=\R^V$ the real vector space of functions $q\colon V \to \R$. For each $q\in H$ and subset $I\subseteq V$, set $q(I) = \sum_{v\in I}q(v)$. It is easy to see that the pair $(q,q)$ is modular. Other examples are given by matroids. For each matroid $M$ with $V$ as ground set, the rank function $\rho!_M \colon \ssub{2}!^V \to \R$ is submodular. Setting $\mu!_M \coloneqq \rho!_M^*$ and $\mu!_M^* = \rho!_M$, we get a modular pair $(\mu!_M, \mu!_M^*)$. 

To each modular pair $(\mu, \mu!^*)$, we associate the polytope $\P_{\mu} = \P_{(\mu, \mu!^*)}$ in $H$ defined by
 \[
\P_{\mu}\coloneqq \left\{q\in H\,\st\,\mu(I)\leq q(I)\leq \mu!^*(I)\text{ for each } I\subseteq V\right\}.
\]
(If $q(V)=\mu(V)$, the inequalities $\mu(I)\leq q(I)$, $I\subseteq V$, are equivalent to $q(I)\leq \mu!^*(I)$, $I\subseteq V$.)

The polytope $\P_{\mu}$ is the base polytope of the polymatroid on the ground set $V$ defined by the submodular function $\mu!^*$. Denoting by $H_d$ the set of points $q\in H$ with $q(V) = d$, it lives in $H_{\mu(V)}$.

In this paper, we introduce and study two families of modular pairs arising from tropical degeneration of linear series on  one-parameter families of curves.

\subsection{Divisor theory on metric graphs and admissible divisors}  Let $G=(V,E)$ be a finite graph and $\ell \colon E \to (0, +\infty)$ be an edge length function associating to each edge $e$ a positive real $\ell_e$. We denote by $\Gamma$ the corresponding metric graph, obtained by plugging an interval $[0,\ell_e]$ of length $\ell_e$ between the two extremities of each edge $e$. We say $(G,\ell)$ is a model for $\Gamma$.

 A divisor $D$ on $\Gamma$ is an element of the free Abelian group generated by the points on $\Gamma$. Writing $(x)$ for the generator associated to each $x\in \Gamma$, we can write $D = \sum_{x\in \Gamma} D(x)(x)$ with integer numbers $D(x)$. The support of $D$, the set of $x\in \Gamma$ for which $D(x)\neq 0$, is finite. The degree of $D$ is $\sum D(x)$.

We say that a divisor $D$ is \emph{$G$-admissible} if each edge $e$ of $\Gamma$ contains at most one point $x$ of the interior $\ring e$ of $e$ in its support, and then, only if $D(x)=1$. An example of a $G$-admissible divisor is depicted in Figure~\ref{fig:admissible}.

To each $G$-admissible divisor $D$ on $\Gamma$, we associate the pair $(\mu!_D, \mu!_D^*)$ of functions on $\ssub{2}!^V$ defined as follows. For a subset $I \subseteq V$, let $\ssGamma_I$ be the subset of $\Gamma$ consisting of all points $x\in\Gamma$ that are either in $I$ or on an edge that connects a pair of vertices of $I$. We set 
\[\mu!_D(I) \coloneqq \sum_{x\in\ssGamma_I} D(x) - \sse_D(I, \ssI^c)+ \frac {\sse(I, \ssI^c)}2
  \qquad \textrm{and} \qquad
  \mu!_D^*(I) \coloneqq \sum_{x\in\ssGamma_I} D(v) + \frac {\sse(I, \ssI^c)}2\]
where $\sse(I, \ssI^c)$ is the number of edges in $G$ that join a vertex in $I$ to one in $\ssI^c$, and $\sse_D(I,\ssI^c)$ is the number of those among them that do not have any point of the support of $D$ in their interior.

We show in Section~\ref{sec:adjoint-modular-admissible-divisor} that $(\mu!_D, \mu!_D^*)$ is a modular pair. It thus gives rise to the base polytope $\P_{\mu!_D}$. Moreover, we provide an alternate description of this polytope as follows. 

For each element $q\in H_d$, a $G$-admissible divisor $D$ is called \emph{$q$-semistable} if 
\begin{equation}
\Big|\sum_{x\in\ssGamma_I}D(x)-q(I)\Big| \leq \frac{\sse(I, \ssI^c)}{2} \quad \textrm{for each subset } I \subseteq V.
\end{equation}
The \emph{semistability polytope} of $D$ denoted $\P_D$ is the polytope in $H$ consisting of all the points $q$ such that $D$ is $q$-semistable. It lives in $H_d$ with $d$ the degree of $D$.

We prove in Section~\ref{sec:semistability-polytope} that for a $G$-admissible divisor $D$ on $\Gamma$, the base polytope $\P_{\mu!_D}$ coincides with the semistability polytope $\P_D$. 

Denoting by $\ssG_D$ the spanning subgraph of $G$ whose edges are those $e\in E$ such that the interior $\ring e$ contains no point in the support of $D$, we show  furthermore in Section~\ref{Voronoi} that, if this subgraph is connected, then $\P_D$ is congruent to the Voronoi cell in $H_0$ associated to $\ssG_D$.

\subsection{Tilings induced by admissible divisors}  Recall that a rational function on $\Gamma$ is a function $f \colon \Gamma \to \R$ that is continuous and whose restriction on each interval $[0,\ell_e]$ in $\Gamma$, $e\in E$, is piecewise affine with integral slopes. The divisor of $f$, denoted $\div(f)$, is defined by
\[\div(f) \coloneqq \sum_{x\in \Gamma} \ord_x(f)(x),\]
where the order of vanishing of $f$ at a point $x$, denoted $\ord_x(f)$, is the sum of the incoming slopes of $f$ at $x$. Its degree is zero. Two divisors $\ssD_1$ and $\ssD_2$ on $\Gamma$ are called linearly equivalent if their difference is the divisor of a rational function. 

Let $\cl$ be a linear equivalence class of divisors on the metric graph $\Gamma$, and denote by  $\stable!_G(\cl)$ the set of all $G$-admissible divisors in $\cl$. More generally, consider an additive subgroup $\valgroup\subseteq\R$ containing the edge lengths $\ell_e$, $e\in E$. The smallest is the subgroup generated by the $\ell_e$; the largest is $\R$ itself. A divisor $D$ in $\cl$ is called \emph{$\valgroup$-rational} if the support of $D$ consists only of $\valgroup$-rational points, that is, points whose distances to vertices of $G$ are all in $\valgroup$. We denote by $\cl_{\valgroup}$ the set of $\valgroup$-rational divisors in the linear equivalence class $\cl$, and denote by $\stable!_G(\cl_{\valgroup})$ the set of those which are also $G$-admissible. (In the context of degenerations of curves, $\valgroup$ is the value group of the underlying valuation.) 

Our first main result, proved in Section~\ref{sec:tiling-semistability}, is the following.
\newtheorem*{thm:semistability-tiling}{\cref{thm:semistability-tiling}}
\begin{thm:semistability-tiling}[Tilings by semistability polytopes of admissible divisors]  Let $\cl$ be a linear equivalence class of divisors of degree $d$ on a metric graph $\Gamma$ with model $(G,\ell)$ and $G$ connected. For a subgroup $\valgroup\subseteq\R$, if $\cl_{\valgroup}$ is nonempty and $\valgroup$ is dense in $\R$, then the set of semistability polytopes $\P_D$ for $D\in\stable!_G(\cl_{\valgroup})$ gives a tiling of $H_d$, that does not depend on the choice of $\valgroup$.
\end{thm:semistability-tiling}

\begin{figure}[!t]
\centering
    \scalebox{.3}{\input{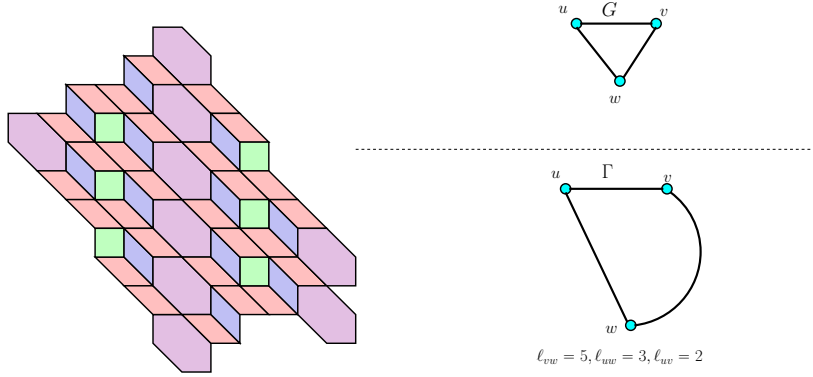}}
\caption{Tiling arising from admissible divisors in the zero class on a metric graph $\Gamma$. The graph $G$ is a $3$-cycle and the edge lengths are depicted in the picture, on the right. The picture is actually that of the projection of the tiling in $H_0\subset\R^V$ to $\R^{\{u,v\}}$.}
\label{fig:tiling3}
\end{figure}

Note that for $\valgroup =\R$, the conditions that $\cl_{\valgroup} \neq \emptyset$ and $\valgroup$ is dense in $\R$ are automatically verified. In principle, one might expect to get more tiles for $\valgroup=\R$ than for a smaller choice. However, since the smaller choice of $\Lambda$ tiles $H_d$ as well, the collection of tiles will be the same. An example of such tiling is depicted in Figure~\ref{fig:tiling3}. Others are depicted in Figures~\ref{fig:tiling2}~and~\ref{fig:tiling4}.

\subsection{Tropicalization of linear series} \label{sec:tropalization-intro} Let $\K$ be an algebraically closed complete non-Archimedean field with a nontrivial valuation. Denote  by $\varR$ its valuation ring, by $\varm$ the maximal ideal of $\varR$, by $\k$ its residue field, and by $\valgroup$ its value group. 

Let $\varX$ be a smooth connected proper curve over $\K$, and denote by $\varX^{\an}$ the Berkovich analytification of $\varX$. Recall that the points on $\varX^{\an}$ are in bijection with the closed points on $\varX$ and valuations on the function field $\K(\varX)$ that extend the valuation of $\K$.
In particular, each rational function $\varf$ on $\varX$ gives rise to an evaluation map $\ev_{\varf} \colon \varX^{\an} \to \R\cup\{\pm \infty\}$.

Let $\Gamma$ be a skeleton for $\varX^{\an}$. It is canonically a metric graph.  Let $(G,\ell)$ be an underlying model, $V$ its vertex set, which we suppose to consist of type 2 points of $\varX^{\an}$, and $E$ its edge set. Denote by $\tau\colon \varX^{\an} \to \Gamma$ the retraction map from $\varX^{\an}$ to $\Gamma$. When $\varX$ is defined over a discretely valued field $\K_0$, by the stable reduction theorem, there exist semistable models for $\varX$ over finite extensions $\K_1$ of $\K_0$.  In this case, $G$ is the dual graph of a semistable reduction of $\varX$, with vertices in bijection with the irreducible components of the semistable curve, and edges in bijection with the nodes. Also, the length of an edge is the singularity degree of the corresponding node, normalized by the degree of the field extension $\rquot{\K_1}{\K_0}$. The metric graph $\Gamma$ is the one associated to $(G, \ell)$; it is independent of the choice of the extension $\rquot{\K_1}{\K_0}$. 

Let $\varD$ be divisor of degree $d$ on $\varX$. Let $\varH$ be a vector subspace of dimension $r+1$ of the space of global sections of the line bundle $\varL = \mathcal O(\varD)$, which  we view in $\K(\varX)$:
 \[\varH \subseteq \bigl \{\varf \in \K(\varX) \,\st\, \div(\varf) + \varD \geq 0\bigr\}.\]
The pair $(\varL, \varH)$ defines a linear series $\g!_{d}^r$ on $\varX$.

The tropicalization of $(\varD, \varH)$ gives rise to a space of rational functions on $\Gamma$.  Let 
\[\tau(\varD) = \sum_{x\in \varX(\K)} \varD(x)(\tau(x))\]
be the divisor on $\Gamma$, the tropicalization of $\varD$ with respect to the skeleton $\Gamma$. For each nonzero rational function $\varf \in\K(\varX)$, we denote by $\trop(\varf)$ the restriction to $\Gamma$ of the evaluation map $\ev_{\varf} \colon \varX^{\an} \to \R\cup\{\pm \infty\}$. This is a rational function on $\Gamma$. Furthermore, by the Specialization Lemma, see e.g.~\cite[Theorem 4.5]{AB15}, we get 
\[\tau(\varD) + \div(\trop(\varf)) = \tau(\varD + \div(\varf)).\]

 We enrich the tropicalization by a collection of $\k$-vector spaces $\VS_v$ associated to the vertices $v$ of $G$. For each vertex $v$, the space $\VS_v$ consists of reductions of functions $\varf \in \varH$ at $v$~\cite[Section 4.4]{AB15}. It was proved in~\cite[Lemma 4.3]{AB15} that $\VS_v$ has dimension $r+1$ over $\k$.  We associate to the linear series $(\varD, \varH)$ the $\k$-vector space 
	 \[\VS \coloneqq \bigoplus_{v\in V} \VS_v.\]

\subsection{Subspaces of $\VS$ associated to admissible divisors and tilings of simplices} Notation as in the previous section, denote by $\cl$ the linear equivalence class of $\tau(\varD)$ on $\Gamma$. Notice that $\tau(\varD)$ is $\valgroup$-rational, and $\valgroup$ is divisible, whence dense in $\R$.

In Section~\ref{sec:tropicalization}, generalizing \cite[Lemma 4.3]{AB15}, we associate an $(r+1)$-dimensional subspace $\ssW_h$ of $\VS$ to each $G$-admissible divisor $D\in \cl_{\valgroup}$ as follows. We write $D=\div(h)+\tau(\varD)$ for a rational function $h$, and denote by $\varM_h$ the space of all functions $\varf\in \varH$ such that $\trop(\varf)(v) \geq h(v)$ for all $v\in V$, that is, 
\[\varM_h \coloneqq \bigl\{f\in \varH \,\st\, \trop(\varf)(v)  \geq h(v)\quad \forall v\in V\bigr\}.\]
As we will show, this is a free $\varR$-module of rank $r+1$. There is a natural $\varR$-linear map $\varM_h\to\VS$ defined by \emph{reduction of functions relative to $h$}, we refer to Section~\ref{sec:tropicalization} for the definition. We let $\ssW_h$ be its image. We will prove that the kernel of this map is $\varm\varM_h$ so that we get an isomorphism $\rquot{\varM_h}{\varm\varM_h}\to\ssW_h$. The proof of this result is based on an important domination property for $G$-admissible divisors, proved in Lemma~\ref{lem:varH}, which stipulates that the inequalities $\trop(\varf)(v)  \geq h(v)$ for vertices $v$ imply the inequalities $\trop(\varf)(x) \geq h(x)$ for all points $x$ on $\Gamma$.

 Changing $h$ by a constant results in the same subspace of $\VS$, so $\ssW_h$ depends only on $D$ and we can thus denote it by $\ssW_D$.  To each $G$-admissible divisor $D$ in $\cl_{\valgroup}$, we associate the  function $\nu!_D^* \colon \ssub{2}!^V \to \R$ defined on each subset $I \subseteq V$ by $\nu!_D^*(I)\coloneqq \dim_{\k} (\ssW_{D,I})$ where $\ssW_{D,I}$ is the projection of $\ssW_{D}$ to $\oplus_{v \in I} \VS_v$. We denote by $\nu!_D$ its adjoint. We show that $(\nu!_D, \nu!_D^*)$ modular pair, that is, $\nu!_D$ is supermodular and $\nu!_D^*$ is submodular. Denote by $\Q_D$ the base polytope associated to $\nu!_D$. Notice that $\Q_D$ lives in the standard simplex $\Delta_{r+1} \subset \R^V$ consisting of all the points $q\in\R_{\geq0}^V$ with $q(V)=r+1$. 
 
 \smallskip

 The following theorem is proved in Section~\ref{sec:tilings-simplex-mixed}.

\newtheorem*{thm:reduction-admissible-tiling}{\cref{thm:reduction-admissible-tiling}}
\begin{thm:reduction-admissible-tiling} 
The collection of polytopes $\Q_D$ associated to $\Lambda$-rational $G$-admissible divisors $D$ in the linear equivalence class of the divisor $\tau(\varD)$ gives a tiling of the standard simplex $\Delta_{r+1}$.
\end{thm:reduction-admissible-tiling}

We refer to Figure~\ref{fig:tiling-simplex} for an example of such a tiling. In Section~\ref{sec:reduced-divisors}, we show that for each vertex $v\in V$, the unique full-dimensional polytope in the tiling that contains the vertex $q_v\in \Delta_{r+1}$ with coordinate $q_v(v)=r+1$ (whence $q_v(u)=0$ for $u\in V$ with $u\neq v$) is the polytope $\Q_D$ associated to the $v$-reduced divisor $D$ in the tropicalization of $\varH$. We refer to Theorem~\ref{thm:reduced-divisors} for the precise statement. In this way, we can view other tiles, and their corresponding $G$-admissible divisors, interpolating between the reduced divisors  with respect to the vertices of $G$ in the linear equivalence of $D$ given by tropicalization of $\varH$.

\subsection{Two-parameter families of tilings associated to tropicalization of linear series} \label{sec:tilings-mixed-intro}
From the proofs of Theorems~\ref{thm:semistability-tiling}~and~\ref{thm:reduction-admissible-tiling}, we obtain the following remarkable result, proved in Section~\ref{sec:tilings-simplex-mixed}.

\newtheorem*{thm:mixed-tiling}{\cref{thm:mixed-tiling}}
\begin{thm:mixed-tiling} 
Let $\alpha,\beta$ be two positive real numbers and $n\coloneqq\alpha d+ \beta(r+1)$. For each $G$-admissible divisor $D$ in $\cl_{\valgroup}$, let $\P_{\alpha, \beta, D}$ be the polytope associated to the supermodular function $\alpha\mu!_D+\beta\nu!_D$. Then, $\P_{\alpha,\beta,D} = \alpha \P_D+\beta \Q_D$, and the collection of polytopes $\P_{\alpha, \beta, D}$ gives a tiling of $H_{n}$. 
\end{thm:mixed-tiling}

\begin{figure}[!t]
\centering
    \scalebox{.19}{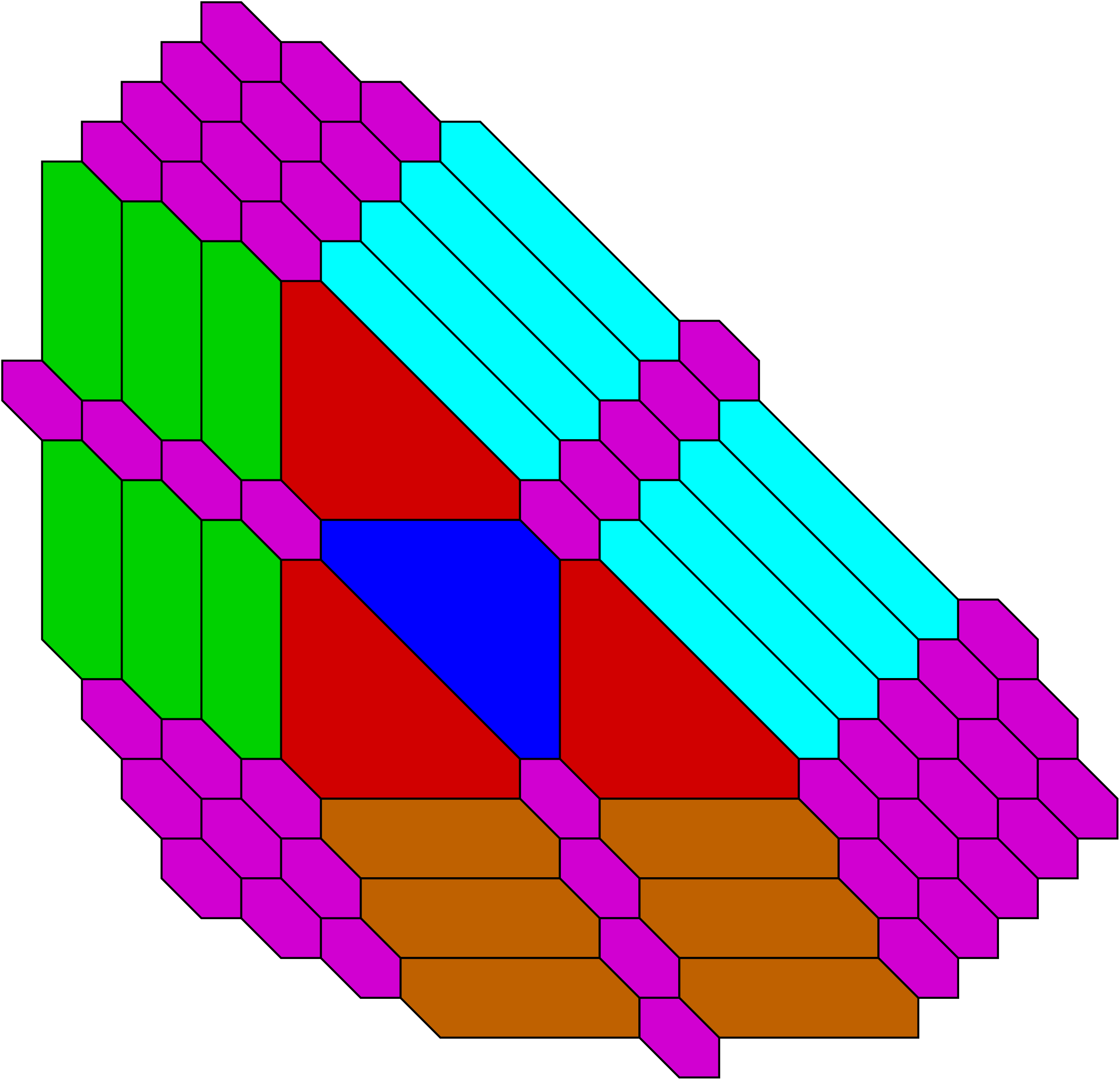}
\caption{An example of tiling $\P_{\alpha, \beta, D}$. The graph $G$ is the $3$-cycle dual to a curve with three components, as in Figure~\ref{fig:tiling3}. The edge lengths are all equal to one. Here $d=6$ and $r=1$. And $\alpha=1$ and $\beta=4$. The tiles $\P_D$ are hexagons, similar to the one depicted in Figure~\ref{fig:tiling4}. The tiling $\Q_D$ of the simplex $\Delta_2$ is the one depicted in Figure~\ref{fig:tiling-simplex}. 
As in Figure~\ref{fig:tiling3}, we depict a projection of the tiling to $\R^2$.}
\label{fig:mixed-tiling}
\end{figure}

 The tilings produced by this theorem play a key role in our approach to address the problem of constructing a moduli space of limit linear series. A more thorough discussion is beyond the scope of this paper and is left to our forthcoming work. An example of such a tiling is given in Figure~\ref{fig:mixed-tiling}.

\subsection{Fundamental collection and finite generation property} 
The \emph{fundamental collection} $\FC_{(\varD, \varH)}$ associated to the pair $(\varD, \varH)$ is the collection of spaces $\ssW_D$ with  full-dimensional $\Q_D$:  
\[
\FC_{(\varD, \varH)}\coloneqq \left\{\ssW_D \st  D \in \stable_G(\cl_{\valgroup}) \, \textrm{  with }\Q_D \textrm{ full-dimensional}\right\}.
\]

Using the domination property, alluded to earlier, we show that for each pair of distinct $G$-admissible divisors $\ssD$ and $\ssD'$ in the linear equivalence class $\cl$ of $\tau(\varD)$, there are certain maps $\ssphi^{D}_{D'}\colon \ssW_{D} \to \ssW_{D'}$ and $\ssphi^{D'}_{D}\colon \ssW_{D'} \to \ssW_D$ such that $\ssphi^{D'}_{D} \circ \ssphi^{D}_{D'} = 0$ and $\ssphi^{D}_{D'} \circ \ssphi^{D'}_D = 0$. Moreover, we show that there is an abundance of pairs that are \emph{exact} in the sense that 
the sequences 
\[\ssW_{D}\to \ssW_{D'} \to \ssW_{D} \quad \textrm{and} \quad \ssW_{D'}\to \ssW_{D} \to \ssW_{D'}\]
are both exact.  We give a geometric meaning to the exactness by translating it  into a measure of proximity of the two tiles associated to $D$ and $D'$. This enables us to move from one tile to the other, and establish the geometric properties needed to establish the results stated previously, see Section~\ref{sec:tilings-families-intro}.
 
Using the above maps,  we prove in Section~\ref{sec:finiteness-lls} the following  finite generation property, stipulating that the fundamental collection recovers all the other reduction spaces. 

\newtheorem*{thm:finiteness-lls}{\cref{thm:finiteness-lls}} 
\begin{thm:finiteness-lls}[Finite generation property] Each $G$-admissible divisor $D \in \stable_G(\cl_{\valgroup})$ whose polytope $\Q_D$ is full-dimensional is necessarily effective.  Furthermore, the fundamental collection $\FC_{(\varD, \varH)}$ of spaces $\ssW_D$ with $\Q_D$ full-dimensional generates all the other spaces in the following sense: For each $G$-admissible divisor $D'\in \cl_{\valgroup}$, the space $\ssW_{D'}$ is generated by the images of the maps $\ssW_D \to \ssW_{D'}$ for $D\in \mathcal S$.
\end{thm:finiteness-lls}

 In our forthcoming work~\cite{AEG23}, we give a geometric description of the fundamental collections in the case of canonical linear series, assuming the limit nodal curve is general in its Deligne--Mumford stratum and approached under a general direction. Using this description, we construct a parameter space for canonical limit linear series. 
 
\subsection{Regularity of the tilings and Chow quotients} Recall from~\cite{GKZ} that a tiling by polytopes of a convex domain $\Omega$ in a real vector space $H$ is called regular if there is a convex function on $\Omega$ which has domains of linearity given by the tiles in the tiling. In Section~\ref{sec:regularity} we prove that the tilings that arise from degeneration of linear series are all regular.

\newtheorem*{thm:regularity-tiling}{\cref{thm:regularity-tiling}}
\begin{thm:regularity-tiling}  The tilings in Theorem~\ref{thm:semistability-tiling}, in Theorem~\ref{thm:reduction-admissible-tiling}, and in Theorem~\ref{thm:mixed-tiling} are all regular. 
\end{thm:regularity-tiling}

In equal characteristic, when the valued field $\K$ and its residue field $\k$ have the same characteristic, we can assume that the valuation ring $\varR$ contains a copy of the residue field $\k$. In this case, we establish in Section~\ref{sec:chow-quotients} a connection between the tilings produced by Theorem~\ref{thm:reduction-admissible-tiling} and the works~\cite{Kap-chow, Laf99, GW22}. We briefly discuss this here.

Notation as in Section~\ref{sec:tropalization-intro}, recall that $\varH\subset \K(\varX)$ is a $\K$-vector space of dimension $r+1$ of the space of global sections of the line bundle $\varL=\mathcal O(\varD)$.  Consider the Grassmannian $\varGr$ of $(r+1)$-planes in the vector space $\left(\k^{r+1}\right)^{V}$.  The torus 
$\varG^V_{\varm}$ acts naturally on the Grassmannian, each factor acting on the corresponding copy of $\k^{r+1}$. 

We show in Section~\ref{sec:chow-quotients} that each rational function $h$ on $\Gamma$ such that $\ssD=\div(h)+\tau(\varD)$ is $G$-admissible gives rise to a morphism (after appropriate bases of $\varH$ and the $\varM_h$ are chosen) 
\[
\beta_h\colon\text{Spec}(\varR)\longrightarrow   \varGr.
\]
The morphism sends the generic point of $\textrm{Spec}(\varR)$ to a point on the orbit of $\delta(\varH) \subset \varH^{V}$, where $\delta\:\varH \hookrightarrow \varH^V$ is the diagonal embedding, and the special point to a point on the orbit of $\ssW_{\ssD}$.

Since the action of the diagonal $\varG_{\varm}\to\varG^V_{\varm}$ on $\varGr$ is trivial, we can consider the action by the quotient $\varT\coloneqq\varG^V_{\varm}/\varG_{\varm}$.  Let $\varY$ be the open subscheme of the Grassmannian where the action of $\varT$ is free. There is an inclusion $\rquot{\varY}{\varT}\hookrightarrow \mathrm{Chow}(\varGr)$ of the quotient $\rquot{\varY}{\varT}$ to the Chow variety of $\varGr$, taking an orbit to its closure, viewed as a prime cycle of 
$\varGr$. The \emph{Chow quotient} is the closure of $\rquot{\varY}{\varT}$ in $\mathrm{Chow}(\varGr)$. The diagonal $\delta(\varH)$ gives a $\K$-point on $\varGr$ whose orbit closure corresponds to a subscheme of $\varGr\times_{\k}\K$. As shown in \cite{GW22},  its closure $\mathcal Z\subseteq \varGr\times_{\k}\mathrm{Spec}(\varR)$ has closed fiber $\varZ$ over $\k$ with reduced associated cycle of the form $\sum_i \varZ_i$, where each prime cycle $\varZ_i$ is the closure of the $\varT$-orbit of a subspace $L_i\subseteq(\k^{r+1})^V$. Furthermore, the associated polytope $\Q_{L_i}$ to $L_i$, as defined in Section~\ref{sec:modular-subspace}, is full-dimensional, and their collection forms a polyhedral decomposition of $\Delta_{r+1}$. 

Let $\cl$ be the linear equivalence class of $\tau(\varD)$ and let $\cl_{\valgroup}\subseteq\cl$ be the subset of $\Lambda$-rational divisors. We show the following theorem. 

\newtheorem*{thm:comparison-chow-quotient}{\cref{thm:comparison-chow-quotient}}
\begin{thm:comparison-chow-quotient}
  Notation as above, the $\ssW_{D}$ for $D\in \stable!_G(\cl_{\valgroup})$ are in the orbit closures of the $L_i$. Conversely, each $L_i$ is in the orbit of $\ssW_D$ for some $D\in \stable!_G(\cl_{\valgroup})$. In particular, the tiling by the $\Q_D$ in Theorem~\ref{thm:reduction-admissible-tiling} coincides with the tiling by the $\Q_{L_i}$ and their faces.
\end{thm:comparison-chow-quotient}  

\begin{center}
***
\end{center}

 \smallskip

To obtain the above results, we rely on general theorems proved in Section~\ref{sec:admissible-families-tilings} that ensure certain collections of modular pairs give rise to tilings, and use interesting features of admissible divisors to ensure the properties required in these theorems are verified. We discuss briefly these results here.

\subsection{Tilings induced by families of modular pairs}  \label{sec:tilings-families-intro}

Let $n\in\Z$. Let $V$ be a finite nonempty set and denote by $\M_n(V)$ the set of supermodular functions $\mu$ on $\ssub{2}!^V$ with $\mu(V)=n$. Let $\C \subset \M_n(V)$ be a subfamily. We denote by $\P_{\C}$ the union of the polytopes $\P_{\mu}$ for $\mu \in \C$ and call it the \emph{support} of the family $\C$.

An ordered bipartition of $V$ is a pair $(I,J)$ consisting of disjoint nonempty subsets $I$ and $J$ of $V$ with $V = I \sqcup J$. An ordered bipartition $\pi$ induces a \emph{splitting} $\mu!_{\pi}$ of each supermodular function $\mu$, defined by 
\[\mu!_{\pi}(S) \coloneqq \mu(S\cup I) - \mu(I) + \mu(S\cap J)\quad\text{for each }S\subseteq V.\]
The splitting $\mu!_\pi$ is again supermodular and verifies $\mu!_{\pi}(V) =\mu(V)$.

 We say $\C$ is \emph{separated} if each two distinct elements $\mu, \nu \in \C$ admit a \emph{nontrivial separation}. This means there exists an ordered bipartition $\pi=(I,J)$ such that $\mu(I) +\nu(J) \geq n$, and either this inequality is strict, or we have $\mu!_{\pi} \neq \mu$ or $\nu!_{\pi}\neq \nu$.

We say $\C$ is \emph{closed} if for each $\mu\in \C$ and each bipartition $\pi$ of $V$, we have $\mu!_{\pi}\in \C$.

\newtheorem*{thm:tiling}{\cref{thm:tiling}}
\begin{thm:tiling}[Tilings induced by separated and closed families] If $\C \subseteq \M_n(V)$ is separated and closed, then the collection of polytopes $\P_{\mu}$ for $\mu \in \C$ gives a tiling of $\P_{\C}$.
\end{thm:tiling}

In order to characterize the support of a family, we introduce the following two properties.

 Let $\Omega\subseteq H_n$ be a subset. We say $\C$ is \emph{complete for $\Omega$} if for each $\mu\in\C$ and each bipartition $\pi=(I,J)$ of $V$ such that $\P_{\mu_\pi}$ intersects $\Omega$ and $\mu!_{\pi}=\mu$, there is $\mu'\in\C$ distinct from $\mu!_{\pi}$ such that $\mu!_{\sspi^c}'=\mu!_{\pi}$, where $\sspi^c\coloneqq (J,I)$. 

Given an element $\mu$ of $\M_n(V)$, we call the minimum value among all differences $\mu!^*(I) - \mu(I)$ for nonempty proper subsets $I\subset V$ the \emph{spread} of $\mu$. We say $\C$ is \emph{positive} if the infimum of the spreads which are positive of all elements in $\C$ is positive. In particular, there must exist elements in $\C$ with positive spreads.

\newtheorem*{thm:tiling2}{\cref{thm:tiling2}}
\begin{thm:tiling2}[Tilings induced by complete and positive families] Let $\Omega$ be a connected open subset of $H_n$. Assume $\P_{\C}$ intersects $\Omega$. If $\C$ is complete for $\Omega$ and positive, then $\P_{\C}\supseteq\Omega$. 
\end{thm:tiling2}

\subsection{Characterization of admissible divisor classes} Making use of Theorems \ref{thm:tiling} and \ref{thm:tiling2}, the proofs of Theorems~\ref{thm:semistability-tiling},~\ref{thm:reduction-admissible-tiling}~and~\ref{thm:mixed-tiling} reduce in each case to proving that the corresponding family of supermodular functions is separated, closed, positive and complete. In establishing these properties, we rely on an interesting feature of admissible divisors on a metric graph, proved in Section~\ref{sec:admissible-characterization} and of independent interest, given by the following description of $G$-admissible divisors in linear equivalence classes. 

\newtheorem*{thm:existence-uniqueness-canonical-extension}{\cref{thm:existence-uniqueness-canonical-extension}}
\begin{thm:existence-uniqueness-canonical-extension}
Let $D$ be a divisor on a metric graph $\Gamma$. For each function $f \colon V \to \R$, there exists a unique rational extension $\ssfhat \colon \Gamma \to \R$
such that $D+ \div(\ssfhat)$ is $G$-admissible. 
\end{thm:existence-uniqueness-canonical-extension}

The theorem thus provides a bijection between real-valued functions on vertices, up to addition by constants, and $G$-admissible divisors in a given linear equivalence class.

\subsection{Related work}\label{sec:related-work}

\subsubsection{Degeneration of linear series} The literature on degenerations of line bundles in families of curves is rich. We refer to \cite{DS, OS79, AK1, AK2, Cap, Esteves01, MRV, KP}, to name a few references. That of degenerations of linear series has fewer references. It was systematized in the pioneering work by Eisenbud and Harris \cite{EH86}, but only for limit curves of compact type. As it led to many results, they posed the problem of extending the study to more general stable curves in \cite{EHBull}. The problem remains unresolved to this day.

Almost twenty years later, Osserman advocated looking at all limits of linear series in \cite{Os-moduli} for degenerations to two-component curves of compact type. A few years later, Esteves and Osserman \cite{EO} made it precise, by explaining the need to consider all limits of linear series to account for the limits of their associated divisors. For a more general nodal curve, Osserman described later in 
\cite{Os-pseudocompact} a certain quiver representation that
captures all the limits of linear series, but did not dwell much on it, preferring to consider only the ``extreme" limits of linear series (in our terminology, they correspond to those $\ssW_D$  associated to the $v$-reduced divisors in the linear equivalence class), more in line with the approach by Eisenbud and Harris \cite{EH86}. Those had been considered by Amini and Baker in \cite{AB15}. Osserman explains the connection in \cite{Os-AB}.

 Recently, Esteves, Santos and Vital \cite{ESV1, ESV2} studied the quiver representation described by Osserman in \cite{Os-pseudocompact}, though only for degenerations with regular total space, which correspond to metric graphs $\Gamma$ whose edge lengths are all equal to 1. They found many properties of the representations, illustrating them with a few examples. Our examples, depicted in Figures~\ref{fig:mixed-tiling}~and~\ref{fig:tiling-simplex}, borrow the example of degeneration of linear series worked out in \cite[Ex.~10.2]{ESV1}, and give for it the tiling $\Q_D$ corresponding to the vector spaces and the mixed tiling $\P_{\alpha,\beta,D}$ corresponding to the admissible divisors and the vector spaces.

\subsubsection{Tropical point of view} Divisor theory for graphs was developed in the work by Baker and Norine~\cite{BN07}, extended to metric graphs by Mikhalkin and Zharkov \cite{MZ08} and Gathman and Kerber \cite{GK08}, and to the hybrid setting of metrized complexes by Amini and Baker~\cite{AB15}. Tropical methods have been used in the study of the geometry of curves and their moduli spaces, with notable application in the recent work by Farkas, Jensen and Payne~\cite{FJP} to prove that the moduli spaces of curves of genera 22 and 23 are of general type. We refer to Baker and Jensen~\cite{BJ} and Jensen and Payne~\cite{JP21} for a survey of these results.

Semistability conditions for divisors on metric graphs have appeared before from a different perspective. We direct the reader to the recent works \cite{AP20}, \cite{AAPT} and \cite{MMUV}, \cite{CPS23}, and the references therein. In these references, the focus is on the construction of tropical counterparts to Caporaso's compactification of the universal Picard variety over the moduli of stable curves \cite{Cap}. There, certain polyhedral decompositions appear from the stratification of the tropical Jacobians by the combinatorial types of the objects considered. Compared to the polytopes studied in our paper, those appearing in the above context are all parallelepiped, they lie in a different space, and capture aspects related to the asymptotic geometry of Picard varieties.

\subsubsection{Matroidal subdivisions}
Matroidal subdivisions appear in the work by Kapranov on Chow quotients~\cite{Kap-chow} and its sequel by Keel and Tevelev~\cite{KT06}. Polymatroidal subdivisions of the simplex have been studied in the work by Lafforgue~\cite{Laf99}. A recent work by Giansiracusa and Wu~\cite{GW22} provides a link between the work by Lafforgue, Kapranov and Thaddeus~\cite{Tad99}. 
We will make the connection between our framework and these works precise in Section~\ref{sec:chow-quotients}.

Other results on the connection between matroidal subdivisions and the geometry of Grassmannians can be found in the works by Gelfand, Goresky, MacPherson and Serganova \cite{GGMS87}, Lafforgue~\cite{Laf03}, and Keel and Tevelev \cite{KT06}.  The work by Ardilla and Billey~\cite{AB07} relates matroidal subdivisions of a product of two simplices to the geometry of flag varieties.  A discussion of features related to the geometry of positive Grassmannians can be found in the survey papers by Postnikov~\cite{Pos18} and Williams~\cite{Wil22}. 
Matroidal tilings are also related to rigidity properties for matroids themselves; see the work by Lafforgue \cite{Laf03} and a recent work by Baker and Lorscheid~\cite{BL23}. A nice introduction to the topic can be found in the lecture notes by Fink \cite{Fink15}.

\subsection{Organization of the paper} In the preliminary Section~\ref{sec:modular-pairs}, we recall basic results about polymatroids and discuss their face structure. In Section~\ref{sec:admissible-families-tilings}, we prove the tiling Theorems~\ref{thm:tiling}~and~\ref{thm:tiling2} associated to families of modular pairs, see Section~\ref{sec:tilings-families-intro}.

In Section~\ref{sec:admissible-divisor-semistability-polytope}, we study semistability polytopes of admissible divisors and show they are examples of polymatroids. 
Section~\ref{sec:admissible-characterization} describes the admissible divisors in linear equivalence classes, Theorem~\ref{thm:existence-uniqueness-canonical-extension}. 
Section~\ref{sec:tiling-semistability} is devoted to the proof of Theorem~\ref{thm:semistability-tiling} on tilings induced by semistability polytopes of $G$-admissible divisors. 

In Section~\ref{sec:tropicalization}, we establish fundamental results about reductions of spaces of sections of line bundles whose associated divisors are admissible. In Section~\ref{sec:tilings-simplex-mixed}, we prove Theorems~\ref{thm:reduction-admissible-tiling} and~\ref{thm:mixed-tiling} on tilings induced by reduction of linear series, and tilings obtained by mixing semistability polytopes and reduction of linear series, see Section~\ref{sec:tilings-mixed-intro}. 

Section~\ref{sec:geometric-features} is devoted to discussing geometric features of the tilings introduced in this paper. We prove Theorem~\ref{thm:regularity-tiling} on regularity of these tilings, Theorem~\ref{thm:comparison-chow-quotient} on the connection between Theorem~\ref{thm:reduction-admissible-tiling} and Chow quotients, Theorem~\ref{thm:reduced-divisors}, which gives a description of the tiles containing the vertices of the simplex in terms of reduced divisors in tropical linear series, and the finite generation property stated in Theorem~\ref{thm:finiteness-lls}. 

Finally, Section~\ref{sec:discussion} discusses complementary related results, in particular, the question of periodicity of the tilings issued from Theorem~\ref{thm:semistability-tiling}.

\smallskip

 Sections~\ref{sec:modular-pairs} through~\ref{sec:tiling-semistability} are written to be  accessible both to a combinatorial and an algebro-geometric audience. Sections~\ref{sec:tropicalization} and~\ref{sec:tilings-simplex-mixed} require some knowledge of non-Archimedean geometry and the theory of Berkovich curves, that we briefly recall in  Section~\ref{sec:tropicalization}. 

\subsection*{Acknowledgment} We are grateful to the hospitality of CMLS at École Polytechnique and IMPA. We thank Carolina Araujo, Dan Corey, Federico Ardila, Felipe Léon, Lucas Gierczak, and Oliver Lorscheid for discussions on related topics.

Part of this research was conducted while the first author was visiting scholar at the Berlin Mathematics Research Center MATH+. He thanks MATH+ for support, and the mathematics institutes at TU and HU for warm hospitality. 

The second author was supported by CNPq, Proc.~307675/2019, CAPES-COFECUB, Project 932/19, and FAPERJ, Proc.~E-26/202.992/2017.


\section{Adjoint modular pairs and their base polytopes}\label{sec:modular-pairs}

In this section, we define modular pairs and their base polytopes. The material in this section is essentially known and we refer to~\cite{Sch} for more details.  

Let $V$ be a finite nonempty set and denote by $H$ the vector space $\R^V$. Elements of $H$ are maps $q \colon V \to \R$. For each integer $n$, we denote by $H_n$ the affine subspace of $H$ consisting of all $q$ with $\sum_{v\in V} q(v)=n$.

We denote by $\ssub{2}!^V$ the family of subsets of $V$. For a subset $I\subset V$, we denote by $\ssI^c$ the complement of $I$ in $V$, that is, $\ssI^c =V \setminus I$.

\subsection{Adjoint modular pairs} We say a function $\mu\colon \ssub{2}!^V\to\R$
is \emph{supermodular} if $\mu(\emptyset)=0$ and we have the inequalities
\[
\mu(\ssI_1)+\mu(\ssI_2)\leq \mu(\ssI_1\cup \ssI_2)+ \mu(\ssI_1\cap \ssI_2)
\quad\text{for each }\ssI_1,\ssI_2\subseteq V.
\]
 The quantity $\mu(V)$ is called the \emph{range} of $\mu$. Similarly, we say a function  $\eta\colon \ssub{2}!^V \to \R$ is \emph{submodular} if $\eta(\emptyset)=0$ and the inequalities above are all reversed, and call $\eta(V)$ the range of $\eta$. Functions that are both supermodular and submodular are called \emph{modular}.

Denote by $\M(V)$ the set of all supermodular functions on $\ssub{2}!^V$ and by
$\M_n(V)$ its subset consisting of those of range $n$, for
$n\in\mathbf R$. Similarly, we use the notation $\subM(V)$ and $\subM_n(V)$ for the set of submodular functions and its subset consisting of those of range $n$, respectively. 

For each $\mu\colon \ssub{2}!^V\to\R$, let $\mu!^* \colon \ssub{2}!^V\to\R$ be the
\emph{function adjoint to $\mu$} defined by
\[
\mu!^*(I)\coloneqq\mu(V)-\mu(V\setminus I)\quad\text{for each }I\subseteq V.
\]
It is easy to see that $\mu$ is supermodular,
resp.~submodular, if and only if $\mu!^*$ is submodular,
resp.~supermodular. Furthermore, we have the equality $\ssub{(\mu!^*)}!^* = \mu$, and $\mu$ and $\mu!^*$ are of the same range. We will refer to $\mu!^*$ as the \emph{submodular function adjoint to $\mu$} and call the ordered pair $(\mu,\mu^*)$
an \emph{adjoint modular pair} for $V$. Note that in an adjoint modular pair, the first component is always supermodular, while the second is always submodular. 

For $\mu$ supermodular, applying the supermodularity
inequality to $I$ and $\ssI^c=V\setminus I$,
and using that $\mu(\emptyset) =0$, we get that $\mu(I)\leq\mu!^*(I)$ for each $I\subseteq V$. Moreover, the equality holds for
$I=\emptyset$ and $I=V$. 

 Modular functions on $\ssub{2}!^V$ are in bijection with $H=\R^V$: each element $q\in H$ can be viewed as a modular function $q\colon\ssub{2}!^V\to\R$ by putting  $q(I)\coloneqq \sum_{v\in I}q(v)$ for each $I\subseteq V$, with the convention that $q(\emptyset)\coloneqq 0$. Moreover, the pair $(q,q)$ is a modular pair. 

\subsection{Base polytope associated to an adjoint modular pair}

For each adjoint modular pair $(\mu,\mu!^*)$, we define the subset in $H$
\[
\P_{\mu}=\P_{(\mu,\mu!^*)}\coloneqq \left\{q\in H\,\st\,\mu(I)\leq q(I)\leq \mu!^*(I)\text{ for each } I\subseteq V\right\}.
\]

Clearly, $\P_{\mu}$ is a convex polytope if it is nonempty. Also,  $\P_{\mu}\subset H_n$, where $n\coloneqq \mu(V)$, where, we recall, $H_n\coloneqq \big\{q\in\R^V\,\big|\, q(V)=n\big\}$.

If $q\in H_n$, then $q(I)\leq\mu!^*(I)$ if and only if $\mu(\ssI^c)\leq q(\ssI^c)$ for each $I\subseteq V$. This means that one set of inequalities in the definition of $\P_{\mu}$ implies the other. 

We call $\P_{\mu}$ the \emph{base polytope} associated to $\mu$.  The following proposition justifies the terminology. 

\begin{prop} Let $(\mu,\mu!^*)$ be an adjoint modular pair. The polytope $\P_{\mu}$ is the base polytope of the polymatroid defined by $\mu!^*$.
\end{prop}

\begin{proof} The polymatroid defined by $\mu!^*$ is the set of all
  points $q \in H$ that verify $q(I) \leq \mu!^*(I)$ for all $I
  \subseteq V$~\cite[Chapter 44]{Sch}. The base polytope is the
  face of the polymatroid consisting of those $q$ that satisfy $q(V) = \mu!^*(V)$. As we noted in the previous paragraph, the inequalities $\mu(I) \leq q(I)$, $I \subseteq V$, are implied by the inequalities involving $\mu!^*$. This proves the claim. 
\end{proof}

\subsection{Ordered partitions} An \emph{ordered partition} of $V$ is an ordered sequence $\sspi = (\sspi_1, \dots, \sspi_s)$ of
pairwise disjoint subsets of $V$  whose union is equal to $V$. It is called \emph{nontrivial} if $\sspi_i\neq\emptyset$ for every
$i$. If $\pi$ is nontrivial, we call it a \emph{bipartition} if $s=2$, a
\emph{tripartition} if $s=3$, and an \emph{$s$-partition} in
general. If $\pi$ is a bipartition, we put $\sspi^c\coloneqq (\sspi_2,\sspi_1)$.

The data of an ordered partition $\pi$ as above is equivalent to the data of a filtration 
\[\filt_\bullet\coloneqq(\filt_0,\dots,\filt_s), \text{ where } \emptyset=\filt_0 \subseteq \filt_1 \subseteq \filt_2 \subseteq \dots \subseteq \filt_s = V,\]
via the correspondence $\filt_j \coloneqq \sspi_1 \cup \dots \cup
\sspi_j$ for $j=1, \dots, s$, in one direction, and
$\sspi_j \coloneqq \filt_j\setminus \filt_{j-1}$ for $j=1, \dots, s$, in the other direction.

\subsection{Restriction, contraction and splitting}
For each supermodular function $\mu\in\M(V)$ and each
pair of subsets $\ssJ_1,\ssJ_2\subseteq V$ with $\ssJ_1\subseteq \ssJ_2$, we define 
\[
\mu!_{J_2/J_1}(I)\coloneqq\mu\bigl((I\cap \ssJ_2)\cup \ssJ_1\bigr)-\mu(\ssJ_1)\quad\text{for each }
I\subseteq V.
\]
Then, we have $\mu!_{J_2/J_1}\in\M(V)$. We say that $\mu!_{J_2/J_1}$ is obtained from $\mu$ by \emph{restricting to $\ssJ_2$ and contracting $\ssJ_1$.}

Given a subset $J\subseteq V$, we also have the restriction map $\cdot \rest{J}\colon \M(V) \to \M(J)$, that sends $\mu$ to $\mu\rest{J}$ given by $\mu\rest{J}(I)=\mu(I)$ for each $I\subseteq J$. Clearly, 
$\mu!_{J/\emptyset}\rest{J}=\mu\rest{J}$, and more generally, $\mu!_{J_2/J_1}\rest{J_2}=\ssub{\bigl(\mu\rest{J_2}\bigr)}!_{J_2/J_1}$ for a pair $\ssJ_1 \subseteq \ssJ_2$ of subsets of $V$. Moreover, $\mu!_{J_2/J_1}$ is obtained from $\ssub{\bigl(\mu\rest{J_2}\bigr)}!_{J_2/J_1}$ by composing with the map $\ssub{2}!^V \to \ssub{2}!^{J_2}$ that sends $I\subseteq V$ to $I\cap \ssJ_2 \subseteq \ssJ_2$.

For each ordered partition $\pi = (\pi_1, \dots, \pi_s)$ of $V$ with
 $\filt_\bullet$ the corresponding filtration on $V$, we define the
 function $\mu!_\pi$ on $\ssub{2}!^V$ as the sum:
\begin{equation}\label{muV}
\mu!_\pi\coloneqq \sum_{i=1}^s\mu!_{\filtind_i/\filtind_{i-1}}.
\end{equation}
Being a sum of supermodular functions, $\mu!_\pi$ is itself
supermodular and so belongs to $\M(V)$. Moreover, it has the same range as $\mu$. 

We say that $\mu!_\pi$ is the \emph{splitting of $\mu$ with respect to the ordered partition $\pi$}. We call a supermodular function $\nu \in \M(V)$ a \emph{splitting of $\mu$} if it coincides with $\mu!_\pi$ for some ordered partition $\pi$ of $V$, \emph{nontrivial} if $\nu\neq\mu$.  We denote by $\mu!_\pi^*$ the corresponding submodular function, so that $(\mu!_\pi, \mu!_\pi^*)$ is an adjoint modular pair.

The following proposition gives basic properties of splittings. 

\begin{prop}\label{muprod} Let $\mu\in\M(V)$, and $\pi = (\sspi_1, \dots, \sspi_s)$
  be an ordered partition of $V$. Then, for each  $I\subseteq V$, we have:
  \begin{itemize}
      \item $\mu!_\pi(I)=\mu!_\pi(I\cap \sspi_1)+\cdots+\mu!_\pi(I\cap\sspi_s)$,
      \item $\mu(I)\leq\mu!_\pi(I)\leq\mu!^*_\pi(I)\leq\mu!^*(I)$.
  \end{itemize}
Furthermore, the following  statements
  are equivalent:
\begin{enumerate}[label=(\arabic*)]
  \item \label{equi1} $\mu=\mu!_\pi$.
  \item \label{equi2} $\mu(V)=\mu(\sspi_1)+\cdots+\mu(\sspi_s)$.
  \item \label{equi3} $\mu(I)=\mu(I\cap \sspi_1)+\cdots+\mu(I\cap \sspi_s) $
    for each subset $I \subseteq V$.  
  \item \label{equi4} $\mu=\mu!_{\pi'}$ for each reordering $\pi'$ of $\pi$.
  \end{enumerate}
  If they hold, we have 
\[
\P_{\mu}=\P_{\mu\rest{\pi_1}}\times\cdots\times\P_{\mu\rest{\pi_s}}
\]
under the natural identification $\R^V = {\mathbb R}^{\sspi_1} \times \dots \times {\mathbb R}^{\sspi_s}$.
\end{prop}

\begin{proof} Let $\filt_\bullet$ be the filtration of $V$
  corresponding to $\pi$. Consider the first statement. Since
  $\mu!_\pi(I\cap\sspi_i)=\mu!_{\filtind_i/\filtind_{i-1}}(I)$ for each $i$ and
  $I\subseteq V$, the equality follows
  directly from \eqref{muV}.
  Then, 
\begin{align*}
\mu!_\pi(I)=\sum\mu!_{\pi}(I\cap\sspi_i)=&\sum\big(\mu((I\cup\filt_{i-1})\cap\filt_i)-\mu(\filt_{i-1})\big)\\
\geq&
\sum\big(\mu(I\cup\filt_{i-1})+\mu(\filt_i)-\mu(I\cup\filt_i)-\mu(\filt_{i-1})\big)=\mu(I).
\end{align*}
 The
  second inequality $\mu!_\pi(I)\leq\mu!^*_\pi(I)$ is a property of adjoint modular pairs. The third
  follows from the first and adjunction, as, for each $I\subseteq V$,
\[
\mu!_\pi(\ssI^c)+\mu!_\pi^*(I)=\mu!_\pi(V)=\mu(V)=\mu(\ssI^c)+\mu!^*(I).
\]

 We now prove the equivalence of \ref{equi1}-\ref{equi2}-\ref{equi3}-\ref{equi4}. By the first part of the proposition, \ref{equi1} implies \ref{equi2}. In addition, by supermodularity,  for each $I\subseteq V$, we have
  \[
  \sum_{i=1}^s\mu(I\cap \sspi_i)\leq\mu(I)\leq
  \sum_{i=1}^s\mu(I\cap \sspi_i)+\mu(V)-\sum_{i=1}^s\mu(\sspi_i).
  \]
Thus, \ref{equi2} implies \ref{equi3}.

  Assume \ref{equi3}. To prove \ref{equi1}, it is enough to show that $\mu(I)=\mu!_\pi(I)$ for each $I\subseteq \sspi_i$ and each $i=1, \dots, s$. But
  then, 
  \[
  \mu!_\pi(I)=\mu(I\cup \filt_{i-1})-\mu(\filt_{i-1})=\mu(I)
  \]
  from the definition and \ref{equi3}.

  Finally, \ref{equi1} and \ref{equi4} are equivalent because they are both equivalent to \ref{equi2}.
  
  The last statement about base polytopes is clear. 
\end{proof}

\subsection{Refinements}    
Given two ordered partitions $\pi = (\sspi_1, \dots, \sspi_s)$ and $\pi' =
(\sspi_1', \dots, \sspi_l')$ of $V$, we say that $\pi'$ is a
\emph{refinement} of $\pi$ and write $\pi \preceq \pi'$ if there is an
increasing sequence $\ell_0,\dots,\ell_s$ of integers with $\ell_0=0$
and $\ell_s=l$ such that 
\[\sspi_{j} = \sspi_{\ell_{j-1}+1}' \cup \dots \cup \sspi_{\ell_j}',\qquad j=1, \dots, s.\]
Equivalently, $\pi \preceq \pi'$ if there are ordered partitions $(\sspi_{i,1},\dots,\sspi_{i,\ell_i})$ of $\sspi_i$, $i=1,\dots,s$, such that
\[
\pi' = (\sspi_{1,1},\dots,\sspi_{1,\ell_1},\sspi_{2,1},\dots,\sspi_{2,\ell_2},\dots,\sspi_{s,1},\dots,\sspi_{s,\ell_s}).
\]
Denoting by $\filt_\bullet$ and $\filt'_\bullet$ the filtrations corresponding to $\pi$ and $\pi'$, respectively, this is further equivalent to requiring that each $\filt_i$ in $\filt_\bullet$ appears at least as many times in $\filt'_\bullet$ as in $\filt_\bullet$.

At any rate, if $\pi'$ is an arbitrary ordered partition as above, it induces a refinement $\pi''$ of $\pi$, where
\[
\pi''\coloneqq (\sspi_1'\cap\sspi_1,\dots,\sspi_l'\cap\sspi_1,\sspi_1'\cap\sspi_2,\dots,\sspi_l'\cap\sspi_2,\dots,\sspi_1'\cap\sspi_s,\dots,\sspi_l'\cap\sspi_s).
\]
The associated filtration to $\pi''$ is
\[
\filt_1\cap\filt_1'\subseteq\cdots\subseteq\filt_1\cap\filt_l'\subseteq
\filt_1\cup (\filt_2\cap\filt_1')\subseteq\cdots\subseteq
\filt_1\cup (\filt_2\cap\filt_l')\subseteq
\filt_2\cup (\filt_3\cap\filt_1')\subseteq\cdots\subseteq\filt_s.
\]

One special case of a refinement of $\pi$ is the one obtained by adding empty sets in any ordering. In this case, we say as well that $\pi$ is obtained from its refinement by removing (certain) empty sets. If $\pi'$ is a refinement of $\pi$, then $\pi''$ is obtained from $\pi'$ by adding empty sets.

Each ordered partition $\pi$ is a refinement of the 1-partition consisting of a single part $V$. Alternatively, $\pi$ is obtained from the 1-partition by the succession of refinements induced by the bipartitions among  $(\filt_{\ell},\filt_{\ell}^c)$ for $\ell=1,\dots,s-1$, after empty sets are removed. If $\pi'$ is obtained from the 1-partition by a succession of refinements induced by bipartitions, after empty sets are removed, so is $\pi''$ obtained from $\pi$ by the same succession. These facts will allow us to reduce certain arguments below.

If $\pi'$ is obtained from $\pi$ by adding empty sets, then 
$\mu!_{\pi}=\mu!_{\pi'}$. At any rate, the following proposition holds.

\begin{prop}\label{refineV} Let $\mu\in\M(V)$. Let $\pi$ and $\pi'$ be
  ordered partitions of $V$. Then
  $\ssub{(\mu!_\pi)}!_{\pi'}=\mu!_{\pi''}$, where $\pi''$ is the
  refinement of $\pi$ induced by $\pi'$.
\end{prop}

\begin{proof} Proceeding by induction, we only need to prove the statement in
  the case where $\pi'$ is a bipartition of $V$. Thus, $\pi'=(\sspi_1',\sspi_2')$. 
  
  Let $\pi = (\sspi_1, \dots, \sspi_s)$ and consider the filtration $\filt_\bullet$ associated to $\pi$.
  The filtration associated to $\pi''$ has an extra term  $\filt_{\ell}^-\subseteq V$ 
  satisfying $\filt_{\ell-1}\subseteq \filt^- _{\ell}\subseteq
  \filt_{\ell}$ for each $\ell=1,\dots,s$. Put $\varpi!_{\ell,1}\coloneqq\filt^-_{\ell}\setminus \filt_{\ell-1}$ and
  $\varpi!_{\ell,2}\coloneqq \filt_\ell\setminus \filt^-_{\ell}$ so
  that $(\varpi!_{\ell,1}, \varpi!_{\ell,2})$ forms a partition of
  $\sspi_\ell$ for each $\ell$.  Indeed,
  $\varpi!_{\ell,i}\coloneqq\sspi_\ell\cap\sspi_i'$ for each $\ell=1,\dots,s$ and
  $i=1,2$.

  Now, by supermodularity and Proposition~\ref{muprod}, 
  \[
  \mu(V)=\ssub{(\mu!_\pi)}!_{\pi'}(V)\geq\sum_{\ell=1}^s\sum_{i=1,2}\ssub{(\mu!_\pi)}!_{\pi'}(\varpi!_{\ell,i})=\sum_{\ell=1}^s
  \ssub{(\mu!_\pi)}!_{\pi'}(\sspi_\ell)\geq
  \sum_{\ell=1}^s\mu!_\pi(\sspi_\ell)=\mu!_\pi(V)=\mu(V),
  \]
  whence equalities hold. By Proposition~\ref{muprod} again, it is
  enough to show that $\ssub{(\mu!_\pi)}!_{\pi'}(I)=\mu!_{\pi''}(I)$
  for each subset $I\subseteq\varpi!_{\ell,i}$ for each $\ell$ and $i$.

  Now, for each $\ell=1,\dots,s$ and $I\subseteq \sspi_\ell$, if $I\subseteq
  \varpi!_{\ell,1}$, then
  \[
  \ssub{(\mu!_\pi)}!_{\pi'}(I)=\mu!_\pi(I)=
  \mu!_{\filtind_\ell/\filtind_{\ell-1}}(I)=
  \mu!_{\filtind_{\ell}^-/\filtind_{\ell-1}}(I)=\mu! _{\pi''}(I).
  \]
On the other hand, if $I\subseteq \varpi!_{\ell,2}$, then, using
Proposition~\ref{muprod} for the second equality below,
\[
  \ssub{(\mu!_\pi)}!_{\pi'}(I)=\mu!_{\pi}(I\cup\sspi_1')-\mu!_{\pi}(\sspi_1')=
  \mu!_{\pi}(I\cup\varpi!_{\ell,1})-\mu!_{\pi}(\varpi!_{\ell,1})
  =\mu(I\cup \filt^-_{\ell})-\mu(\filt^-_{\ell})
  =\mu!_{\pi''}(I).
\]
The result follows.
\end{proof}

\subsection{Capcup lemma} The following key lemma will be quite useful.
\begin{lemma}\label{capcup} Let $\mu\in\M(V)$ be a supermodular function and $\ssI_1,\ssI_2\subseteq V$. If there is $q\in\P_{\mu}$ such that $q(\ssI_1)=\mu(\ssI_1)$ and $q(\ssI_2)=\mu(\ssI_2)$, then 
\begin{equation}\label{mu=}
  \mu(\ssI_1)+\mu(\ssI_2)= \mu(\ssI_1\cup \ssI_2)+ \mu(\ssI_1\cap \ssI_2).
\end{equation}
Conversely, if the above equality holds, then, for each 
$x\in\P_{\mu}$, the following statements are equivalent:
\begin{enumerate}[label=(\arabic*)]
  \item\label{capcup1} $x(\ssI_1)=\mu(\ssI_1)$
    and $x(\ssI_2)=\mu(\ssI_2)$.
  \item\label{capcup2} $x(\ssI_1)+x(\ssI_2)=\mu(\ssI_1)+\mu(\ssI_2)$.
    \item\label{capcup3} $x(\ssI_1\cup \ssI_2)+x(\ssI_1\cap \ssI_2) =\mu(\ssI_1\cup \ssI_2)+\mu(\ssI_1\cap \ssI_2)$.
    \item\label{capcup4} $x(\ssI_1\cup \ssI_2)=\mu(\ssI_1\cup \ssI_2)$ and
      $x(\ssI_1\cap \ssI_2)=\mu(\ssI_1\cap \ssI_2)$.
      \end{enumerate}
\end{lemma}

\begin{proof} For each $x\in\P_{\mu}$,
\[
x(\ssI_1)+x(\ssI_2)=x(\ssI_1\cup \ssI_2)+x(\ssI_1\cap \ssI_2).
\]
Assume there is $q\in\P_{\mu}$ such that $q(\ssI_1)=\mu(\ssI_1)$
and $q(\ssI_2)=\mu(\ssI_2)$. Then,
\[
q(\ssI_1\cup \ssI_2)+q(\ssI_1\cap \ssI_2) \leq
\mu(\ssI_1\cup \ssI_2)+\mu(\ssI_1\cap \ssI_2).
\]
Also, $q(\ssI_1\cap \ssI_2)\geq\mu(\ssI_1\cap \ssI_2)$
and $q(\ssI_1\cup \ssI_2)\geq\mu(\ssI_1\cup \ssI_2)$.
Thus, equalities hold everywhere.

Let $x\in\P_{\mu}$. If \eqref{mu=} holds, then \ref{capcup2} and
\ref{capcup3} are equivalent. Clearly, \ref{capcup1} implies \ref{capcup2}, and \ref{capcup4} implies \ref{capcup3}. But
$x(J)\geq\mu(J)$ for each $J\subseteq V$. Therefore, \ref{capcup2} implies \ref{capcup1}, and
\ref{capcup3} implies \ref{capcup4}. The lemma follows.
\end{proof}
 
\subsection{Face structure}  In this section, we give a
characterization of the face structure of $\P_{\mu}.$

\begin{prop}\label{splitpoly}  Let $\mu\in\M(V)$ be a supermodular
  function and $\pi=(\sspi_1, \dots, \sspi_s)$ an ordered partition of $V$. Let $\filt_\bullet$ be the corresponding filtration. Then, we have
  \[
  \P_{\mu_\pi}=\P_{\mu}\cap\left\{q\in\R^V\,\st\,q(\filt_j)=\mu(\filt_j)\,\text{\rm
  for }j=1,\dots,s-1\right\}.
  \]
\end{prop}

\begin{proof} Proceeding by induction, using
  Proposition~\ref{refineV}, we may assume $\pi=(\sspi_1,\sspi_2)$. For each $I\subseteq V$, we have
\begin{equation}\label{mumuV}
\mu!_\pi(I)=\mu!_{\pi_1/\emptyset}(I)+\mu!_{V/\pi_1}(I)=
\mu(I\cap \sspi_1)+\mu(I\cup \sspi_1)-\mu(\sspi_1).
\end{equation}
It follows from Proposition~\ref{muprod} that $\mu!_{\pi}(\sspi_1)=\mu!_{\sspi}^*(\sspi_1)$, whence
\[
\P_{\mu_\pi}\subseteq\P_{\mu}\cap\left\{q\,\st\,q(\sspi_1)=\mu(\sspi_1)\right\}.
\]
So we only need to prove the reverse inclusion. Let thus $q\in \P_{\mu}$ be such that
$q(\sspi_1)=\mu(\sspi_1)$. For each $I\subseteq V$, since $q\in \P_{\mu}$, it
follows from \eqref{mumuV} that 
\[
q(I)+q(\sspi_1)=q(I\cap \sspi_1)+q(I\cup \sspi_1)\geq\mu(I\cap
\pi_1)+\mu(I\cup \sspi_1)=\mu!_\pi(I)+\mu(\sspi_1),
\]
and so $q(I)\geq\mu!_\pi(I)$. We conclude that $q\in\P_{\mu_\pi}$. 
\end{proof}

\begin{defi}[Codimension of a supermodular function] Let $\mu\in\M(V)$. Denote by $\cd!_{\mu}$ the largest integer $s$ for
which there is an $s$-partition $\pi$ of $V$ such that
$\mu=\mu!_\pi$. We call $\cd!_{\mu}$ the \emph{codimension} of $\mu$.
\end{defi}

The following proposition justifies the name given to $\cd!_{\mu}$.

\begin{prop}\label{splitpoly2} Let $\mu\in\M(V)$ be a supermodular
  function. Then, the base polytope $\P_{\mu}$ is nonempty and we can
  recover $\mu$ (and $\mu!^*$) from $\P_{\mu}$ as follows:
\begin{equation}\label{Ptomu}
  \mu(I)=\min\left\{q(I)\,\st\,q\in \P_{\mu}\right\}\quad\text{and}\quad
  \mu!^*(I)=\max\left\{q(I)\,\st\,q\in \P_{\mu}\right\}\quad\text{for each }I\subseteq V.
\end{equation}
In particular, we have $\P_{\mu_1}=\P_{\mu_2}$
if and only if $\mu!_1=\mu!_2$. In addition, the base polytopes $\P_{\mu_\pi}$ for ordered partitions $\pi$ of $V$ are
all the faces of $\P_{\mu}$.   Finally, we have $\dim\P_{\mu}+\cd!_{\mu}=|V|$.
\end{prop}

\begin{proof} That $\P_{\mu}\neq\emptyset$ follows from
  Proposition~\ref{splitpoly}. In fact,
the latter states that $\P_{\mu_\pi}\subseteq\P_{\mu}$ for each ordered partition
$\pi$ of $V$. And, by Proposition~\ref{muprod}, $\P_{\mu_\pi}$ is a point in the case $\pi$ is maximal among nontrivial partitions.

Furthermore, Proposition~\ref{splitpoly} yields the inequality $q(I)\geq\mu(I)$ for each
$I\subseteq V$ and $q\in
\P_{\mu}$, with equality if $q\in \P_{\mu_{(I,\ssI^c)}}$. Since
$\P_{\mu_{(I,\ssI^c)}}\neq\emptyset$, we conclude that $\mu(I)=\min\left\{q(I)\,\st\,q\in
  \P_{\mu}\right\}$ for each $I\subseteq V$.  The fact that $\mu!^*(I)$ is the maximum of $q(I)$ for $q \in \P_{\mu}$
  follows from the equality $\mu(\ssI^c)+\mu!^*(I)=\mu(V)$.

   The statement that each $\P_{\mu_\pi}$ is a face of $\P_{\mu}$ follows from
  Proposition~\ref{splitpoly}, as $\P_{\mu_\pi}\neq\emptyset$
  and  $q(\filt_j)\geq\mu(\filt_j)$ for each $j$ and $q\in\P_{\mu}$.
  Conversely, let $\F$ be a face of
  $\P_{\mu}$.  Let $\pi = (\sspi_1, \dots, \sspi_s)$ be a maximal nontrivial ordered partition
 for the property that $\F\subseteq\P_{\mu_\pi}$. 
 We claim the equality holds. Indeed, if not, then there is a
 nonempty proper subset $I\subset V$ which is not a (partial) union of $\sspi_j$ such that $q(I)=\mu(I)$ for each  $q\in\F$. Let $\filt_\bullet$ be the filtration associated to $\pi$. Let
$i\in\mathbb N$ be the smallest integer for which
$I\cap \filt_i$ is not a partial union of the $\sspi_j$. Then, $I\cap \sspi_i$ is
neither $\emptyset$ nor $\sspi_i$. Thus
$\filt_{i-1}\subsetneq \filt_{i-1}\cup (I\cap \filt_i)\subsetneq  \filt_i$.
Since  by Proposition~\ref{capcup}, $q(\filt_{i-1}\cup (I\cap  \filt_i))=\mu(\filt_{i-1}\cup (I\cap  \filt_i))$ for
each $q\in\F$, we get a nontrivial refinement of $\pi$
satisfying the same property as $\pi$. This contradiction proves the equality $\F = \P_{\mu_\pi}$.

Finally, assume $\F =\P_{\mu}$. Let $\pi$ be a $\cd!_{\mu}$-partition of $V$
such that $\mu=\mu!_\pi$. Then, $\pi$ is maximal for the property that
$\P_{\mu_\pi} = \P_{\mu}$.  Let $\dd!_{\mu} \coloneqq\dim\P_{\mu}$. We have
$\dd!_{\mu}\leq |V|-\cd!_{\mu}$, or $\cd!_{\mu}\leq|V|-\dd!_{\mu}$. We claim the equality $\dd!_{\mu}=|V|-\cd!_{\mu}$ holds. Indeed, if
not, there would be a
 (nonempty proper) subset $I\subset V$ which is not a partial union of
the $\pi_j$ such that $q(I)=\mu(I)$ for each  $q\in\F$. As before, this
contradicts the maximality of $\pi$. We conclude that $\cd!_{\mu}=|V|-\dd!_{\mu}$. 
  \end{proof}

\subsection{Base polytope of a positive combination of modular pairs} Given two modular pairs $(\mu, \mu!^*)$ and $(\nu,\nu!^*)$, and two positive real numbers $\alpha$ and $\beta$, the pair $(\alpha\mu+\beta\nu, \alpha\mu!^*+\beta\nu!^*)$ is easily seen to be modular. 

\begin{prop}\label{prop:positive-combination} The base polytope $\P_{\alpha\mu+\beta\nu}$ coincides with $\alpha\P_{\mu} + \beta\P_{\nu}$.
\end{prop}

\begin{proof} The inclusion $\alpha\P_{\mu} + \beta\P_{\nu} \subseteq\P_{\alpha\mu+\beta\nu}$ is obvious.

For the reverse inclusion, we will proceed by induction on the dimension of $\P_{\alpha\mu+\beta\nu}$. If that dimension is zero, then clearly the inclusion 
$\alpha\P_{\mu} + \beta\P_{\nu} \subseteq\P_{\alpha\mu+\beta\nu}$ must be an equality. 

Suppose the dimension is nonzero. Since $\alpha\P_{\mu}+\beta\P_{\nu}$ is convex, and $\P_{\alpha\mu+\beta\nu}$ is the convex hull of its proper faces, it is enough to show that 
\[
\P_{\alpha\mu_{\pi}+\beta\nu_{\pi}}\subseteq
\alpha\P_{\mu}+\beta\P_{\nu}
\]
for each ordered partition $\pi$ of $V$ such that 
$\alpha\mu+\beta\nu$ is not $\pi$-split. This follows from the induction hypothesis, which yields
\[
\P_{\alpha\mu_{\pi}+\beta\nu_{\pi}}=
\alpha\P_{\mu_{\pi}}+\beta\P_{\nu_{\pi}}\subseteq\alpha\P_{\mu}+\beta\P_{\nu}. \qedhere
\]
\end{proof}


\section{Families of modular pairs and their tilings} \label{sec:admissible-families-tilings}

In this section, we introduce a general framework that produces tilings induced by collections of modular pairs.

For a polytope $\P$ in the real vector space $H=\R^V$, we denote by $\mathring \P$ the
\emph{relative interior} of $\P$ defined as the set of all points $q
\in \P$ that do not belong to any proper face of $\P$.

\subsection{Simple supermodular functions}

\begin{defi}[Simple supermodular function]
An element $\mu\in\M(V)$ is called \emph{simple} if there is no nontrivial ordered
partition $\pi$ of $V$ for which $\mu=\mu!_\pi$.
\end{defi}

 Notice that the function $\mu$ is simple if and only if
 $\mu(I)+\mu(J)<\mu!(V)$ for each bipartition $(I,J)$ of $V$, or equivalently by Proposition~\ref{muprod}, $\mu(I)<\mu!^*(I)$ for each nonempty proper subset $I$ of $V$, or equivalently by Proposition~\ref{splitpoly2}, $\P_{\mu}$ has (maximum) dimension $|V|-1$.

 \begin{defi}[Spread] For $\mu \in \M(V)$, the minimum of $\mu!^*(I)-\mu(I)$ for nonempty proper subsets $I$ of $V$ is called the \emph{spread} of $\mu$. By convention, the spread is zero if $V$ consists of a single element. 
 \end{defi}

 Note that the spread is positive if and only if $\mu$ is simple. 

\subsection{Separation}

Let $\mu,\nu\in \M_n(V)$ be two supermodular functions of range $n$,
and let $\pi=(I,J)$ be a bipartition of $V$.

\begin{defi} We say that the bipartition $\pi$ is a \emph{separation}
  for the ordered pair $(\mu, \nu)$ if $\mu(I) + \nu(J)  \geq n$.
  We say that the separation is \emph{strict} if the inequality is
strict. We say that the separation is \emph{nontrivial} provided that it is either strict, or 
we have $\mu!_\pi \neq \mu$ or $\nu!_\pi \neq \nu$.
\end{defi}

Alternatively, $\pi$ is a separation for $(\mu, \nu)$ if and only if $\mu(I)\geq\nu!^*(I)$, or equivalently, $\nu(J)\geq\mu!^*(J)$. It follows that $\pi$ is a separation for $(\mu, \nu)$ if and only if $\ssq_1(I)\geq \ssq_2(I)$, equivalently, $\ssq_2(J)\geq
\ssq_1(J)$, for each $\ssq_1\in\P_{\mu}$ and $\ssq_2\in\P_{\nu}$. Furthermore,
the separation is strict if and only if these inequalities are strict for each $\ssq_1\in\P_{\mu}$ and
$\ssq_2\in\P_{\nu}$, and it is nontrivial if and only if the inequalities are strict for
certain $\ssq_1\in\P_{\mu}$ and $\ssq_2\in\P_{\nu}$.

Note that $\pi=(I,J)$ is a (strict, resp.~nontrivial) separation for $(\mu,\nu)$  if and only if $\sspi^c=(J,I)$ is a (strict, resp.~nontrivial) separation for $(\nu, \mu)$. Therefore, the existence of a (strict, resp.~nontrivial) separation for a pair of elements in $\M_n(V)$ does not depend on the order of the pair.

Given two simple elements $\mu,\nu\in\M_n(V)$, we say that $\P_{\mu}$ and
$\P_{\nu}$ \emph{intersect in
codimension at most 1} if there are faces $\F_1$ and $\F_2$ of $\P_{\mu}$
and $\P_{\nu}$ respectively of codimension at most 1 such that
$\ssub{\mathring\F}!_1\cap\ssub{\mathring\F}!_2\neq\emptyset$.  

\begin{prop}\label{uniquebipart}
  Let $\mu,\nu\in \M_n(V)$ be two supermodular functions of
  range $n$. 
  \begin{enumerate}[label=(\Roman*)]
  \item \label{prop:part1-sep}For each bipartition $\pi=(I,J)$ of $V$ the following
  statements hold:
  \begin{enumerate}[label=(\arabic*)]
    \item\label{prop:part1-1} If $\pi$ is a separation for $(\mu, \nu)$, then
  \begin{equation}\label{capmunu}
  \P_{\mu}\cap\P_{\nu}=\P_{\mu_\pi}\cap\P_{\nu_{\sspi^c}}.
\end{equation}
\item\label{prop:part1-2} If $\P_{\mu_\pi}\cap\P_{\nu_{\sspi^c}}\neq\emptyset$, then $\pi$
  is a nonstrict separation for $\mu$ and $\nu$.
  \item\label{prop:part1-3} If $\pi$ is a nontrivial separation for $\mu$ and $\nu$, then
    $\mathring\P_{\mu}\cap\mathring\P_{\nu}=\emptyset$.
\end{enumerate}
   \item \label{prop:part2-sep} Finally, if $\mu$ and $\nu$ are simple, then $\P_{\mu}$ and 
$\P_{\nu}$ intersect in
codimension at most $1$ only if there is at most one separation for
$\mu$ and $\nu$. 
\end{enumerate}
\end{prop}

\begin{proof} As argued before, if $\pi$ is a 
  separation for $\mu$ and $\nu$, then $\ssq_1(I)\geq\mu(I)\geq\nu!^*(I)\geq \ssq_2(I)$ for
  each $\ssq_1\in\P_{\mu}$ and $\ssq_2\in\P_{\nu}$, whence \eqref{capmunu}
  holds and \ref{prop:part1-1} follows. On the other hand, if $\P_{\mu_\pi}\cap\P_{\nu_{\sspi^c}}$ is
  nonempty, then for each $q$ in the intersection,  we have $\mu(I)=q(I)=\nu!^*(I)$. Therefore, $\pi$ is a nonstrict separation for
  $\mu$ and $\nu$ and \ref{prop:part1-2} follows. To prove \ref{prop:part1-3}, if $\pi$ is a strict separation for $\mu$ and
  $\nu$, then it follows from \ref{prop:part1-2} that $\P_{\mu}\cap\P_{\nu}=\emptyset$. In the remaining case, since $\pi$ is a nontrivial separation for $\mu$ and $\nu$, we have either
  $\mu!_\pi\neq\mu$ or $\nu!_\pi\neq\nu$. Hence, \eqref{capmunu}
  implies that $\mathring\P_{\mu}\cap \P_{\nu}=\emptyset$ or
  $\P_{\mu}\cap\mathring\P_{\nu}=\emptyset$, and thus
  $\mathring\P_{\mu}\cap \mathring\P_{\nu}=\emptyset$. This proves 
  Part~\ref{prop:part1-sep} of the proposition.
  
  Assume now that $\mu$ and $\nu$ are simple. Then, every separation
  for $\mu$ and $\nu$ is nontrivial. Assume that $\P_{\mu}$ and
$\P_{\nu}$ intersect in codimension at most $1$. Since $\mu$ and $\nu$ are simple, if either $\mathring\P_{\mu}\cap \P_{\nu}\neq\emptyset$ or  $\P_{\mu}\cap\mathring\P_{\nu}\neq\emptyset$, then
$\mathring\P_{\mu}\cap \mathring\P_{\nu}\neq\emptyset$. In this case
there is no separation for $\mu$ and $\nu$ by Part~\ref{prop:part1-sep}.

Assume $\mathring\P_{\mu}\cap \mathring\P_{\nu}=\emptyset$. Then, there
are facets $\F_1$ and $\F_2$ of $\P_{\mu}$
and $\P_{\nu}$ respectively such that
$\ssub{\mathring\F}!_1\cap\ssub{\mathring\F}!_2\neq\emptyset$. There is a bipartition
$\pi$ of $V$ such that $\F_1=\P_{\mu_\pi}$. Since
$\mathring\P_{\mu}\cap \mathring\P_{\nu}=\emptyset$, we must have that
$\F_2=\P_{\nu_{\sspi^c}}$, and hence $\pi$ is a nonstrict separation
for $\mu$ and $\nu$ by Part~\ref{prop:part1-sep}.

Let $\lambda\coloneqq (I',J')$ be a
  separation for $\mu$ and $\nu$. Then,
  $\P_{\mu}\cap\P_{\nu}=\P_{\mu_\lambda}\cap\P_{\nu_{\lambda!^c}}$,
  and hence $\mathring \P_{\mu_{\pi}}\cap \mathring
  \P_{\nu_{\sspi^c}}\subseteq \P_{\mu_\lambda}$. Then, the span of
  $\P_{\mu_{\pi}}$
  is contained
  in the subspace of
  $\R^V$ given by $q(I')=\mu(I')$ and
  $q(J')=\mu!^*(J')$. As the span has codimension 2, the two subspaces are
  the same. If follows that $I=I'$ or $I=J'$. But $I=J'$ yields
  $\mu(I)=\mu!^*(I)$, contradicting the assumption that $\mu$ is
  simple. Thus, $I=I'$ and hence $\lambda=\pi$, as required. 
  \end{proof}

\subsection{Families}
Let $\C \subset \M_n(V)$ be a family of supermodular functions of
range $n$. We say that:
\begin{itemize}
    \item $\C$ is \emph{separated} if any two distinct elements
$\mu, \nu \in \C$ admit a nontrivial separation,
    \item $\C$ is \emph{closed} if each splitting of each element of $\C$ is in $\C$,
    \item $\C$ is \emph{simple} if each element of $\C$ is a splitting of a simple element of $\C$,
    \item $\C$ is \emph{complete} if for each $\mu\in\C$ and each bipartition $\pi$ of $V$, there is $\mu'\in\C$ distinct from $\mu!_\pi$ such that $\mu!_{\sspi^c}'=\mu!_{\pi}$,
    \item $\C$ is \emph{positive} if the infimum of the spreads of its simple elements is positive.
\end{itemize}
We denote by $\P_{\C}$ the union of the polytopes $\P_{\mu}$ for $\mu \in \C$ and call it the \emph{support of the collection  $\C$}.

\smallskip

Here is a relative variant. Let $\Omega\subseteq H_n$ be a subset. We say that:
\begin{itemize}
    \item $\C$ is \emph{simple for $\Omega$} if each element $\mu\in\C$ such that $\P_{\mu}$ intersects $\Omega$ is a splitting of a simple element of $\C$,
    \item $\C$ is \emph{complete for $\Omega$} if for each $\mu\in\C$ and each bipartition $\pi$ of $V$ such that $\P_{\mu_\pi}$ intersects $\Omega$, there is $\mu'\in\C$ distinct from $\mu!_\pi$ such that $\mu!_{\sspi^c}'=\mu!_{\pi}$.
\end{itemize}
We keep the other definitions unchanged in the relative setting. 

Clearly, $\C$ is simple (resp.~complete) if and only if $\C$ is simple (resp.~complete) for $H_n$. Also, if $\C$ is closed, then $\C$ is complete for $\Omega$ if for each $\mu\in\C$ and each bipartition $\pi$ of $V$ such that $\mu=\mu!_{\pi}$ and  $\P_{\mu}$ intersects $\Omega$, there is $\mu'\in\C$ distinct from $\mu$ such that $\mu!_{\pi}'=\mu$. Finally, notice that $\C$ can only be positive if it contains simple elements.

\begin{prop}\label{prop:fullissimple} Let $\Omega\subseteq H_n$ be a subset. If $\C$ is complete for $\Omega$, then $\C$ is simple for $\Omega$. 
\end{prop}

\begin{proof} Let $\mu\in\C$ such that $\P_{\mu}$ intersects $\Omega$. Assume $\mu$ is not simple. Then, there is a bipartition $\pi$ of $V$ such that $\mu=\mu!_\pi$. If $\C$ is complete for $\Omega$, there is $\mu'\in C$ distinct from $\mu!_{\pi}$ such that $\mu!_{\sspi^c}'=\mu$. Since $\mu'\neq\mu$, we have that $\mu$ is a nontrivial splitting of $\mu'$. 

If $\mu'$ is simple, we stop. If not, since $\P_{\mu'}$ intersects $\Omega$, we may repeat the above process for $\mu'$ instead of $\mu$. We will eventually end up with a simple $\nu\in\C$ of which $\mu$ is a splitting.
\end{proof}

Notice that the converse does not hold: the collection of splittings of a simple supermodular function is closed, simple and positive but is not complete. 

\subsection{Tilings} 
 Let $\C \subset \M_n(V)$ be a family of supermodular functions of range $n$. 
\begin{thm}[Tilings induced by separated and closed families]\label{thm:tiling} If $\C \subseteq \M_n(V)$ is separated and closed, then the collection of polytopes $\P_{\mu}$ for $\mu \in \C$ gives a tiling of $\P_{\C}$.
\end{thm}

 \begin{proof} We need to prove
   that for each pair $\mu,\nu \in
   \C$, the base polytopes $\P_{\mu}$ and $\P_{\nu}$ are either disjoint, or they intersect in a common face. We can suppose the intersection is nonempty. Let $q$ be
   a point in the intersection of $\P_{\mu}$ and
   $\P_{\nu}$. Let $\F$ and $\F'$ be the faces of $\P_{\mu}$ and
   $\P_{\nu}$, respectively, which contain $q$ in their relative
   interiors. Denote by $\pi$ and $\pi'$ the two ordered partitions of
   $V$ such that $\F=\P_{\mu_\pi}$ and $\F'=\P_{\nu_{\pi'}}$. Since
   $\C$ is closed, both $\mu!_\pi$ and $\nu!_{\pi'}$ belong to $\C$.
   Since $\C$ is separated, we should have $\mu!_\pi=\nu!_{\pi'}$, otherwise, the
   relative interiors of $\F$ and $\F'$ would be disjoint by
   Proposition~\ref{uniquebipart}. We infer that $\F =\F'$. 
 
 We have shown that the two polytopes $\P_{\mu}$ and $\P_{\nu}$ intersect in the union of their common faces. Since their intersection is convex, this implies that they intersect in a common face.  
\end{proof}

We need the following lemma for our next result.

\begin{lemma}\label{lem:fillHn} Let $\Omega$ be a connected open subset of $H_n$. Let $\mathcal P$ be a nonempty collection of full-dimensional polytopes in $H_n$ whose union intersects $\Omega$. Assume that $\mathcal P$ is locally finite, that is, each ball of $H_n$ contains finitely many polytopes of $\mathcal P$. Also, assume that each facet intersecting $\Omega$ of a polytope in $\mathcal P$ is shared with another polytope in $\mathcal P$. Then, the support of $\mathcal P$ contains $\Omega$.  
\end{lemma}

\begin{proof} Let $\Theta$ be the union of the polytopes $\P$ in $\mathcal P$, and $q\in\Omega$. By hypothesis, $\Theta\cap\Omega\neq\emptyset$. Let $\gamma$ be a piecewise linear path in $\Omega$ that connects a point of $\Theta \cap\Omega$ to $q$. Let $\mathcal P'\subseteq\mathcal P$ be the collection of polytopes that intersect $\gamma$, and denote by $\Theta'$ the union of polytopes in $\mathcal P'$. Clearly, $\Theta'\neq\emptyset$. Also, since $\mathcal P$ is locally finite, $\mathcal P'$ is finite. Furthermore, since $\Omega$ is open, modifying $\gamma$ if necessary, we may assume that it intersects the boundary of $\P$ away from its faces of codimension 2 for each $\P\in\mathcal P'$. Since $\mathcal P'$ is finite, there is a point $q'\in \Theta'\cap\gamma$ closest to $q$ along $\gamma$. Suppose by contradiction that $q'\neq q$. Then, $q'$ is on a facet of some $\P\in\mathcal P'$. By hypothesis, there is $\P'\in\mathcal P$ sharing that facet with $\P$. Since $q'\in\gamma$, also $\P'\in\mathcal P'$. But then there is an open ball of $q'$ in $H_n$ contained in $\Theta'$. This contradicts the choice of $q'$ as closest point to $q$.
\end{proof}

\begin{thm}[Tilings induced by complete and positive families]\label{thm:tiling2} Let $\Omega$ be a connected open subset of $H_n$. Assume $\P_{\C}$ intersects $\Omega$. If $\C$ is complete for $\Omega$ and positive, then $\P_{\C}\supseteq\Omega$.
\end{thm}

\begin{proof} Let $\mathcal P$ be the collection of polytopes $\P_{\mu}$ for simple $\mu\in\C$ intersecting $\Omega$. Since $\C$ is complete for $\Omega$, then $\C$ is simple for $\Omega$ by Proposition~\ref{prop:fullissimple}. Since $\P_{\C}\cap\Omega\neq\emptyset$, it follows that $\mathcal P$ is nonempty and the union of its polytopes intersect $\Omega$. 

Since $\C$ is positive, the minimum of the volumes of the (full-dimensional) polytopes $\P\in\mathcal P$ is positive. This implies that $\mathcal P$ is locally finite.

Finally, let $\mu\in\C$ such that $\P_{\mu}\in\mathcal P$. Consider a facet of $\P_{\mu}$ intersecting $\Omega$. It coincides with $\P_{\mu_{\pi}}$ for a bipartition $\pi$ of $V$. Since $\C$ is complete for $\Omega$, there is $\mu'\in\C$ distinct from $\mu!_{\pi}$ such that $\mu!_{\pi}=\mu!_{\sspi^c}'$. Thus, $\P_{\mu_{\pi}}$ is a face of $\P_{\mu'}$. Now, $\P_{\mu_{\pi}}\neq\P_{\mu'}$, because $\mu'\neq\mu!_{\pi}$. Then, $\P_{\mu_{\pi}}$ is a proper face of $\P_{\mu'}$. Therefore, $\P_{\mu'}$ is full dimensional, and $\P_{\mu_{\pi}}$ is a facet of $\P_{\mu'}$. In particular, $\P_{\mu'}\in\mathcal P$. The theorem follows now by applying Lemma~\ref{lem:fillHn}.   
\end{proof}


\section{Admissible divisors on metric graphs and their semistability polytopes} \label{sec:admissible-divisor-semistability-polytope}

 We refer to the survey paper~\cite{BJ} for basic results on algebraic geometry of graphs and metric graphs, and their link to tropical and non-Archimedean geometry.

\subsection{Admissible divisors}\label{sec:admissible}

 Let $G=(V,E)$ be a finite graph and $\ell \colon E \to (0, +\infty)$ be an edge length function. We denote by $\ell_e$ the real number $\ell(e)$. We denote by $\Gamma$ the metric graph obtained by the metric realization of the pair $(G, \ell)$: this is obtained by plugging an interval of length $\ell_e$ between the two extremities of the edge $e$. Endowing $\Gamma$ with the path metric gives the metric space $\Gamma$. The pair $(G,\ell)$ is called a model for $\Gamma$.

 The group of divisors on $\Gamma$, denoted $\Div(\Gamma)$, is by definition the free Abelian group generated by the points on $\Gamma$. Writing $(x)$ for the generator associated to the point $x\in \Gamma$, 
 \[\Div(\Gamma) \coloneqq\Bigl\{\sum_{x\in A \subset\Gamma} n_x (x)\,\, \st \,\, n_x\in \mathbb Z\, \textrm{ and } A \textrm{ a finite set} \Bigr\}.\]
 
For a divisor $D\in \Div(\Gamma)$ and $x\in \Gamma$, the coefficient of $(x)$ in $D$ is
denoted $D(x)$. Its support, $\mathrm{Supp}(D)$, is the set of points $x$
with $D(x)\neq 0$. And its degree, $\mathrm{deg}(D)$, is defined as 
 \[ \mathrm{deg}(D)\coloneqq\sum_{x\in \Gamma} D(x).\]
 The subset of $\Div(\Gamma)$ consisting of divisors of degree $d$ is denoted by $\Div^d(\Gamma)$.
 
 If $I$ is a subset of $V$, we denote by $E[I]$ the set of edges
 of $G$ with both extremities in $I$. The cardinality of $E[I]$ will
 be denoted by $e[I]$. We call the subgraph $G[I] \coloneqq (I, E[I])$ the induced subgraph on $I$. We denote by $\ssGamma_I$ the metric graph with model $(G[I], \ell\rest{E[I]})$. For a divisor $D$ on $\Gamma$, and a subset $I \subset V$, we denote by $\ssD_I$ the divisor on $\ssGamma_I$ defined by
 \[
 \ssD_I\coloneqq \sum_{x\in \ssGamma_I} D(x)(x).
 \]
\begin{defi}[$G$-admissible divisor]\label{def: admissible-divisor}
 A divisor $D$ on $\Gamma$  is called \emph{$G$-admissible} if for each
 edge $e$ of $G$, the coefficient of $D$ at any point in the interior of $e$ is $0$ with at most one possible exception $x$ for which $D(x)=1$. 
\end{defi}

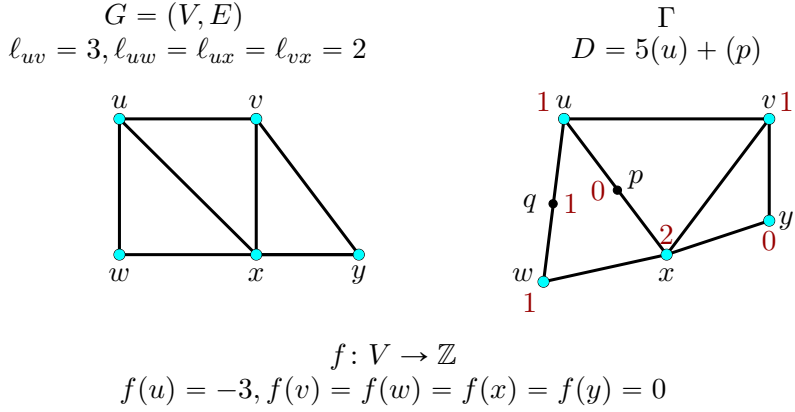
\begin{figure}
    \centering
\begin{tikzpicture}[scale=0.9]
\draw[line width=0.4mm] (5,-2) -- (7,-2);
\draw[line width=0.4mm] (5,0) -- (7,0);
\draw[line width=0.4mm] (5,0) -- (7,-2);
\draw[line width=0.4mm] (7,-2) -- (7,0);
\draw[line width=0.4mm] (5,-2) -- (5,0);
\draw[line width=0.4mm] (8.5,-2) -- (7,0);
\draw[line width=0.4mm] (7,-2) -- (8.5,-2);
\draw[line width=0.4mm] (8.5,-2) -- (7,-2);

\draw[line width=0.1mm] (5,0) circle (0.8mm);
\filldraw[aqua] (5,0) circle (0.7mm);
\draw[line width=0.1mm] (7,0) circle (0.8mm);
\filldraw[aqua] (7,0) circle (0.7mm);
\draw[line width=0.1mm] (5,-2) circle (0.8mm);
\filldraw[aqua] (5,-2) circle (0.7mm);
\draw[line width=0.1mm] (7,-2) circle (0.8mm);
\filldraw[aqua] (7,-2) circle (0.7mm);
\draw[line width=0.1mm] (8.5,-2) circle (0.8mm);
\filldraw[aqua] (8.5,-2) circle (0.7mm);

\draw[black] (5,0.25) node{$u$};
\draw[black] (7,0.25) node{$v$};
\draw[black] (5,-2.3) node{$w$};
\draw[black] (7,-2.3) node{$x$};
\draw[black] (8.5,-2.3) node{$y$};
\draw[black] (5.8,1.5) node{$G=(V,E)$};
\draw[black] (6,1) node{$\ell_{uv}=3, \ell_{uw}=\ell_{ux}=\ell_{vx}=2$};

\draw [line width=0.4mm](11.5,0) -- (14.5,0);
\draw [line width=0.4mm](11.5,0) -- (11.2,-2.4);
\draw [line width=0.4mm](13,-2) -- (11.2,-2.4);
\draw [line width=0.4mm](13,-2) -- (11.5,0);
\draw [line width=0.4mm](14.5,-1.5) -- (13,-2);
\draw [line width=0.4mm](14.5,-1.5) -- (14.5,0);
\draw [line width=0.4mm](13,-2) -- (14.5,0);

\draw[line width=0.1mm] (11.5,0) circle (0.8mm);
\filldraw[aqua] (11.5,0) circle (0.7mm);

\filldraw[black] (12.28,-1.05) circle (0.6mm);
\filldraw[black] (11.34,-1.25) circle (0.6mm);

\draw[line width=0.1mm] (14.5,0) circle (0.8mm);
\filldraw[aqua] (14.5,0) circle (0.7mm);
\draw[line width=0.1mm] (14.5,-1.5) circle (0.8mm);
\filldraw[aqua] (14.5,-1.5) circle (0.7mm);
\draw[line width=0.1mm] (13,-2) circle (0.8mm);
\filldraw[aqua] (13,-2) circle (0.7mm);
\draw[line width=0.1mm] (11.2,-2.4) circle (0.8mm);
\filldraw[aqua] (11.2,-2.4) circle (0.7mm);

\draw[black] (11.5,0.25) node{$u$};
\draw[darkred] (11.2,0.25) node{$1$};

\draw[black] (14.5,0.25) node{$v$};
\draw[darkred] (14.75,0.25) node{$1$};

\draw[black] (10.9,-2.3) node{$w$};
\draw[darkred] (11,-2.7) node{$1$};

\draw[black] (13,-2.3) node{$x$};
\draw[darkred] (13,-1.7) node{$2$};

\draw[black] (12.55,-.9) node{$p$};
\draw[darkred] (12,-1.05) node{$0$};

\draw[black] (11,-1.25) node{$q$};
\draw[darkred] (11.6,-1.25) node{$1$};

\draw[black] (14.75,-1.5) node{$y$};
\draw[darkred] (14.5,-1.8) node{$0$};

\draw[black] (13,1.5) node{$\Gamma$};
\draw[black] (13,1) node{$D=5(u)+(p)$};

\draw[black] (9,-3.5) node{$f\colon V \to \mathbb Z$};
\draw[black] (9,-4) node{$f(u)=-3, f(v)=f(w)=f(x)=f(y)=0$};
\end{tikzpicture} 
\caption{A graph $G=(V,E)$ and an edge length function $\ell$ are described on the left; $\ell$ assigns length 1 to the three edges $wx$, $xy$ and $vy$. The corresponding metric graph $\Gamma$ is depicted on the right. The divisor $D$ has coefficient $5$ at $u$ and coefficient $1$ at the point $p$ in the middle of the edge $ux$. The function $f\colon V\to \R$ described in the bottom gives the $G$-admissible divisor $D' = D+\div_\ell(f;D)$ on $\Gamma$ depicted in red, with coefficient 1 at $u, v, w$, coefficient 1 at the point $q$ in the middle of $uw$, and coefficient 2 at $x$.}
\label{fig:admissible}
\end{figure}

An example of a $G$-admissible divisor is depicted in Figure~\ref{fig:admissible}. The $v$-reduced divisors are also examples of $G$-admissible divisors; see Proposition~\ref{prop:v-reduced-adm}.

\begin{defi}[Spanning subgraphs of $G$ associated to $G$-admissible divisors] Given a $G$-admissible divisor $D$ on $\Gamma$, we denote by $\ssG_D$ the spanning subgraph of $G$ whose edge set $\ssE_D$ consists of all the edges $e\in E$ whose interiors contain no point on the support of $D$.
\end{defi}
  
Given subsets $I,J \subseteq V$, we denote by $\ssE_D(I,J)$ the set of edges in $\ssG_D$ with one extremity in $I$ and the other  in $J$. The cardinality of $\ssE_D(I,J)$ will be denoted $\sse_D(I,J)$.  For $D=0$, we remove the subscript $D$. Thus, $E(I,J)$ will be the set of edges between $I$ and $J$ in $G$, and $e(I,J)$ will be the number of these edges. Then, we put
\[
 \delta!_D(I,J)\coloneqq \sse(I,J) - \sse_D(I,J),
 \] 
 which is the number of edges of $G$ between $I$ and $J$ that contain
 points on the support of $D$ in their interiors.  

For each $e\in \ssE_D(I,\ssI^c)$, let $v_e$ be its extremity in
$I$. Denote by $\ssD^I$ the divisor on $\ssGamma_I$ given by
\[
 \ssD^I\coloneqq \ssD_I-\sum_{e\in E_D(I,\ssI^c)}(v_e).
 \]
 Since $D$ is
 $G$-admissible, we get the following
 equalities:
 \begin{equation}\label{DDII}
   \deg(D) = \deg(\ssD_I) + \deg(\ssD_{\ssI^c}) + \delta!_D(I,\ssI^c)=
   \deg(\ssD_I) + \deg(\ssD^{I^c}) + \sse(I,\ssI^c).
   \end{equation}
 
\subsection{Adjoint modular pair associated to a $G$-admissible divisor}\label{sec:adjoint-modular-admissible-divisor} Let $G=(V,E)$ be a finite graph, $\ell\colon E \to (0, +\infty)$ an edge length function, and $\Gamma$ the corresponding metric graph. Let $D$ be a $G$-admissible divisor on $\Gamma$ of degree $d$. We associate to $D$ the pair $(\mu!_D, \mu!_D^*)$ of integral valued functions  on $\ssub{2}!^V$ defined as follows. For each $I \subseteq V$, we set 
\[\mu!_D(I) \coloneqq \deg(\ssD^{I}) + \frac {\sse(I,\ssI^c)}2
  \qquad \textrm{and} \qquad
  \mu!_D^*(I) \coloneqq \deg(\ssD_I) + \frac {\sse(I,\ssI^c)}2.\]

\begin{prop} Let $D$ be a $G$-admissible divisor on $\Gamma$. The pair $(\mu!_D, \mu!_D^*)$ is an adjoint modular pair.
\end{prop}

\begin{proof}
  Let $\ssG_D$ be the spanning subgraph of $G$ whose edge set is $\ssE_D$.
We write
\begin{align*}
\mu!_D^*(I) = \deg(\ssD_I) + \frac {\delta!_D(I,\ssI^c)}2 + \frac {\sse_{D}(I,\ssI^c)}2
\end{align*}
and obtain the decomposition  $\mu!_D^* = \ssq_D + \eta!_D$ for the functions
\begin{align*}
\ssq_D &\colon \ssub{2}!^V \to \Z, \quad \ssq_D(I) \coloneqq \deg(\ssD_I) + \frac {\delta!_D(I,\ssI^c)}2 \quad \textrm{ and }   \\
\eta!_{D} &\colon \ssub{2}!^I \to \Z, \quad  \eta!_D(I) \coloneqq \frac{\sse_{D}(I,\ssI^c)}2 \qquad \forall \, I \subseteq V.
\end{align*} 
By Proposition~\ref{prop:submodularity-boundary} below applied to the graph $\ssG_D$, the function $\eta!_D$ is submodular. Moreover, by Proposition~\ref{prop:centroid-admissible} below, $\ssq_D$ is modular. From these, $\mu!_D^*$ is submodular.

To finish the proof, note that \eqref{DDII} yields
\[\mu!_D(I) + \mu!_D^*(\ssI^c) = \deg(\ssD^I) + \deg(\ssD_{\ssI^c})+
  e(I,\ssI^c)=\deg(D)=d\]
for each $I \subseteq V$. Also, $\mu!_D(\emptyset) =0, \mu!_D^*(V)=d$. So, $(\mu!_D, \mu!_D^*)$ is an adjoint modular pair. 
\end{proof}

\begin{prop}\label{prop:submodularity-boundary} Let $G = (V, E)$ be a finite graph. 
\begin{enumerate}
\item The function $\zeta = \zeta!_G\colon \ssub{2}!^V \to \Z$ defined by
\[\zeta(I) \coloneqq \sse(I,\ssI^c), \qquad \forall \,\, I \subseteq V\]
is submodular.
\item The function $\chi = \sschi_G \colon \ssub{2}!^V \to \Z$ defined by 
\[\chi(I) \coloneqq \sse[I] + \frac{\sse(I,\ssI^c)}2, \qquad \forall \,\, I \subseteq V\]
is modular.
\end{enumerate}
\end{prop}

\begin{proof} Both statements are well-known. For the sake of completeness, we include a short proof. 
Notice that for two  spanning subgraphs $\ssG_1=(V,\ssE_1)$ and $\ssG_2=(V, \ssE_2)$ with disjoint
  edge sets and $\ssE_1\sqcup \ssE_2 = E$, we have $\zeta!_{\ssG_1} +
  \zeta!_{\ssG_2} = \zeta!_G$ and $\sschi_{\ssG_1}+\sschi_{\ssG_2}=\sschi_G$. We
  may thus reduce the proofs of both statements to the case where the
  graph $G$ has only one edge $\{u,v\}$, $u,v\in V$.

  In this case, $\chi(\emptyset)=0$ and $\chi(V)=1$, whereas
  $\chi(\{u\})=\chi(\{v\})=1/2$, whence $\chi$ is modular. Also,
  $\zeta(\emptyset)=\zeta(V)=0$, whereas $\zeta(\{u\})=\zeta(\{v\})=1$.
  Given subsets $I$, $J$ of $V$,
  \[\zeta(I \cup J) +\zeta(I \cap J )=\zeta(I) + \zeta(J)\]
  holds trivially if $I$ or $J$ is $\emptyset$ or $V$, or if $I=J$. In
  the only remaining case, the left-hand side is $0$, whereas the
  right-hand side is $2$, whence $\zeta$ is submodular.
  \end{proof}

\begin{prop}\label{prop:centroid-admissible}
Let $G=(V,E)$, $\ell \colon E \to (0,+\infty)$ an edge length function, and $\Gamma$ the corresponding metric graph. 
Let $D$ be a $G$-admissible divisor of degree $d$ on $\Gamma$. Then, the function $\ssq_D \colon \ssub{2}!^V \to \Z$ defined by
\begin{align}\label{qL}
 \ssq_D(I) &\coloneqq \deg(\ssD_I) + \frac {\delta!_D(I,\ssI^c)}2 \qquad \forall\,\, I\subseteq V
\end{align}
is modular. Furthermore, $\ssq_D$ belongs to $H_d$.
\end{prop}

\begin{proof} Denote by $\ssG_D^c$ the spanning subgraph of $G$ with
  edge set
  $E\setminus\ssE_D$. We write for each subset $I \subseteq V$,
\[\ssq_D(I) = \deg(\ssD_I) + \frac {\delta!_D(I,\ssI^c)}2 = D(I) + \sse_{\ssG_{D}^c}[I] + \frac{\sse_{\ssG_D^c}(I,\ssI^c)}{2} = D(I) + \sschi_{\ssG_{D}^c}(I).\]
Here, $D(I) =\sum_{v\in I} D(v)$, and  $\sschi_{\ssG_{D}^c}$ is the
function defined in Proposition~\ref{prop:submodularity-boundary} for
the graph $\ssG_D^c$. The result now follows by modularity of both functions $I \mapsto D(I)$ and $\sschi_{\ssG_{D}^c}$. The last statement is trivial.
\end{proof}

\subsection{Semistability polytope of a $G$-admissible divisor}\label{sec:semistability-polytope}

Let $G=(V,E)$ be a finite graph, $\ell\colon E \to (0, +\infty)$ an edge length function, and $\Gamma$ the corresponding metric graph.

Let $H =\R^V$. Let $d$ be an integer. A \emph{degree-$d$ polarization} is a modular function on $\ssub{2}!^V$ given by an element $q\in H_d$, that is, an element $q\in \R^V$ with $\sum_{v\in V} q(v)=d$. 
 
 \smallskip
 
 For a polarization $q\in H_d$, a  $G$-admissible divisor $D\in \Div(\Gamma)$ is called
\emph{$q$-semistable} if it has degree $d$ and we have the inequalities 
\begin{equation}
\big|\deg(\ssD_I)-q(I)\big| \leq \frac{\sse(I,\ssI^c)}{2} \qquad \textrm{for each proper nonempty subset } \emptyset \subsetneq I \subsetneq V.
\end{equation}
The $G$-admissible divisor $D$ is called \emph{$q$-stable} if all these inequalities are strict.

\begin{defi}[Semistability polytope] Notation as above, for a $G$-admissible divisor $D$ of degree $d$ on $\Gamma$, the \emph{semistability polytope} of $D$  is the polytope  $\P_{D}$ in $H_d$ defined as
\[\P_{D} \coloneqq \left \{ q \in H_d \,\st \, D \textrm{ is } q\textrm{-semistable}\right\}.\]
In other words, $\P_{D}$ is the set of all $q \in H_d$ that verify
\[\forall \textrm{ proper nonempty subset }\emptyset \subsetneq I \subsetneq V, \qquad \Bigl|\deg(\ssD_I) - q(I)\Bigr|\leq \frac{\sse(I,\ssI^c)}2. \qedhere\]
\end{defi}

\begin{remark}\label{rmk:PD-enough}
    Notice that $\P_D$ does not depend on the location of the points in the supprt of $D$ in the interiors of edges. It depends only on the values of $D$ at vertices of $G$ and on the spanning subgraph $\ssG_D$ of $D$.
\end{remark}

Let $\ssq_D$ be the element of $H_d$ defined in Proposition~\ref{prop:centroid-admissible}, that is, 
\begin{equation*}
\ssq_{D}(I)\coloneqq \deg(\ssD_I)+\frac{\delta!_D(I,\ssI^c)}{2}, \qquad \forall \,\, I \subseteq V. 
\end{equation*}

\begin{thm}\label{thm:semistability-polytope}
Let $D$ be a $G$-admissible divisor of degree $d$ on $\Gamma$.
\begin{enumerate}
\item The semistability polytope  $\P_{D}\subset H_d$ is the set of points $q \in H_d$ that verify the inequalities
\begin{equation*}\label{eq:qqD}
|q(I)-\ssq_{D}(I)|\leq \frac{\sse_{D}(I,\ssI^c)}{2}
\end{equation*}
for each subset $I\subseteq V$.
\item\label{eq:qqD2} Let $(\mu!_D, \mu!_D^*)$ be the adjoint modular pair associated to $D$ in Section~\ref{sec:adjoint-modular-admissible-divisor}.
The semistability polytope $\P_D$ coincides with the base polytope $\P_{\mu!_D}$.
\end{enumerate}
\end{thm}

\begin{proof} We start by proving~\ref{eq:qqD}. Since $q$ and $\ssq_D$ are both modular, and $q(V) =\ssq_D(V)=d$, the set of inequalities in~\ref{eq:qqD} is equivalent to the set of inequalities
\[q(I)- \ssq_D(I)\leq \frac{\sse_{D}(I,\ssI^c)}{2}\qquad \forall \,\,I\subseteq V.\]
Indeed, the inequalities $\ssq_D(I)-q(I)\leq \frac{\sse_{D}(I,\ssI^c)}{2}$ are derived from those above applied to $\ssI^c$.

Using the definition of $\ssq_D$, these inequalities are in turn equivalent to the set of inequalities
\[ q(I) - \deg(\ssD_I)\leq \frac{\sse_{D}(I,\ssI^c) + \delta!_D(I,\ssI^c)}{2} = \frac{\sse(I,\ssI^c)}{2}\qquad \forall \,\,I\subseteq V.\]
To conclude, we note that the other set of inequalities,  $\deg(\ssD_I) - q(I) \leq \frac{\sse(I,\ssI^c)}{2}$ for $I\subseteq V$, follows from the inequalities
\[\ssq_D(I)-q(I)\leq \frac{\sse_{D}(I,\ssI^c)}{2},\]
which imply that 
\[\deg(\ssD_I) - q(I) \leq \frac{\sse_{D}(I,\ssI^c) - \delta!_D(I,\ssI^c)}{2} \leq \frac{\sse(I,\ssI^c)}{2}.\]

We now prove \ref{eq:qqD2}. The base polytope $\P_{\mu!_D}$ is the set
of all $q \in H_d$ which verify $q(I) \leq \mu!_D^*(I)$, $I\subset V$,
as the other set of inequalities, $\mu!_D(I) \leq q(I)$ for $I\subset V$,
is implied. As explained above, the set of inequalities is
equivalent to the ones in~\ref{eq:qqD}, that define the semistability polytope.  
\end{proof}

\subsection{A criterion for simpleness} \label{proof:admissible-divisors}

Let $G=(V,E)$ be a finite graph, $\ell \colon E \to (0, +\infty)$ an edge length function and $\Gamma$ the corresponding metric graph. Let $d\in\Z$.

Let $\cl$ be an element of $\Pic^d(\Gamma)$ and consider $\stable!_G(\cl)$, the set of all the $G$-admissible divisors on $\Gamma$ which are in $\cl$. 
For each ordered partition $\pi=(\sspi_1, \dots, \sspi_s)$ of $V$, we denote by $\ssE_\pi$ the set of all edges with extremities in two different parts of the partition.

\begin{prop}\label{prop:splitG} Let $D$ be a $G$-admissible divisor and $\pi$ be an ordered partition of $V$. Then, $\mu!_{D}$ is $\pi$-split if and only if $\ssG_D$ does not contain any edge of $\ssE_\pi$.
\end{prop}
\begin{proof} This follows immediately from the definition of $\mu!_D$ and Proposition~\ref{muprod}.
\end{proof}
We deduce the following result.

\begin{prop}\label{prop:fulldimconnected} For each $G$-admissible divisor $D$ of degree $d$, the polytope $\P_D$ is full-dimensional in $H_d$ if and only if $\ssG_D$ is connected. 
\end{prop}
\begin{proof} Equivalently, from Proposition~\ref{splitpoly2},  $\mu!_D$ is simple if and only if $\ssG_D$ is connected, which is a consequence of the previous proposition. 
\end{proof}

The $v$-reduced divisors are examples of $G$-admissible divisors whose associated base polytopes are full-dimensional; see Proposition~\ref{prop:v-reduced-adm}

\subsection{Congruence of semistability and Voronoi polytopes}\label{Voronoi} 

Let $G=(V,E)$ be a finite connected graph. In the following, $\E=\E(G)$ is the set of oriented edges of $G$. For each $e\in\E$, let $\te_e$ and $\he_e$ denote its tail and its head. When we write $e=uv$, we mean only that $\te_e=u$ and $\he_e=v$. For each $e\in\E$ and $f\in\R^V$, write
\[
\partial f(e)\coloneqq f(\te_e)-f(\he_e).
\]

The Laplacian 
\[\Delta \colon \R^V \to H_0\]
is the linear map defined by 
\[
\Delta(f)(v) \coloneqq \sum_{\substack{e\in\E\\ \te_e=v}} \partial f(e)\qquad \forall \, v\in V.
\]
Since $G$ is connected, $\Delta$ is surjective. We denote by $\Lambda$ the full rank lattice of $H_0$ given by 
\[
\Lambda\coloneqq \Delta(\Z^V).
\]

For each $h \in H_{0}$, define 
\[\|h\|^2 \coloneqq  \frac{1}{2}\sum_{e\in\E} \partial f(e)^2\] 
for each element $f \in\R^V$ with
$\Delta(f) = h$. It is independent of the choice of $f$ made.

 We denote by $\Vor_{G}$ the Voronoi cell of the origin in $(H_0,\|\cdot\|)$ with respect to the lattice $\Lambda$, that is, 
 \[\Vor_{G} \coloneqq \Bigl\{h\in H_0 \,\, \st\,\, \|h\|^2 \leq \|h-\lambda\|^2 \,\,\,  \forall \lambda \in \Lambda\Bigr\}.\]

\begin{thm}\label{thm:voronoi-semistability} Let $D$ be a $G$-admissible divisor and $\P_D$ the corresponding semistability polytope. Assume $\ssG_D$ is connected. Then, $\P_D$ is a translation of $\Vor_{\ssG_D}$. More precisely, we have 
\[\P_D = \Vor_{\ssG_D} + \ssq_D.\]
\end{thm}

\begin{proof} We work with the spanning subgraph $\ssG_D$ and denote by $\Delta$ the Laplacian of $\ssG_D$ and $\Lambda$ the corresponding lattice. 

Let 
\[\innone{}{\cdot\,, \cdot} \colon H_0 \times H_0 \to \R\]
be the scalar product associated to the quadratic form $\|\cdot\|^2$. It is given by 
\[
\innone{}{\ssq_1, \ssq_2} =\sum_{v\in V} \ssq_1(v) \ssf_2(v) = \sum_{v\in V} \ssf_1(v) \ssq_2(v) = \frac{1}{2}\sum_{e\in\E} \partial f_1(e)\partial f_2(e),
\]
with $\ssf_i\in\R^V$ verifying $\Delta(\ssf_i) =\ssq_i$ for $i=1,2$.

The Voronoi cell is equivalently given as the set of those $h\in H_0$ that verify the set of inequalities 
\begin{align}\label{ineq:voronoi}
 \innone{}{h,\lambda} \leq \frac{\|\lambda\|^2}2\qquad \forall  \lambda\in \Lambda.
\end{align}
The translation $\Vor_{\ssG_D} + \ssq_D$ is the set of all $q \in H_d$ which verify the inequalities 
\[
 \innone{}{q-\ssq_D,\lambda} \leq \frac{\|\lambda\|^2}2\qquad \forall  \lambda\in \Lambda,
\]
equivalently, 
\[
 \bigl|\innone{}{q-\ssq_D,\lambda}\bigr| \leq \frac{\|\lambda\|^2}2\qquad \forall  \lambda\in \Lambda.
\]
Let $\one!_I$ be the characteristic function of the subset $I \subset V$. For $\lambda = \Delta(\one!_I)$, a direct verification shows that $\|\lambda\|^2 = \sse_D(I,\ssI^c)$, and $\innone{}{q-\ssq_D,\lambda} =q(I) - \ssq_D(I)$.  The above inequality applied to $\lambda = \Delta(\one!_I)$ thus gives
\[\bigl|q(I) - \ssq_D(I)\bigr| \leq \frac{\sse_D(I,\ssI^c)}2.\]
Applying Theorem~\ref{thm:semistability-polytope}, we infer that $\Vor_{\ssG_D} + \ssq_D \subseteq \P_D$.

To prove the equality in the above inclusion, it will be enough to show that for $h\in H_0$, the set of inequalities 
\[
 \innone{}{h,\lambda} \leq \frac{\|\lambda\|^2}2
\]
for $\lambda =\Delta(\one!_I)$ for all $I \subseteq V$ implies all the inequalities in \eqref{ineq:voronoi}. 

Let $\lambda\in\Lambda$. Let $f\in C^0(G,\Z)$ such that $\lambda=\Delta(f)$. 
Denote by $\ssn_1<\ssn_2< \dots< \ssn_r$ all 
the values taken by $f$. Since
$\Delta(f)=\Delta(f+\one!_{V})$, we may assume
that $\ssn_1\geq 0$.   For each $i=1,\dots,r$, define the subset $\ssI_i\subseteq V$ by
\[\ssI_i \coloneqq \bigl\{v\in V\,\, \st\,\, f(v) \geq \ssn_i \bigr\}.\]
We have 
\[\ssI_1 \supset \ssI_2 \supset \dots \supset \ssI_r.\]
Setting $\ssn_{0}\coloneqq 0$, we can write
$f = \sum_{j=1}^r (\ssn_j -\ssn_{j-1}) \one!_{I_j}$. This gives
\[
\lambda=\Delta(f) = \sum_{j=1}^r (\ssn_j -\ssn_{j-1}) \Delta(\one!_{I_j}).
\]
Note that we have 
\[
\innone{}{\Delta(\one!_{I_j}), \Delta(\one!_{I_k})} = \Bigl|\,\ssE_D(\ssI_j,\ssI_j^c) \cap \ssE_D(\ssI_k,\ssI_k^c)\,\Bigr|\geq 0
\]
for each
$j,k\in [r]$. This implies that 
\[\|\lambda\|^2 \geq \sum_{j=1}^r (\ssn_j -\ssn_{j-1})^2 \|\Delta(\one!_{I_j})\|^2 \geq \sum_{j=1}^r (\ssn_j -\ssn_{j-1}) \|\Delta(\one!_{I_j})\|^2.\] 
We infer that the inequality $\innone{}{h,\lambda} \leq \frac{\|\lambda\|^2}2$ for $\lambda=\Delta(f)$ is implied by the $r$ inequalities 
\[
\innone{}{h,\Delta(\one!_{I_j})} \leq \frac{\|\Delta(\one!_{I_j})\|^2}2 \qquad j=1, \dots, r,
\]
and the theorem follows.
\end{proof}

\begin{remark}\label{rmk:voronoi-semistability}
    If $\ssG_D$ is not connected, then $\P_D$ is not a translate of a Voronoi cell in $H_0$, as $\P_D$ does not have maximum dimension by Proposition~\ref{prop:fulldimconnected}. On the other hand, let $\mu\coloneqq\mu!_D$, and let $\pi=(\sspi_1,\dots,\sspi_s)$ be an ordered $s$-partition such that $\mu=\mu!_{\pi}$, for $s\coloneqq\cd!_{\mu}$. By Proposition~\ref{muprod}, 
    $\P_{\mu}$ decomposes as the product of the $\P_{\mu\rest{\pi_i}}$ inside ${\mathbb R}^{\sspi_i}$. Now, 
    each $\mu\rest{\pi_i}$ is a simple supermodular function on $\pi_i$ because $s=\cd!_{\mu}$, whence $\P_{\mu\rest{\pi_i}}$ has maximum dimension by Proposition~\ref{splitpoly2}. Also,   $\mu\rest{\pi_i}=\mu!_{D_{\pi_i}}$, and therefore, $\ssG_{D_{\pi_i}}$ is a connected spanning subgraph of $G[\sspi_i]$ for each $i$ by Proposition~\ref{prop:fulldimconnected}. Applying Theorem~\ref{thm:voronoi-semistability} to each $\ssD_{\pi_i}$, we conclude that $\P_D$ is the translation by $\ssq_D$ of the product of the Voronoi cells of the $\ssG_{D_{\pi_i}}$.
\end{remark}


\section{Characterization of admissible divisor classes}\label{sec:admissible-characterization}
In this section we prove Theorem~\ref{thm:admissible2}, which gives a characterization of all $G$-admissible divisors that lie in a given linear equivalence of divisors on a metric graph. 

Let $G=(V,E)$ be a finite graph. Recall that $\E=\E(G)$ is the set of oriented edges of $G$, that for each $e\in\E$, we let $\te_e$ and $\he_e$ denote its tail and its head, and that we write $e=uv$ to simply mean that $\te_e=u$ and $\he_e=v$. In $\Gamma$, we may identify each $e\in\E$ with the interval $[0, \ell_e]$, where $0$ corresponds to $\te_e$ and $\ell_e$ to $\he_e$. For each $r\in [0, \ell_e]$, we denote by $x^e_r$ the point at distance $r$ from $u$ on the edge $e$. Note that with this notation we have $x_r^e = x_{\ell_e-r}^{\bar e}$, where $\bar e$ is $e$ with the opposite direction. We call the set of points $x^e_r$ with $0<r<\ell_e$ the interior of $e$, denoted $\ring e$. Finally, for each $I,J\subseteq V$, let $\E(I,J)$ be the subset of $e\in\E$ with $\te_e\in I$ and $\he_e\in J$.

\subsection{Linear equivalence of divisors}  Let $\ell\colon E \to (0, +\infty)$ be an edge length function, and $\Gamma$ the corresponding metric graph. A rational function on $\Gamma$ is a continuous function $f\colon \Gamma \to \R$ whose restriction to each edge of $\Gamma$ is piecewise affine  with integral slopes.  Denote by $\Rat(\Gamma)$ the set of rational functions on $\Gamma$.

The order of vanishing, $\ord_x(f)$, of a rational function $f$ on $\Gamma$ at $x\in\Gamma$ is the sum of the incoming slopes of $f$ at $x$. The divisor $\div(f)$ associated to $f\in \Rat(\Gamma)$ is then defined by
 \[\div(f)\coloneqq \sum_{x\in \Gamma} \ord_x(f) (x).\]
   Elements of this form belong to $\Div^0(\Gamma)$ and are called principal. They form a subgroup of $\Div^0(\Gamma)$ that we  denote by $\Prin(\Gamma)$. A divisor $D_1$ is called linearly equivalent to a divisor $D_2$, and we write $D_1 \sim D_2$, if the difference $D_1 -D_2$ is principal.

 Let $\Pic(\Gamma)$ be group of divisors in $\Gamma$ modulo linear equivalence, that is,
 \[\Pic(\Gamma)\coloneqq \rquot{\Div(\Gamma)}{\sim}.\]
 For each $d\in\Z$, let $\Pic^d(\Gamma)\subseteq \Pic(\Gamma)$ be the collection of classes of divisors in $\Div^d(X)$. Then, $\Pic^0(\Gamma)\subseteq \Pic(\Gamma)$ is a subgroup and the $\Pic^d(\Gamma)$ are its cosets.

In the following, we use calligraphic letters such as $\cl$ to denote elements of $\Pic(\Gamma)$. We write $D \in \cl$ to say that $D$ is a divisor on $\Gamma$ whose class in $\Pic(\Gamma)$ is $\cl$. Let $\stable!_G(\cl)$ denote the set of all $G$-admissible divisors $D \in \cl$.

\subsection{Domination property of $G$-admissible divisors}  Before we proceed, we formulate a key property of $G$-admissible divisors. This lemma will be of fundamental importance both in this section and later.

\begin{lemma}[Domination property]\label{lem:varH} Let $D$ be a divisor on $\Gamma$ and let $h\in \Rat(\Gamma)$ so that $\div(h)+D$ is $G$-admissible. Let $h'\in \Rat(\Gamma)$ be a rational function with the property that $\div(h') +D$ is effective outside the set of vertices of $G$. 
The following are equivalent:
\begin{enumerate}[label=(\arabic*)]
\item\label{ineq:vert} $h'(v) \geq h(v)$ for all $v\in V$.
\item\label{ineq:all} $h'(x)\geq h(x)$ for all $x \in \Gamma$.
\end{enumerate} Moreover, the inequalities in \ref{ineq:vert} are all strict if and only if so are the inequalities in \ref{ineq:all}.
\end{lemma}
\begin{proof} Implication \ref{ineq:all} $\Rightarrow$ \ref{ineq:vert} is obvious. As for the converse, put $D'\coloneqq \div(h) + D$. By hypothesis, $D'$ is $G$-admissible. Also, $\div(h'-h) + D'$ is effective outside the set of vertices.  We may thus assume $h=0$ and $D$ is $G$-admissible. For the sake of a contradiction, assume $h'(x)<0$ for a point $x \in \Gamma$. By assumption, $x\notin V$. Let $e=uv$ be the edge of $G$ that contains $x$ in its interior.  Let $S$ be the closed subset of $e$ where $h'\rest{e}$ takes its minimum.  Since $h'(u),h'(v) \geq 0$, and $h'(x)<0$, we get $S \subset \ring e$. For each point $y$ on the boundary $\partial S$ of $S$ in $\ring e$, we get $\div(h')(y) \leq -1$. Moreover, if $S=\{y\}$ then $\div(h')(y) \leq -2$. Since $D$ is $G$-admissible, its support contains at most one point on $\ring e$ and the coefficient of $D$ at this point is one. We infer that $D+\div(h')$ has a negative coefficient in $\ring e$, a contradiction.   

The same proof applies to the last statement in the lemma. 
\end{proof}

\subsection{Twister function $\twistd{D}$ associated to a divisor} To each divisor $D\in \Div(\Gamma)$, we associate 
its \emph{twister function} 
\[\twistd{D}\colon \E \rightarrow \R\] 
which on any oriented edge  $e$ of $G$ takes the value 
\[\twistd{D}(e) \coloneqq\sum_{t\in (0, \ell_e)} (\ell_e-t) D(x^e_t).\]
Note that the above sum concerns only finitely many points, those in the support of $D$ which lie in the interior of $e$. This means that $\twistd{D}(e)$ is well-defined.

\begin{prop}\label{tde}
 For each oriented edge $e$ of $G$, we have
 \[
 \twistd{D}(e) + \twistd{D}(\ol e) = \ell_e \sum_{t\in (0, \ell_e)} D(x_t^e).
 \]
\end{prop}

\begin{proof}
This is clear from the identity $x_t^e = x_{\ell_e-t}^{\ol e}$.
\end{proof}

\subsection{Twisted integer slopes} Let $\twist\colon \E \rightarrow \Z$
be an integer valued function. 

\begin{defi}\label{twdr}\rm  For each real valued function $f\colon V \rightarrow \R$
and $e\in\E$, define 
\[
\slztwist{\ell}{\twist}f(e)\coloneqq \Big\lfloor \frac{f(\te_e)-f(\he_e)+\twist(e)}{\ell_e}\Big\rfloor,
\]
that we call the \emph{(incoming) twisted integer slope of $f$ at} $\te_e$ \emph{along the edge $e$}.
\end{defi}
The name given to $\slztwist{\ell}{\twist}f(e)$ will be justified in Theorem~\ref{thm:existence-uniqueness-canonical-extension}. Note that if the edge lengths are all equal to one, and the twist function is 0, then we have $\slztwist{\one}{0}f(e)=\partial f(e)$, so the notation is consistent with the one 
in Section~\ref{sec:admissible-divisor-semistability-polytope}. 

The following statement follows directly  from Proposition~\ref{tde}.

\begin{prop} For a divisor $D$ on $\Gamma$ 
and $\tau=\twistd{D}$, we have the equality for each $e\in\E$:
\[\slztwist{\ell}{\twist}f(e) + \slztwist{\ell}{\twist}f(\ol e) = \epsilon + \sum_{t\in (0, \ell_e)} D(x_t^e),\]
 with $\epsilon$ equal to $0$ or $-1$ depending on whether $\frac{f(\te_e)-f(\he_e)+ \twist(e)}{\ell_{e}}$ is an integer or not, respectively.
\end{prop}

\subsection{Admissible extensions of functions with respect to a divisor} Let $D$ be a divisor on $\Gamma$.
 For each real valued function $f\colon V \rightarrow \R$, we seek a rational extension $\hat f \colon \Gamma \to \R$ for which $D+\div(\hat f)$ is $G$-admissible. The following theorem guarantees the existence and uniqueness of such an extension.

We will abuse the notion of Euclidean division by saying that $y-\lfloor y/x\rfloor x$ is the remainder of the Euclidean division of any real number $y$ by any positive real number $x$.

\begin{thm}[Admissible extension of a function with respect to a divisor]\label{thm:existence-uniqueness-canonical-extension}
Notation as above, there is a unique rational extension 
$\hat f\colon \Gamma \rightarrow \R$ such that 
$D+\div(\hat f)$ is $G$-admissible. Moreover, letting $\tau = \twistd{D}$, the twister function associated to $D$, the function $\hat f$ 
is  characterized by the following properties:
\begin{enumerate}
\item[(1)] The incoming slope of $\hat f$ at each $u\in V$ along an oriented edge $e$ with $\te_e=u$ is $\slztwist{\ell}{\twist}f(e)$. 
\item[(2)] For each oriented edge $e$
 of $G$ and $t\in (0, \ell_e)$, the divisor $D+\div(\hat
 f)$ takes value $0$ at $x_t^e$, unless $\ell_e-t$ is the remainder of 
the Euclidean division of $f(\te_e) -f(\he_e)+ \tau(e)$ by $\ell_e$, in which case it takes
value $1$.
\end{enumerate}
\end{thm}

\begin{proof} We first prove the uniqueness. Let $\hat f, \hat f' \colon \Gamma \to \R$ be two rational functions extending $f$ such that both $D+ \div(\hat f)$ and $D+ \div(\hat f')$ are $G$-admissible. We apply Lemma~\ref{lem:varH} to the divisor $D$ and rational functions $\hat f,\hat f'$, and deduce that $\hat f' \geq\hat f$ and $\hat f \geq \hat f'$, hence $\hat f = \hat f'$, as required. 

We now prove the existence. Identify each oriented edge $e=uv$ with the interval $[0,\ell_e]$, and enumerate the points in $\mathrm{Supp}(D)$ lying in $\ring e$ by $0<t_1< t_2<\dots<t_m<\ell_e$. Also, let $r_e \in (0,\ell_e]$ be so that $\ell_e-r_e$ is the remainder of the Euclidean division of $f(u) -f(v)+ \tau(e)$ by $\ell_e$, so that
\[f(u) - f(v) +\tau(e) =  \slztwist{\ell}{\twist}f(e)\ell_e + (\ell_e-r_e).\] 

Consider the function $g_e\colon [0,\ell_e] \to \R$ defined by 
\[
g_e(t)\coloneqq f(u)-\slztwist{\ell}{\twist}f(e)t+ \left (\sum_{i=1}^m D(x^e_{t_i})\max(0,t-t_i)\right)-
\max(0,t-r_e) \qquad \forall\,\, t\in[0, \ell_e].
\]
 Note that we have $g_e(0) = f(u)$ by definition. Also, $g_e(\ell_e) = f(v)$, since
 \begin{align*}
g_e(\ell_e)&=f(u)-\slztwist{\ell}{\twist}f(e)\ell_e+ \left (\sum_{i=1}^m D(x^e_{t_i})(\ell_e-t_i)\right)- (\ell_e-r_e) \\
&=f(u) - f(v) +\tau(e)- \slztwist{\ell}{\twist}f(e)\ell_e -(\ell_e-r_e) +f(v) = f(v).
\end{align*}
This means that there exists a rational function $\hat f\colon \Gamma\to\R$ which on each oriented edge $e$ of $G$ coincides with $g_e$. By the definition of the functions $g_e$, Properties (1) and (2) in the statement of the theorem are satisfied. This means that $D+\div(\hat f)$ is $G$-admissible.  \qedhere
\end{proof}

\begin{defi}\rm Let $D$ be a divisor on $\Gamma$ and $f\colon V \to\R$ a real-valued function on the vertices of $G$. The function $\hat f$ given by the above theorem will be called the
\emph{$G$-admissible extension} of $f$ with respect to $D$. The principal divisor 
$\div(\hat f)$ will be denoted $\div_\ell(f; D)$.
\end{defi}

\begin{prop}\label{prop:adm}
Let $D\in\Div(\Gamma)$ and $f\in\R^V$. Set $\twist\coloneqq\twistd{D}$ and $D'\coloneqq D+\div_\ell(f;D)$. Let $z\in\Gamma$. Then:
\begin{itemize}
 \item[(1)] If $z$ is a vertex of $G$, then 
\[D'(z) = D(z) + \sum_{\substack{e\in\E\, \\ \text{\rm t}_e=z}}\,\slztwist{\ell}{\tau}f(e).\] 
 \item[(2)] If $z=x_r^{\ol e}$, where $e=uv$ is an oriented edge, and $r$ is the remainder of the Euclidean division of $f(u)-f(v)+\tau(e)$ by $\ell_e$, such that $r$ is nonzero, then $D'(x_{r}^{\ol e}) =1$.
 \item[(3)] In any other case,  we have $D'(z)=0$.
\end{itemize}
\end{prop}

\begin{proof} This is a reformulation of Properties (1)
  and (2) in Theorem~\ref{thm:existence-uniqueness-canonical-extension}.
\end{proof}

\subsection{Characterization of $G$-admissible divisors}
The following is our theorem on characterization of $G$-admissible divisors in $\stable!_G(\cl)$, for a class $\cl \in \Pic(\Gamma)$.

\begin{thm}\label{thm:admissible2} Let $D$ be a divisor on $\Gamma$. 
\begin{itemize}
 \item[$(i)$] For each $f\colon V \to\R$, 
the divisor $D+\div_\ell(f;D)$ is $G$-admissible.
 \item[$(ii)$] Each $G$-admissible divisor $D'$ on $\Gamma$ that is linearly 
equivalent to $D$ 
 is of the form $D+\div_\ell(f;D)$ for some $f\colon V\to\R$.
\end{itemize}
\end{thm}

\begin{proof}
 By Proposition~\ref{prop:adm}, each divisor of the form
 $D+\div_\ell(f;D)$ is $G$-admissible, whence the first statement. 
Conversely, to prove the second statement, 
write $D' = D +\div(h)$ for a rational function $h\colon \Gamma \to\R$, and denote
by 
$f\colon V\to\R$ the restriction of $h$ to the vertices of
$G$. Then, $h$ is an extension of $f$ such that $D+\div(h)$ is
$G$-admissible. It follows from Theorem~\ref{thm:existence-uniqueness-canonical-extension}
that $h$ is the $G$-admissible
extension $\hat f$ of $f$ with respect to $D$, and thus $D'=D+\div_\ell(f;D)$.
\end{proof}

\subsection{A consequence} We state the following consequence of
Theorem~\ref{thm:existence-uniqueness-canonical-extension}.

\begin{prop}\label{cor:sumcan} Let $D$ be a divisor on $\Gamma$ and $\ssf_1,
  \ssf_2\colon V\to\R$. Let $f \coloneqq \ssf_1+\ssf_2$.
Denote by $\hat f$ and $\ssfhat_1$ the $G$-admissible
extensions of $f$ and $\ssf_1$ with respect to $D$.
Denote by $\ssfhat_2$ the $G$-admissible
extension of $\ssf_2$ with respect to the divisor $\ssD_1\coloneqq D +\div_\ell(\ssf_1; D)$.
Then,
\[
\ssfhat_1+\ssfhat_2=\hat f.
\]
In particular, 
 \[
\div_\ell(f;D) = \div_\ell(\ssf_1; D)+\div_\ell(\ssf_2; \ssD_1).
 \]
\end{prop}

\begin{proof} 
Note that
\begin{align*}
D+ \div(\ssfhat_1+ \ssfhat_2) =
D+ \div(\ssfhat_1)+ \div(\ssfhat_2) =
\ssD_1+ \div(\ssfhat_2).
\end{align*}
It follows that $D+ \div(\ssfhat_1+ \ssfhat_2)$ is
$G$-admissible. 
Since $\ssfhat_1+ \ssfhat_2$ restricts to
$\ssf_1+\ssf_2$ on the vertices of $G$, it follows from
Theorem~\ref{thm:existence-uniqueness-canonical-extension} that $\ssfhat_1+
\ssfhat_2=\hat f$. The last statement is immediate.
\end{proof}

\subsection{$\valgroup$-Rationality}\label{rationality} 
Let $\valgroup\subseteq\R$ be an additive subgroup containing $\ell_e$ for each $e\in E$. We say that a divisor $D$ on $\Gamma$ is \emph{$\valgroup$-rational} if for each edge $e=uv\in\E$, each point $x$ of the support of $D$ on $e$ is of the form $x_t^e$ for $t\in\valgroup$. 
The $\Lambda$-rational divisors form a subgroup of $\Div(\Gamma)$. A rational function $h\colon\Gamma\to\R$ is called \emph{$\valgroup$-rational} if $\div(h)$ is $\valgroup$-rational.

For each $\cl\in\Pic(\Gamma)$, we denote by  $\cl_{\valgroup}\subseteq\cl$ the subset of $\valgroup$-rational divisors, and by $\stable!_G(\cl_{\valgroup})\subseteq\stable!_G(\cl)$ the subset of $G$-admissible, $\valgroup$-rational divisors.

\begin{prop}\label{prop:rational-admissible extension} Let $D$ be a $\valgroup$-rational divisor on $\Gamma$. Let $f\colon V\to\R$ be a real-valued function and $\hat f$ the $G$-admissible extension of $f$ with respect to $D$. Then,
$\hat f$ is $\valgroup$-rational if and only if $f(u)-f(v)\in\valgroup$ for each $e=uv\in\E$.
\end{prop}

\begin{proof} Put $\twist\coloneqq\twistd{D}$. Since $D$ is $\valgroup$-rational, $\twist(e)\in\valgroup$ for every $e\in\E$. Then, $f(u)-f(v)\in\valgroup$ for each $e=uv\in\E$ if and only if the remainder of the Euclidean division of $f(u)-f(v)+\twist(e)$ by $\ell_e$ is in $\valgroup$ for each $e=uv\in\E$, if and only if $D+\div(\hat f)$ is $\valgroup$-rational, by Theorem~\ref{thm:existence-uniqueness-canonical-extension}, if and only if $\hat f$ is  $\valgroup$-rational.
\end{proof}

\subsection{Reduced divisors, I}\label{sec:v-reduced-1} Recall that a divisor $D\in\Div(\Gamma)$ is $v$-reduced for a vertex $v\in V$ if $\ssD(x)\geq 0$ for all $x\neq v$ in $\Gamma$, and there is no closed subset $S \subset \Gamma \setminus\{v\}$ that can chip-fire and respects the effectiveness of $D$ outside $v$.  In each linear equivalence class $\cl\in\Pic(\Gamma)$, there is a unique $v$-reduced divisor for each $v\in V$.

\begin{prop}\label{prop:v-reduced-adm} Let $v\in V$ and $D\in\Div(\Gamma)$ be a $v$-reduced divisor of degree $d$. Then, $D$ is $G$-admissible, $\ssG_D$ is connected and $\P_D$ is full-dimensional in $H_d$. Furthermore, if $D$ is linearly equivalent to a $\valgroup$-rational divisor, then $D$ is $\valgroup$-rational.
\end{prop}

\begin{proof} Indeed, the chip-firing property of $D$ yields that $\ssD$ has at most degree one in the interior of each edge of $G$, or equivalently, that $\ssD$ is $G$-admissible. 

Furthermore, the subgraph $\ssG_D$ is connected as otherwise, we can find a closed subset of $\Gamma \setminus \{v\}$ that can chip-fire for $D$ keeping effectivity outside $v$. This implies that the polytope $\P_{D}$ associated to $D$ is full-dimensional in $H_d$.

Let $\cl$ be the class of $D$ and $\ssD_0\in\cl$ a $\valgroup$-rational divisor. Let $h$ be an element of $\Rat(\Gamma)$ such that $\ssD=\ssD_0+ \div(h)$. Now, $\ssD_0$ is $\valgroup$-rational, and the edge lengths are all in $\valgroup$, whence $h(u)-h(v)\in\valgroup$ for each $e=uv$ in $\ssG_D$. Since $\ssG_D$ is connected, this means that $h(u)-h(v)\in\valgroup$ for each $u,v\in V$. Since $D$ is $G$-admissible, Proposition~\ref{prop:rational-admissible extension} yields that $h$ is $\Lambda$-rational, and thus $D$ is $\valgroup$-rational.
\end{proof}


\section{Tilings by semistability polytopes of admissible divisors} \label{sec:tiling-semistability}

Let $G=(V,E)$ be a graph, $\ell \colon E \to (0, +\infty)$ an edge length function and $\Gamma$ the corresponding metric graph.  Let $\valgroup\subseteq\R$ be an additive subgroup containing $\ell_e$ for each $e\in E$. Recall that for each class $\cl$ in $\Pic^d(\Gamma)$, we denote by $\stable!_G(\cl_{\valgroup})$ the set of all $G$-admissible, $\valgroup$-rational divisors $D \in \cl$. Let $\C!_G(\cl_{\valgroup}) \subset \M_d(V)$ be the collection of all supermodular functions $\mu!_D$ for $D\in \stable!_G(\cl_{\valgroup})$.
Define 
\[\P_{\cl} \coloneqq \P_{\C!_G(\cl_{\valgroup})} = \bigcup_{D\in \stable!_G(\cl_{\valgroup})} \P_D.\]
Here is the main theorem of this section.

\begin{thm}\label{thm:semistability-tiling-sccp} If $G$ is connected, $\cl_{\valgroup}$ is nonempty and $\valgroup$ is dense in $\R$, then the collection $\C!_G(\cl_{\Lambda})$ is separated, closed, complete and positive in the sense of Section~\ref{sec:admissible-families-tilings}.
\end{thm}

As a consequence, we deduce the following theorem.

\begin{thm}[Tilings by semistability polytopes of admissible divisors] \label{thm:semistability-tiling} Let $\cl$ be a linear equivalence class of divisors of degree $d$ on a metric graph $\Gamma$ with model $(G,\ell)$ and $G$ connected. For a subgroup $\valgroup\subseteq\R$, if $\cl_{\valgroup}$ is nonempty and $\valgroup$ is dense in $\R$, then the set of semistability polytopes $\P_D$ for $D\in\stable!_G(\cl_{\valgroup})$ gives a tiling of $H_d$, that does not depend on the choice of $\valgroup$.
\end{thm}

\begin{figure}[t]
  \centering
      \scalebox{.3}{\input{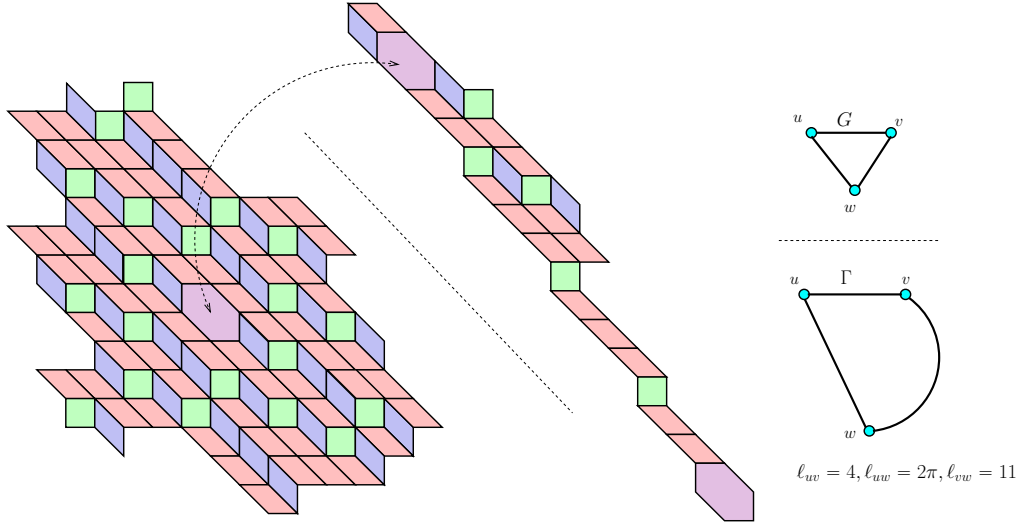}}
  \caption{Tiling arising from admissible divisors in the zero class on a metric graph $\Gamma$. The graph $G$ is a $3$-cycle and the edge lengths are depicted in the picture, on the right.  Note that the ratios of edge lengths are not all rational. The middle figure shows periodicity of the tiling in one direction. As in Figure~\ref{fig:tiling3}, we depict a projection of the tiling to $\R^{\{u,v\}}$. 
  }
  \label{fig:tiling2}
  \end{figure}

Figures \ref{fig:tiling3},  \ref{fig:tiling2}, and \ref{fig:tiling4} are three examples of such tilings of $H_0\subseteq\R^3$ for the triangle graph $G$. Notice that, while each $\P_D$ does not depend on the location of points of the support of $D$ in the interiors of the edges, and in particular, on the lengths of the edges, in general, the tiling given by all of the $\P_D$ for $D\in\stable!_G(\cl_{\valgroup})$ does. The aforementioned figures hint at that, as the class $\cl$ is that of zero in all examples, but the tilings are very different. The first and the last are periodic but the second is periodic only in one direction. We will discuss periodicity in Section~\ref{sec:admissible-sum}.

In the rest of this section we will prove the properties of 
$\C!_G(\cl_{\Lambda})$ asserted in Theorem~\ref{thm:semistability-tiling-sccp}, assuming only what is necessary for each property, as stated in each proposition below, and then combine the results to prove the above two theorems. 

\subsection{$\C!_G(\cl_{\valgroup})$ is separated}\label{sec:separated-admissible} 

\begin{prop}\label{prop:Cseparated} The collection $\C!_G(\cl_{\valgroup})$ is separated.
\end{prop} 
\begin{proof} Let $\ssD_1$ and $\ssD_2$ be distinct $\valgroup$-rational $G$-admissible divisors in $\cl$, and put $\mu!_1\coloneqq\mu!_{D_1}$ and $\mu!_2\coloneqq\mu!_{D_2}$. Then, $\ssD_2=\ssD_1+\div(\hat f)$ for a function $f \colon V \to \R$, with $\hat f$ the $G$-admissible extension of $f$ with respect to $\ssD_1$. Let $\epsilon!_1,\dots,\epsilon!_r$ be the increasing sequence of values of $f$, taken on subsets $\sspi_1,\dots,\sspi_r$ of $V$, respectively. By Proposition~\ref{prop:rational-admissible extension}, we may assume the $\epsilon_j$ are in $\Lambda$. Since $\ssD_1\neq\ssD_2$, we have $r>1$. 

Let $I\coloneqq \sspi_1$ and $J\coloneqq \ssI^c$. Put $\varpi\coloneqq (I,J)$. We claim that $\varpi$ is a separation for $\mu!_1$ and $\mu!_2$. We need to show that 
\[\mu!_1(I) + \mu!_2(J) \geq d,\]
that is,
\[\deg(\ssD_{1,I}) -\sse_{D_1}(I,J) + \frac{\sse(I,J)}{2} + \deg(\ssD_{2,J}) -\sse_{D_2}(J,I) + \frac{\sse(J,I)}{2} \geq d.\]
This is equivalent to 
\[
\deg(\ssD_{1,I}) + \delta!_{D_1}(I,J) + \deg(\ssD_{2,J}) -\sse_{D_2}(I,J) \geq d,
\]
whence to 
\[\deg(\ssD_{2,J}) -\deg(\ssD_{1,J}) - \sse_{D_2}(I,J) \geq 0.\]
Using Proposition~\ref{prop:adm}, the left-hand side in the above sum is
\[
\sum_{e\in\E(J,I)}\big(-\slztwist{\ell}{\tau} f(e)-\epsilon!_{D_2}(e)\big),
\]
where $\tau\coloneqq \twistd{\ssD_1}$, and  
$\epsilon!_{D_2}(e)\coloneqq 1$ if $e\in\E(\ssG_{D_2})(I,J)$, else $\epsilon!_{D_2}(e)\coloneqq 0$.

\begin{figure}[t]
  \centering
      \scalebox{.3}{\input{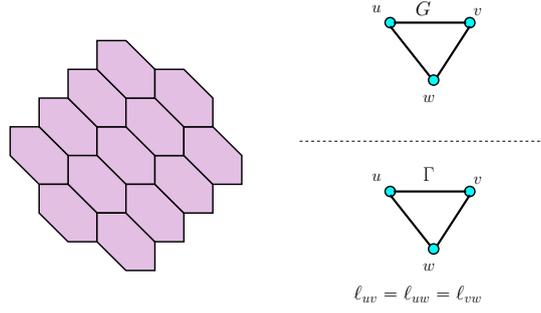}}
  \caption{Tiling arising from admissible divisors in the zero class on a metric graph $\Gamma$. The graph $G$ is a $3$-cycle and the edge lengths are all equal. As in Figure~\ref{fig:tiling3}, we depict a projection of the tiling to $\R^{\{u,v\}}$.
  }
  \label{fig:tiling4}
  \end{figure}

Consider $e \in\E(I, J)$. Notice that $\slztwist{\ell}{\tau} f(e)\leq 0$ from our choice of $\varpi$. Thus, the contribution of $e$ to the above sum is non-negative unless $e\in\E(\ssG_{D_2})(I, J)$. Suppose $e\in\E(\ssG_{D_2})(I,J)$. By Proposition~\ref{prop:adm} again, $f(\te_e)-f(\he_e)+\tau(e)$ is divisible by $\ell_e$. Since $f(\te_e)<f(\he_e)$ from our choice of $\varpi$ and 
$\tau(e)\geq 0$, we must have $\slztwist{\ell}{\tau} f(e)\geq 1$. Therefore, also in this case, the contribution of $e$ to the above sum is non-negative, finishing the proof of our claim.

If the separation is trivial, then $\mu!_1$ and $\mu!_2$ are $\varpi$-split, and the inequalities above are equalities. The first assertion yields $\ssE_{D_1}(I,J)=\ssE_{D_2}(I,J)=\emptyset$. This and the second assertion yield $\slztwist{\ell}{\tau} f(e)=0$, or equivalently, $f(\te_e)-f(\he_e)+\tau(e)<\ell_e$ for each $e\in\E(J,I)$. Now, put $\ssf_1\coloneqq \epsilon!_2\one!_{J}+\epsilon!_1\one!_{I}$ and 
$\ssD'\coloneqq \ssD_1+\div(\ssfhat_1)$, with $\ssfhat_1$ the $G$-admissible extension of $\ssf_1$ with respect to $\ssD_1$. Then, $\ssD'\in\cl_{\valgroup}$ by Proposition~\ref{prop:rational-admissible extension}. Since $\epsilon!_2-\epsilon!_1+\tau(e)<\ell_e$ for each $e\in\E(J,I)$, 
by Proposition~\ref{prop:adm} again, $\mu!_{D'}=\mu!_1$. On the other hand, Proposition~\ref{cor:sumcan} yields $\ssD_2=\ssD'+\div(\ssfhat_2)$, with $\ssfhat_2$ the $G$-admissible extension of $\ssf_2$ with respect to $\ssD'$, for $\ssf_2\coloneqq f-\ssf_1$. 

Notice that $\ssf_2$ takes only $r-1$ values. If $r=2$ then $\ssD_2=\ssD'$, hence $\mu!_2=\mu!_{D'}=\mu!_1$. Thus, if $\mu!_1\neq\mu!_2$, either $r>2$, or $r=2$ and the separation $\varpi$ is nontrivial. Argue now by induction on $r$. If $\mu!_1\neq\mu!_2$ and $\varpi$ is trivial, then $\mu!_2\neq\mu!_{D'}$, and thus by induction there is a nontrivial separation for 
$\mu!_{D'}$ and $\mu!_2$. Since $\mu!_{D'}=\mu!_1$, we conclude.
\end{proof}

\subsection{Admissible chip-firing} \label{sec:admissible-chip-firing}
For each $G$-admissible divisor $D$ and each subset $I \subset V$, let $X\coloneqq \ssX_D(I)$ be the closed subset of $\Gamma$ that consists of the union of all the edges that join a pair of vertices of $I$, and all the edge segments that join a vertex of $I$ to a point of the support of $D$ on the interior of an edge $e\in E(I,\ssI^c)-\ssE_D(I,\ssI^c)$. On each edge of $E(I,\ssI^c)$, we thus get an entire segment consisting of the points that are not included in $X$. Assume $E(I,\ssI^c)\neq\emptyset$. 

\begin{defi}\label{defi:minimum-length} Notation as above, we denote by $\ml_D(I)$ the \emph{minimum length} of the segments of the form $e \setminus \ssX_D(I)$ for $e \in E(I,\ssI^c)$.
\end{defi}

\begin{figure}[!t]
\centering
    \scalebox{.3}{\input{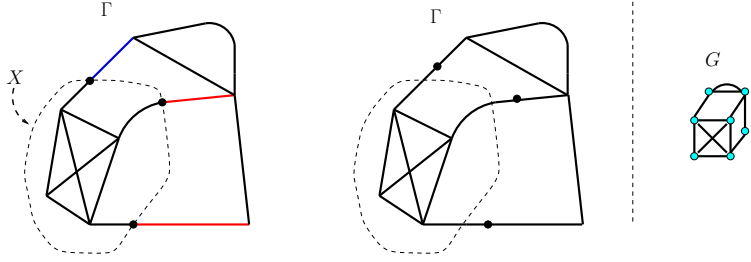}}
\caption{The set $X$ in the metric graph $\Gamma$ depicted in the figure is a cut set. The three dotted points are in the support of a $G$-admissible divisor. When $X$ fires at distance $\epsilon>0$, the chip placed on each dotted point moves to the point at distance $\epsilon$ on each of the three outgoing branches from $X$. Here the move is small.}
\label{fig:admissible-chip-firing2}
\end{figure}

In Figure~\ref{fig:admissible-chip-firing2}, $I$ is a subset of four vertices which form a complete graph in $G$, and $X=\ssX_D(I)$ is depicted on the left. Here, $\ssE_D(I, \ssI^c)=\emptyset$, and the segments $e \setminus \ssX_D(I)$, $e\in E(I, \ssI^c)$, are colored in red or blue. The length of the blue segment gives $\ml_D(I)$.

\begin{prop}\label{prop:minimum-length} We have $\ml_D(I)>0$. If $D$ is $\valgroup$-rational, then $\ml_D(I)$ belongs to $\valgroup$. 
\end{prop}
\begin{proof} This is straightforward.
\end{proof}

Let $\epsilon \in [0, \ml_D(I)]$ be a real number and $f\coloneqq \ssf_{\epsilon}^I\colon V\to\R$ be the function taking value $0$ on $I$ and value $\epsilon$ on $\ssI^c$. In this case, the divisor $D +\div(f; D)$ coincides with the chip-firing move of $D$ by distance $\epsilon$ from the cut $X$ of $\Gamma$.  We call this divisor the \emph{$G$-admissible chip-firing move of $D$ by distance $\epsilon$ induced by $I$.} If $\epsilon$ is smaller than $\ml_D(I)$, we say the move is \emph{small}; otherwise, when $\epsilon =\ml_D(I)$, we call it \emph{flipping}. We refer to Figure~\ref{fig:admissible-chip-firing2} for an example.

\subsection{$\C!_G(\cl_{\valgroup})$ is closed}\label{sec:closed-admissible} 

\begin{prop}\label{prop:Cclosed} If $\valgroup$ is dense in $\R$, then the collection $\C!_G(\cl_{\valgroup})$ is closed.
\end{prop}

\begin{proof} Let $\mu = \mu!_D$ be an element of $\C!_G(\cl_{\valgroup})$ for a $G$-admissible $\valgroup$-rational divisor $D \in \stable!_G(\cl_{\valgroup})$. Let $\pi=(\ssF^c, F)$ be an ordered bipartition of $V$ and $\ssE_\pi=E(\ssF^c,F)$. By Proposition~\ref{refineV}, it is enough to show that $\mu!_{\pi}$ belongs to $\C!_G(\cl_{\valgroup})$.

If $\ssE_\pi=\emptyset$, then $\mu!_{\pi}=\mu$, and thus $\mu!_{\pi}\in\C!_G(\cl_{\valgroup})$. Assume $\ssE_\pi\neq\emptyset$. Choose a positive real number $\epsilon\in\valgroup$ small enough so that the function $\ssf$ on $V$ that takes value $0$ on $\ssF^c$ and value $\epsilon$ on $F$ yields a small $G$-admissible chip-firing move 
$D'=D+\div(f;D)$ of $D$. Clearly, $D'\in\stable!_G(\cl_{\valgroup})$. Since $\epsilon$ is small and nonzero, $\ssG_{D'}$ does not contain any edge of $\ssE_\pi$. From Proposition~\ref{prop:splitG} 
we deduce that $\mu!_{D'}$ is $\pi$-split. It remains to show that $\mu!_{D'} = \mu!_\pi$.

 Since $\mu!_{D'}$ and $\mu!_\pi$ are $\pi$-split, it is enough to show that 
 \[
\mu!_{D'}(I) = \begin{cases}
\mu(I) & \textrm{ if } I\subseteq \ssF^c\\
\mu(I\cup \ssF^c)-\mu(\ssF^c) & \textrm{ if } I\subseteq F.
\end{cases}
\]
Now,
\[
\deg((\ssD')^I)=
\begin{cases}
    \deg(\ssD^I)&\text{if }I\subseteq \ssF^c,\\
    \deg((\ssD)^I)+\sse_{D}(I,\ssF^c)&\text{if }I\subseteq F.
\end{cases}
\]
Thus, if $I\subseteq \ssF^c$, then $\mu!_{D'}(I)=\mu(I)$. On the other hand, assume $I\subseteq F$. By definition,
\[
\mu(I\cup \ssF^c)-\mu(\ssF^c)=\deg(\ssD^{I \cup \ssF^c}) 
+\frac 12 e(I \cup \ssF^c,\ssI^c\cap F)
-\deg(\ssD^{\ssF^c})-\frac 12 e(\ssF^c,F).
\]
Since $I\subseteq F$, we have 
\[
\sse_{D''}(I \cup \ssF^c,\ssI^c\cap F)
-\sse_{D''}(\ssF^c,F)=\sse_{D''}(I,\ssI^c)
-2\sse_{D''}(I,\ssF^c)
\]
for each $G$-admissible divisor $D''$, in particular, for $D''=0$ and $D''=D$. We infer that
\begin{align*}
  \deg(\ssD^{I \cup \ssF^c}) - \deg(\ssD^{\ssF^c}) &=
  \deg(\ssD_{I \cup \ssF^c}) - \deg(\ssD_{\ssF^c})-\sse_{D}(I\cup \ssF^c,\ssI^c\cap F)+\sse_{D}(\ssF^c,F)\\
  &=\deg(\ssD_{I})+\delta!_D(I,\ssF^c)-\sse_D(I,\ssI^c)+2\sse_D(I,\ssF^c)\\
  &=\deg(\ssD^{I}) + e(I,\ssF^c)+ \sse_{D}(I,\ssF^c).
  \end{align*}
Therefore,
\[
\mu(I\cup \ssF^c)-\mu(\ssF^c) = \deg(\ssD^{I}) 
+\frac 12 e(I,\ssI^c)+\sse_{D}(I,\ssF^c)=\deg((\ssD')^I)
+\frac 12 e(I,\ssI^c)=\mu!_{D'}(I),
\]
finishing the proof.
\end{proof}

\subsection{$\C!_G(\cl_{\valgroup})$ is complete and positive}\label{sec:complete-admissible}

\begin{prop}\label{prop:Csimple} If $G$ is connected, then the collection $\C!_G(\cl_{\valgroup})$ is simple. If in addition $\valgroup$ is dense in $\R$, then the collection $\C!_G(\cl_{\valgroup})$ is complete.
\end{prop}

\begin{proof} Let $\mu\coloneqq \mu!_D$ for $D\in\stable!_G(\cl_{\valgroup})$, and let 
$\pi=(F,\ssF^c)$ be a bipartition of $V$ such that $\mu!_\pi = \mu$. To show the first statement, it is enough to show the 
existence of $D'\in\stable!_G(\cl_{\valgroup})$ such that $\mu'\coloneqq \mu!_{D'}$ satisfies $\mu' \neq \mu$ and $\mu!_{\pi}' =\mu$.

Since $G$ is connected, $E(F,\ssF^c)\neq\emptyset$. We use the notation introduced in Section~\ref{sec:admissible-chip-firing}. Let $D'$ be the (flipping) $G$-admissible chip-firing move of $D$ by distance $\ml_{D}(\ssF^c)$ induced by $\ssF^c$, and put $\mu'\coloneqq\mu!_{D'}$. Since $D$ is $\valgroup$-rational, $\ml_{D}(\ssF^c)\in\valgroup$, and thus $D'$ is $\valgroup$-rational and $\mu'\in\C!_G(\cl_{\valgroup})$. Also, 
$E_{D'}(F,\ssF^c)\neq\emptyset$, and thus 
$\mu'$ is not $\pi$-split. In addition, $D$ is a small $G$-admissible chip-firing move of $D'$ by $\ml_{D}(\ssF^c)$ induced by $F$, and thus, as seen in the proof of 
Proposition~\ref{prop:Cclosed}, we have $\mu=\mu!_{\pi}'$. Finally, $\mu'\neq\mu$ because $\mu$ is $\pi$-split and $\mu'$ is not.

As for the second statement, if $\valgroup$ is dense in $\R$, then $\C!_G(\cl_{\valgroup})$ is closed by Proposition~\ref{prop:Cclosed}, and hence the same argument shows that the collection $\C!_G(\cl_{\valgroup})$ is complete.
\end{proof}

\begin{remark}\label{rmk:mlD}
The proofs of Proposition~\ref{prop:Cclosed}~and~\ref{prop:Csimple} yield that, for each $D\in\stable!_G(\cl)$ and bipartition $\pi=(\ssF^c,F)$ of $V$ such that $\mu!_D$ is equal to its $\pi$-splitting, the quantity $\ml_{D}(\ssF^c)$ is the largest positive real number such that $\mu!_{D_t}=\mu!_{\ssD}$ for each $t\in [0,\ml_D(\ssF^c))$ where $\ssD_t\coloneqq D+\div(t\one!_{F};D)$.
\end{remark}

\begin{prop}\label{prop:Cpositive} If $G$ is connected and $\cl_{\valgroup}$ is nonempty, then the collection $\C!_G(\cl_{\valgroup})$ is positive.    
\end{prop}

\begin{proof} Since $\C!_G(\cl_{\valgroup})$ is simple by Proposition~\ref{prop:Csimple}, and $\cl_{\valgroup}\neq\emptyset$, it follows that $\C!_G(\cl_{\valgroup})$ contains simple elements. Now, for each $G$-admissible divisor $D$ such that $\mu!_{D}$ is simple, the spread of $\mu!_{D}$ is at least $1/2$.    
\end{proof}

\subsection{Proof of Theorem~\ref{thm:semistability-tiling-sccp}} It follows from 
Propositions~\ref{prop:Cseparated}, \ref{prop:Cclosed}, \ref{prop:Csimple}, and~\ref{prop:Cpositive} that $\C!_G(\cl_{\valgroup})$ is separated, closed, complete and positive.\qed

\subsection{Proof of Theorem~\ref{thm:semistability-tiling}} We combine Theorem~\ref{thm:semistability-tiling-sccp} with Theorems~\ref{thm:tiling},~\ref{thm:tiling2}~and~\ref{thm:semistability-polytope}, and deduce that the semistability polytopes associated to the $G$-admissible divisors in $\stable!_G(\cl_{\Lambda})$ form a tiling of the full space $H_d$. \qed


\section{Admissible divisors and reduction of linear series} \label{sec:tropicalization}

Let $\K$ be an algebraically closed field with a nontrivial non-Archimedean valuation 
$\val$. We assume that $\K$ is complete with respect to $\val$ and we denote by $\k$ the residue field of $\K$, which is also algebraically closed. Let $\varR$ be the valuation ring of $(\K,\val)$, and denote by $\varm$ the maximal ideal of $\varR$. We denote by $\valgroup$ the value group of $\K$.  Notice that $\valgroup$ is a nontrivial additive divisible subgroup of $\R$, whence dense in $\R$.

 Let $\varX$ be a smooth proper connected curve over $\K$ and denote by $\varX^{\an}$ the Berkovich analytification of $\varX$. We assume familiarity with the theory of Berkovich analytic curves, and refer to~\cite[Section 4]{AB15}, ~\cite[Section 5]{BPR} and~\cite{Duc},  
 which contain what we need here.

 A \emph{semistable vertex set} for $\varX^{\an}$ is a finite set $V$ of type 2 points on $\varX^{\an}$ whose complement $\varX^{\an} \setminus V$ is a disjoint union of finitely many open annuli and infinitely many open disks. A \emph{semistable model} for $\varX$ is an integral proper relative curve $\frakX$ over $\varR$ with generic fiber $\frakX_\eta=\varX$ and special fiber $\frakX_0$ which is reduced, connected and has nodal singularities. Any irreducible component of the special fiber $\frakX_0$ of a semistable model $\frakX$ gives a valuation on $\K(\varX)$ and defines a point of type 2 in $\varX^{\an}$. The set $V$ of points in $\varX^{\an}$ associated to the irreducible components of $\frakX_0$ gives a semistable vertex set for $\varX^{\an}$. Moreover, this process provides a bijection between semistable vertex sets of  $\varX^{\an}$ and semistable models of $\frakX$, see~\cite[Thm.~5.38]{BPR}
 
 A semistable vertex set $V$ gives rise to a \emph{skeleton} $\Gamma$ for $\varX^{\an}$, defined as the union in $\varX^{\an}$ of $V$ and the skeletons of the open annuli in $\varX^{\an} \setminus V$. Using the canonical metric on the skeletons of the open annuli, we can view the skeleton as a metric graph canonically embedded in $\varX^{\an}$. 
  The underlying graph $G=(V,E)$ is connected, and has vertex set $V$ and edge set $E$ in bijection with the  set of open annuli in $\varX^{\an} \setminus V$. There is an edge between vertices $u,v\in V$ for each open annulus in $\varX^{\an} \setminus V$ whose closure in $\varX^{\an}$ contains $u$ and $v$. The edge length function $\ell\colon E \to (0,+\infty)$ associates to each edge the modulus of the corresponding annulus. In particular, all the edge lengths are in $\valgroup$. The graph $G$ is the dual graph of the special fiber $\frakX_0$ of the semistable model $\frakX$  that corresponds to the semistable vertex set $V$, with vertices in bijection with the irreducible components of $\frakX_0$ and with edges in bijection with the nodes of $\frakX_0$. There is an edge $e=uv$ in $G$ for each node that lies on both of the irreducible components associated to $u$ and $v$. Via this correspondence, the edge length function associates to each edge the singularity degree in $\frakX$ of the corresponding node. 

Let $\Gamma$ be a metric graph skeleton of  $\varX^{\an}$ with underlying graph $G=(V, E)$ and edge length function $\ell\colon E \to \valgroup$.

For each $v\in V$, let $\k(v)$ be the residue field of $v$ seen as a point of type 2 in $\varX^{\an}$. The field $\k(v)$ is of transcendental degree one over $\k$ and we denote by $\varC_v$ the corresponding smooth proper curve over $\k$. In the semistable model that corresponds to $V$, the curve $\varC_v$ is the normalization of the irreducible component in $\frakX_0$ associated to $v$, so $\k(v)$ is the function field of this component.
  
  For each point $x$ of $\varX^{\an}$, denote by $\val_x$ the corresponding valuation on $\K(\varX)$, extending $\val$. If $x\in\varX^{\an}\setminus \varX(\K)$, then $\val_x$ is a map from $\K(\varX)$ to $\R\cup\{+\infty\}$. Otherwise, that is, if $x\in \varX(\K)$, we view $\val_x$ as a map from $\K(\varX)$ to $\R\cup\{\pm\infty\}$ which takes value $\val_x(\varf)=\val(\varf(x))$  if $\varf$ is defined at $x$, and $-\infty$ otherwise.  If $x$ is a type 2 point, in particular, if $x$ is a vertex of $G$, the image of $\val_x$ coincides with the value group $\valgroup$ of $\val\colon \K \to \R\cup\{+\infty\}$. 

  For each $\varf \in \K(\varX)$, we denote by $\ev_{\varf} \colon \varX^{\an} \to \R\cup\{\pm \infty\}$ the evaluation map that sends each point $x\in \varX^{\an}$ to $\val_x(\varf)$.
  
  Let $\tau \colon \varX^{\an} \to \Gamma$ be the canonical retraction map from $\varX^{\an}$ to $\Gamma$. We call $\tau$ the \emph{tropicalization map}. The restriction of $\tau$ to $\varX(\K)$ is compatible with the specialization map from the generic fiber $\frakX_\eta$ to $\frakX_0$. Extending by linearity, we get a  tropicalization map $\tau \colon \Div(\varX) \to \Div(\Gamma)$ that sends  a divisor $\varD = \sum_{x\in \varX(\K)}\varD(x)(x)$ on $\varX$ to the divisor $\tau(\varD) = \sum_{x\in \varX(\K)}\varD(x)(\tau(x))$. Since $\tau(x)$ is $\valgroup$-rational for each $x\in \varX(\K)$, we have that $\tau(\varD)$ is $\valgroup$-rational.

\subsection{Tropicalization and reduction of linear series} \label{sec:tropicalization2}
 
For each nonzero rational function $\varf \in\K(\varX)$, we denote by $\trop(\varf)$ the \emph{tropicalization of $\varf$}, defined as the restriction to $\Gamma$ of the evaluation map $\ev_{\varf} \colon \varX^{\an} \to \R\cup\{\pm \infty\}$. This gives an element of $\Rat(\Gamma)$.

 Let $\varD$ be divisor of degree $d$ on $\varX$ and let $\tau(\varD) = \sum_{x\in \varX(\K)} \varD(x)(\tau(x))$ be the $\valgroup$-rational divisor on $\Gamma$ obtained by the tropicalization of $\varD$.  Denote by $\cl$ the class of $\tau(\varD)$ in $\Pic^d(\Gamma)$, and notice that $\cl_\valgroup$ is nonempty.

 Let $\varH$ be a vector subspace of dimension $r+1$ of the space of global sections of the line bundle $\varL = \mathcal O(\varD)$, which we view in $\K(\varX)$:
 \[\varH \subseteq \bigl \{\varf \in\K(\varX) \,\st\, \div(\varf) + \varD \geq 0\bigr\}.\]
The pair $(\varD, \varH)$ defines a linear series 
on $\varX$. By the Specialization Lemma~\cite[Thm.~4.5]{AB15}, 
\[\tau(\varD) + \div(\trop(\varf)) = \tau(\varD + \div(\varf)).\]

\begin{convention}\label{conv:section}
In the following, we will choose a section of the valuation $\val \colon \K\setminus \{0\} \to \valgroup$.
In other words, for each $\lambda \in\valgroup$, we fix $\ssa_\lambda\in\K$ with $\val(\ssa_\lambda)=\lambda$ in such a way that the map $\lambda \mapsto \ssa_\lambda$ induces a homomorphism of groups $\Lambda \to \K\setminus \{0\}$. (Since $\Lambda$ is divisible, such a section exists.)
\end{convention}

For each nonzero $\varf \in \K(\varX)$, putting $a_v\coloneqq  a_{\val_v(f)}\in \K$, the element $a_v^{-1}\varf$ has valuation $\val_v(a_v^{-1}\varf)=0$, and reduction in the residue field $\k(v)$ that we denote by $\varf_v$, for each $v\in V$. We call $\varf_v$ the \emph{relative reduction of $\varf$ at $v$}. 

Note that $\varf_v$ is nonzero. Choosing another section of $\val$ would yield a nonzero $\k$-multiple of $\varf_v$. If $\val_v(f)=0$ then $a_v=1$ and $\varf_v$ is the usual (nonrelative) reduction of $\varf$ in $\k(v)$. We will generalize relative reduction in Definition~\ref{def:rel-red}. Under this generalization, $\varf_v$ will be the reduction of $\varf$ at $v$ relative to $\trop(\varf)$.

Denote by $\VS_v$ the $\k$-vector space $\VS_v$ of $\k(v)$ spanned by the relative reductions $\varf_v$ at $v$ of nonzero functions $\varf\in\varH$ \cite[Section 4.4]{AB15}. The space is independent of the choice of the section of $\val$. We have the following result \cite[Lemma 4.3]{AB15}: 
		
\begin{lemma}\label{lem:reduction}
Notation as above, for each vertex $v$  of $G$, the $\k$-vector subspace $\VS_v \subset \k(v)$ has dimension $r + 1$.
\end{lemma}

\begin{defi}\label{defi:reduction}
	We associate to the linear series $(\varD, \varH)$ the $\k$-vector space 
	 \[\VS \coloneqq \bigoplus_{v\in V} \VS_v,\]
	 defined by the reductions of $\varH$ at the vertices $v$ in $G$. 
\end{defi}

\begin{defi}\label{defi:action}
Let $\k^V\coloneqq \prod_{v\in V}\k$. For each element $\varphi\in\k^V$, we denote by $\varphi_v$ the element in $\k$ given by the $v$-component of $\varphi$, for each $v\in V$. Each $\varphi\in\k^V$ corresponds to an endomorphism of $\VS$ that by an abuse of notation we will also denote by $\varphi$. 
\end{defi}
The following is straightforward.
\begin{prop} For each $\varphi \in \k^V$, the endomorphism $\varphi\colon \VS \to \VS$ is an automorphism if an only if $\varphi_v\neq 0$ for every $v$.
\end{prop}

\subsection{Subspace of $\VS$ associated to a $G$-admissible divisor in $\stable!_G(\cl)$}\label{sec:subspace-of-U} 
We generalize Lemma~\ref{lem:reduction} by associating an $(r+1)$-dimensional subspace of $\VS$ to each $G$-admissible divisor in $\cl$. 

We need the following definition (we will be interested in the $\valgroup$-rational divisor $\tau(\varD)$). 

\begin{defi}[$G$-admissible rational functions]
For each divisor divisor $D$ on $\Gamma$, we will denote by $\adfun_{G}(D)$ the set of all rational functions $h\in \Rat(\Gamma)$ such that $\div(h) +D$ is $G$-admissible. We call them \emph{$G$-admissible rational functions associated to $D$}. If $D$ is $\valgroup$-rational, we denote by 
$\adfun_{G}(D;\valgroup)$ the subset of $h\in\adfun_{G}(D)$ such that $h$ is $\valgroup$-rational.
\end{defi}

Notice that, by Proposition~\ref{prop:rational-admissible extension}, since $G$ is connected, if $D$ is 
$\valgroup$-rational, then a rational function $h\in\adfun_{G}(D)$ is in $\adfun_{G}(D;\valgroup)$ if and only if 
$h(u)-h(v)\in\Lambda$ for each $u,v\in V$.

 Let $h \in \adfun_{G}(\tau(\varD);\valgroup)$. Denote by $\varM_h$ the set of all functions $\varf\in \varH$ with $\trop(\varf)(v) \geq h(v)$ for every $v\in V$. By Lemma~\ref{lem:varH}, this is equivalent to $\trop(\varf) \geq h$ pointwise on $\Gamma$, that is,
\[\varM_h \coloneqq \bigl\{\varf\in \varH \,\st\, \trop(\varf)(v)  \geq h(v)\quad \forall v\in V\bigr\} =  \bigl\{\varf\in \varH \,\st\, \trop(\varf) \geq h\bigr\}.\]

Notice that $\val_v(\ssa_{h(v)}^{-1}\varf)\geq 0$
for each $h \in \adfun_{G}(\tau(\varD);\valgroup)$ and $\varf \in \varM_h$. 

\begin{defi}[Relative reduction map $\red!_h$] \label{def:rel-red}The \emph{reduction map relative to $h$} is the map 
\[\red!_h \colon \varM_h \to \VS\]
defined as follows: For each $\varf\in\varM_h$, let $\red!_{h}(\varf)$ be the element of $\VS$ whose component $\ssub{(\red!_h(\varf))}!_v$ is the reduction of 
$\ssa_{h(v)}^{-1}\varf$ to $\k(v)$ for each $v\in V$. In particular, if $\trop(\varf)(v)>h(v)$, then $\ssub{(\red!_h(\varf))}!_v=0$.
\end{defi} 

Each nonzero $\varf$ in $\varH$ is in $\varM_h$ for 
$h\coloneqq \trop(\varf)$. In this case, the relative reduction $\varf_v$ of $\varf$ at $v$ we have defined earlier coincides with the reduction $\ssub{(\red!_h(\varf))}!_v$ of $\varf$ relative to $\trop(\varf)$.

\begin{defi}[Reduction of $\varH$ relative to $h$] We define the reduction of $\varH$ relative to $h$ as the subset $\ssW_h \subset \VS$ consisting of the reductions of $\varf \in \varM_h$ relative to $h$, that is, 
 \[\ssW_h \coloneqq \bigl\{\, \red!_h(\varf) \,\st\, \varf\in \varM_h\,\bigr\}. \qedhere\] 
\end{defi}
 
Different choices of the section $\lambda \mapsto \ssa_\lambda$, $\lambda\in \valgroup$, will produce different subspaces $\ssW_h$, but all the outcomes $\ssW_h$, obtained via this reduction process, lie on the same orbit of $\VS$ by the natural action of $\prod_{v\in V}\k^\times$ defined in Definition~\ref{defi:action}.

Note that for each constant $c$ in the value group $\valgroup$, the space $\ssW_{h+c}$ coincides with $\ssW_h$. This implies that $\ssW_h$ depends only on $D=\tau(\varD)+\div(h)$. We can thus set $\ssW_D\coloneqq\ssW_h$. We use both notations, $\ssW_D$ and $\ssW_h$, according to the context.

 The function $h$ induces a modified norm on rational functions at points on $\varX^{\an}$. 
 More precisely, we define the \emph{$h$-norm}
 \begin{align*}
 \norm{\cdot}{h} \colon \K(\varX) \times \varX^{\an} &\to \R\cup\{\infty\}\\
    (\varf, x) &\mapsto \norm{\varf(x)}{h} \coloneqq e^{h(\tau(x))} \norm{\varf(x)}{}.
 \end{align*} 
 Here, $\norm{\varf(x)}{}$ is the norm of $\varf$ induced by the point $x\in \varX^{\an}$, given by $\norm{\varf(x)}{} = \exp{(-\val_x(\varf))}$.  
 
 With this notation,
 \[
 \varM_h = \bigl\{\varf\in \varH \,\st\, \norm{\varf(x)}{h} \leq 1 \quad \forall \,x\in\Gamma\subset\varX^{\an}\bigr\}.
 \]
 Notice that $\varM_h$ is an $\varR$-submodule of $\varH$, whence torsion-free. Also, $\K\cdot \varM_h =\varH$, as for each $\varf\in \varH$, taking $a \in \K$ with $|a|$ very small, we have $a\varf\in \varM_h$. Finally, since $\K$ is complete and $\varH$ has finite dimension over $\K$, we get that $\varH$ is complete, whence $\varM_h$ is complete as well.

\smallskip
We have the following generalization of Lemma~\ref{lem:reduction} (see Proposition~\ref{prop:v-reduced-iso}).
  			
 \begin{lemma}[Reduction Lemma]\label{lem:reduction_lemma} Notation as above, $\varM_h$ is a free $\varR$-module of rank $r+1$ and the reduction map $\red!_h$ factors through an isomorphism of $\k$-vector spaces $\rquot{\varM_h}{\varm\varM_h} \simeq \ssW_h$. As a consequence, the subspace $\ssW_h$ is a $\k$-vector subspace of $\VS$ of dimension $r+1$.
 \end{lemma}
 
  \begin{proof}  We give $\VS$ the structure of an $\varR$-module via the quotient map $\varR\to\k$. To prove that $\ssW_h$ is a $\k$-vector subspace of $\VS$ and that $\red!_h$ factors through a surjection $\rquot{\varM_h}{\varm\varM_h} \to \ssW_h$, we only need to prove that $\red!_h\colon\varM_h\to\VS$ is a homomorphism of $\varR$-modules. This is indeed the case, as
 \[\ssub{(\red!_h(a\varf+b\varg))}_v = \tilde a\ssub{(\red!_h(\varf))}_v+ \tilde b\ssub{(\red!_h(\varg))}_v 
 \]
 for each $\varf,\varg\in\varM_h$ and $a,b\in\varR$, where $\tilde a$ and $\tilde b$ are the reductions of $a$ and $b$, respectively.

 To show the map induced by $\red!_h$ is an isomorphism, it will be enough to prove that the kernel of $\red!_h$ is $\varm \varM_h$. Consider an element $\varf\in \varM_h$ which is sent to 0 in $\VS$. This means $\trop(\varf)(v)>h(v)$ for all vertices $v\in V$. But then, since $\valgroup$ contains sufficiently small elements, there is $a\in\varm$ such that $\trop(a^{-1}\varf)(v)>h(v)$ for each $v\in V$. This implies that $\varf \in a\varM_h\subseteq\varm \varM_h$, as required.

Finally, we show that $\varM_h$ is a free $\varR$-module of rank $r+1$, from which we deduce that $\ssW_h$ is of dimension $r+1$. 
This will finish the proof of the lemma.
 
We claim first that $\dim_{\k}(\ssW_h)\leq r+1$. Indeed, let $\ssf_1, \dots, \ssf_p$ be $\k$-linearly independent elements of $\ssW_h$, and lift them to $\varf_1, \dots, \varf_p$ in $\varM_h$. Since $\K\cdot \varM_h = \varH$, it is enough to prove that $\varf_1, \dots, \varf_p$ are $\K$-linearly independent. 

But, if there existed an equation $\ssa_1\varf_1+\dots +\ssa_p\varf_p =0$ for $\ssa_j\in\K$ not all zero, we could assume without loss of generality that $\min_{j}(\val(\ssa_j))=0$, and get 
 \[
 \ssub{\tilde a}_1\ssf_1+\dots+ \ssub{\tilde a}_p\ssf_p=\ssub{\tilde a}!_1\red_h(\varf_1)+\dots+\ssub{\tilde a}!_p\red_h(\varf_p)=0,
 \]
 with $\ssub{\tilde a}!_j$,  the reductions of the $\ssa_j$, not all zero. This would be not possible, the elements $\ssf_1, \dots,\ssf_p$ being independent.
 
 Let $p = \dim_{\k}(\ssW_h)$. Choose a basis $\ssf_1, \dots, \ssf_p$ of $\ssW_h$ and lift it to $\varf_1, \dots, \varf_p$ in $\varM_h$. We claim now that $\varM_h = \varR \varf_1 \oplus \dots \oplus \varR \varf_p$. Once this has been proved, the equality $\K\cdot \varM_h = \varH$ implies that $p=r+1$, as required.

 To prove the claim, we only  need to show that $\varf_1, \dots, \varf_p$ generate $\varM_h$ as an $\varR$-module. Let $\varN$ be the $\varR$-module generated by $\varf_1, \dots, \varf_p$, so that we have $\varM_h = \varN+ \varm \varM_h$. Proceeding by induction, this yields $\varM_h= \varN+ \varm^n\varM_h$ 
 for each $n\in \N$. Since $\varM_h$ is complete, this proves that $\varM_h=\varN$, as required. 
 \end{proof}

\subsection{Reduced divisors, II}\label{sec:v-reduced-2} 
Notice that Lemma~\ref{lem:reduction} was not used in the proof of Lemma~\ref{lem:reduction_lemma}. Actually, we show now that the former can be recovered from the latter.

\begin{prop}\label{prop:v-reduced-iso} Let $\cl$ be the class of $\tau(\varD)$ for a divisor $\varD$ on $\varX$. Let $v$ be a vertex of $V$ and $D$ be the $v$-reduced divisor in $\cl$. Then, $D\in\stable!_G(\cl_{\valgroup})$ and the projection map $\ssW_{\ssD} \to \VS_v$ is an isomorphism.
\end{prop}

\begin{proof} Since $\tau(\varD)$ is $\valgroup$-rational, it follows from Proposition~\ref{prop:v-reduced-adm} that $D\in\stable!_G(\cl_{\valgroup})$. Let $h$ be an element of $\Rat(\Gamma)$ such that $\ssD=\tau(\varD)+ \div(h)$ and $h(v)=0$. Then, $h$ is $\valgroup$-rational, and moreover, $h\rest{V}$ takes values in $\valgroup$. 

We first show that $\ssW_{\ssD} \to \VS_v$ is injective. Let $\varf$ be an element of $\varM_h$ with $(\red_h(\varf))_v=0$, that is, $\trop(\varf)(v)>h(v)$. Let $I$ be the set of vertices $u$ in $V$ such that $\trop(\varf)(u)>h(u)$. We have $v\in I$. We will show that $I=V$, which implies that $\red_h(\varf)=0$.  

Let $\epsilon\coloneqq \min_{u\in I}(\trop(\varf)(u)-h(u))$. We have $\epsilon\in\valgroup$. Let $h'\in\Rat(\Gamma)$ be so that $h'-h$ is the $G$-admissible extension of $\epsilon\one!_{I}$ with respect to $D$. Then, 
$\trop(\varf)(u)\geq h'(u)$ for each $u\in V$, and hence 
$\trop(\varf)\geq h'$ by Lemma~\ref{lem:varH}. Assume for the sake of a contradiction that $I\neq V$. Then, $D'\coloneqq\tau(\varD)+\div(h')$ is the $G$-admissible chip-firing move of $D$ by $\epsilon$ induced by $V\setminus I$. Since $v\in I$ and $D$ is $v$-reduced, $D'$ cannot be effective outside $v$.  Therefore, since $D'$ is $G$-admissible, and $D$ is effective outside $v$, we should have $D'(w)<0$ for some $w\in V\setminus I$. Now, using $\trop(\varf)(w)= h(w)=h'(w)$ and $\trop(\varf)\geq h'$, we deduce that $\div(\trop(\varf))(w)\leq \div(h')(w)$. We infer that
\[0 \leq \div(\trop(\varf))(w)+\tau(\varD)(w)\leq D'(w)<0.\] 
This contradiction proves the claim.

We now show that $\ssW_{\ssD} \to \VS_v$ is surjective. Let $\varf$ be an element of $ \varH$ with $\trop(\varf)(v)=h(v)=0$. Then, we have $\trop(\varf) \geq h$; see~\cite[Lemma 20]{Ami-red}. It follows that $\varf \in \varM_h$ and $(\red_h(\varf))_v$ coincides with the relative reduction of $\varf$ at $v$. Since $\VS_v$ is generated by the relative reduction at $v$ of $\varf \in \varH$ with $\trop(\varf)(v)=0$, and the map $\ssW_{\ssD} \to \VS_v$ is $\k$-linear, the surjectivity follows.
\end{proof}


\section{Tilings induced by tropicalization of linear series} \label{sec:tilings-simplex-mixed}
We introduce in this section a second class of tilings that arise from degenerations of linear series on curves.

Let $V$ be a finite nonempty set. Let $\k$ be a field. For each $v\in V$, let $\VS_v$ be a finite-dimensional vector space over $\k$. The example of interest to us is the ones arising from the reduction of linear series discussed in Section~\ref{sec:tropicalization}. Define
\[\VS\coloneqq \bigoplus_{v\in V} \VS_v.\]

For each subset $I\subseteq V$, let 
\[\VS_I\coloneqq \bigoplus_{v\in I}\VS_v\]
 and denote by 
 \[
 \iota!_{I}\colon\VS_{I}\to\VS\quad\text{and}\quad
 \proj{I}\colon \VS\to \VS_{I}
 \]
 the corresponding insertion and projection 
 maps, respectively. Clearly, the composition $\proj{I}\circ \iota!_I$ is the identity. By convention, we set 
\[\VS_\emptyset\coloneqq (0).\]  
 
 We have a natural exact 
 sequence,
 \begin{equation} \label{eq:U-short-exact-sequence}
 0 \to \VS_I\to \VS \to \VS_{I^c} \to 0,
 \end{equation}
 where the second map is $\iota!_I$ and the third is $\proj{I^c}$.

As in Definition~\ref{defi:action}, we view each element $\varphi$ of $\k^V$ as an endomorphism of $\VS$ that by an abuse of notation we will also denote by $\varphi$. This endomorphism is an automorphism if and only if $\varphi_v\neq 0$ for every $v$.

\subsection{Adjoint modular pairs arising from subspaces} \label{sec:modular-subspace}  Notation as above, let $\ssW\subseteq \VS$ be a vector subspace.

For each $I\subseteq V$, we put 
$\ssW_I\coloneqq \iota!_I\circ \proj{I}(\ssW)$, and denote by $\nu!_{\ssW}^*(I)$ its dimension, that is,
\[
\ssW_I \coloneqq \iota!_I\circ \proj{I}(W)\subseteq\VS \qquad \textrm{and}\qquad \nu!_{\ssW}^*(I) \coloneqq \dim_{\k}(\ssW_I).
\]
Similarly, we denote by $\ssW^I$ the intersection of $\ssW$ with $\iota!_I(\VS_I)$, and by $\nu!_{\ssW}(I)$ its dimension, that is, 
\[
\ssW^I\coloneqq \ssW \cap\iota!_I(\VS_I) \subseteq \VS \qquad \textrm{and} \qquad \nu!_{\ssW}(I) \coloneqq \dim_{\k}(\ssW^I).
\]
We have a short exact sequence 
\begin{equation}\label{eq:short-exact}
0 \to \ssW^{I^c} \to \ssW \to \ssW_{I} \to 0,
\end{equation}
where the second map is the inclusion and the third map is $\iota!_I\circ \proj{I}$.

The construction is functorial: given $\varphi\in\k^V$ and subspaces $\ssW_1,\ssW_2\subseteq\VS$ with $\varphi\ssW_1\subseteq\ssW_2$, we have 
$\varphi\ssW_1^I\subseteq\ssW_2^I$ and $\varphi\ssW_{1,I}\subseteq\ssW_{2,I}$ for each $I\subseteq V$.

Note that with these notations, for $\ssW=\VS$, we have $\VS^I=\VS_I$, and the exact sequence \eqref{eq:short-exact} is identical to \eqref{eq:U-short-exact-sequence}.

\begin{prop} Notation as above, let $\ssW \subseteq \VS$ be a vector subspace of dimension $r+1$. Then, $(\nu!_{\ssW}, \nu!_{\ssW}^*)$ is an adjoint modular pair of range $r+1$.
\end{prop}

\begin{proof} From the short exact sequence~\eqref{eq:short-exact}, we get $\ssW_{I} \simeq \rquot{\ssW}{\ssW^{I^c}}$. We deduce the equation 
\[\nu!_{\ssW}(\ssI^c) + \nu!_{\ssW}^*(I) = r+1.\]
It will be thus enough to show that $\nu!_{\ssW}$ is supermodular, that is, for $I, J \subset V$, we have 
\[
\dim_{\k} \ssW^{I}+\dim_{\k} \ssW^{J}\leq\dim_{\k} \ssW^{I\cup J}+\dim_{\k} \ssW^{I \cap J}.
\]
This follows from the fact that
$\ssW^{I}$ and $\ssW^{J}$ are subspaces of $\ssW^{I \cup J}$
with intersection $\ssW^{I \cap J}$.
\end{proof}

\begin{defi}[Base polytope associated to a subspace $\ssW \subset \VS$] \label{def:P_W}
We denote by $\Q_{\ssW}$ the base polytope associated to the pair $(\nu!_{\ssW}, \nu!_{\ssW}^*)$. Thus,
\[
\Q_{\ssW}=\left\{q\in \R^V\,\st\,
\nu!_{\ssW}(I)\leq q(I)\leq \nu!_{\ssW}^*(I)\text{ for each } I\subseteq V\right\}. \qedhere
\]
\end{defi}

Note that the values taken by any point $q\in \Q_{\ssW}$ at any $v\in V$ are all non-negative. Also, the sum $\sum_{v\in V}q(v)$ coincides with the range of $\nu!_{\ssW}$, which is $r+1$. Denote by $\Delta_{r+1}$ the standard simplex in $H_{r+1}$ given by
\[\Delta_{r+1} \coloneqq \bigl\{\,q\in \R^V \,\,\st\,\,  q(v)\geq  0\,\, \textrm{for all} \,\, v\in V, \,\, \textrm{and}\, q(V) =r+1\,\bigr\}.\]
We infer the following.  

\begin{prop} The base polytope $\Q_{\ssW}$ is a subset of the standard simplex $\Delta_{r+1}$ in $H_{r+1}$.
\end{prop} 

Note as well that $\Q_{\varphi\ssW}=\Q_{\ssW}$ for each $\varphi\in\k^V$ with $\varphi_v\neq 0$ for every $v$.

\subsection{Families of modular pairs coming from reductions of linear series} 
Notation as in Section~\ref{sec:tropicalization}, let $\tau\colon \varX^{\an} \to \Gamma$ be the retraction map from $\varX^{\an}$ to $\Gamma$, and let $\tau(\varD)$ be the tropicalization of $\varD$. Let $\cl$ be its class in $\Pic^d(\Gamma)$ and $\cl_{\valgroup}\subseteq\cl$ the subset of $\Lambda$-rational divisors. Denote by $\stable!_G(\cl_{\valgroup})$ the set of $G$-admissible divisors $D$ in $\cl_{\valgroup}$. Each divisor $D$ in $\stable!_G(\cl_{\valgroup})$ is of the form $\div(h)+\tau(\varD)$ for some element $h\in \adfun_G(\tau(\varD);\Lambda)$. Moreover, since $G$ is connected, $h$ is unique up to translation by a constant in $\valgroup$. By Reduction Lemma~\ref{lem:reduction_lemma}, we can associate to the $G$-admissible $\valgroup$-rational divisor $D$ a subspace of $ \VS$ of dimension $r+1$ that we denote by $\ssW_D \coloneqq \ssW_{h}$.  We denote by $(\nu!_D, \nu!_D^*)$ the adjoint modular pair corresponding to $\ssW_D$, defined in the preceding section,  and denote by $\Q_{D}$ the corresponding base polytope. 

Let $\C!_{(\varX, G)}(\varD, \varH)$ be the collection of supermodular functions $\nu!_D$ for $D$ a $G$-admissible divisor in $\cl_{\valgroup}$, that is, $\Lambda$-rational and linearly equivalent to $\tau(\varD)$. 

We have the following result.

\begin{thm}\label{thm:reduction-admissible-tiling-sccp} The collection $\C!_{(\varX, G)}(\varD, \VS)$ is separated, closed, positive and complete for the interior of the simplex $\Delta_{r+1}$, in the sense of Section~\ref{sec:admissible-families-tilings}.
\end{thm}

From the above theorem, we deduce the following result.

\begin{thm} \label{thm:reduction-admissible-tiling} The collection of polytopes $\Q_D$ associated to $\Lambda$-rational $G$-admissible divisors $D$ in the linear equivalence class of the divisor $\tau(\varD)$ gives a tiling of the standard simplex $\Delta_{r+1}$.
\end{thm}

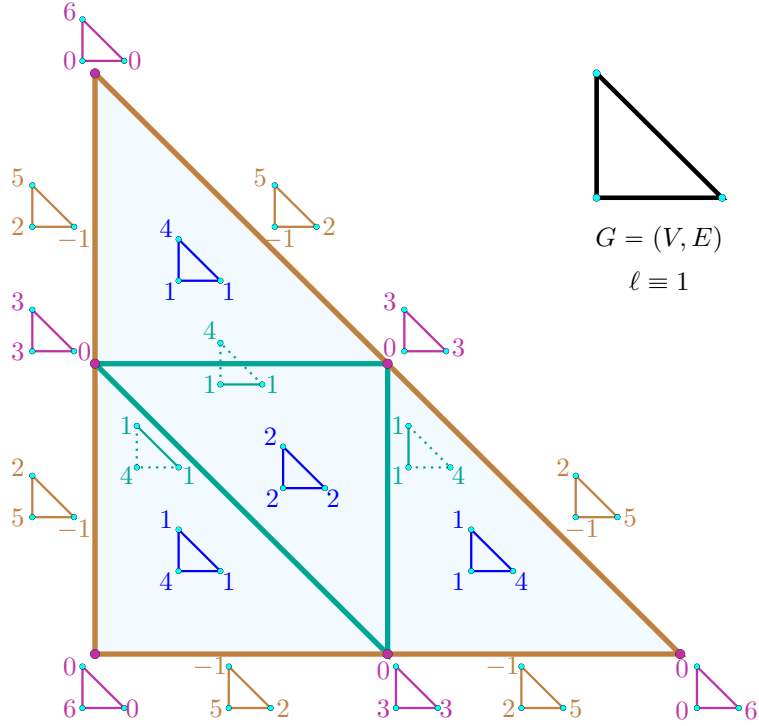
\begin{figure}[t]
\centering
\begin{tikzpicture}[scale=.55]

\draw[line width=0.6mm] (12,0) -- (12,-3) -- (15,-3) -- (12,0);
\draw[line width=0.1mm] (12,0) circle (.8mm);
\filldraw[aqua] (12,0) circle (.7mm);
\draw[line width=0.1mm] (12,-3) circle (.8mm);
\filldraw[aqua] (12,-3) circle (.7mm);
\draw[line width=0.1mm] (15,-3) circle (.8mm);
\filldraw[aqua] (15,-3) circle (.7mm);
\draw (13.5,-4) node{{\small $G=(V,E)$}};
\draw (13.5,-5) node{{\small $\ell \equiv 1$}};

\fill[lblue!5!white] (0,0) -- (0,-7) -- (7,-7) -- (0,0);
\draw[line width=0.3mm, blue] (2,-4) -- (2,-5) -- (3,-5) -- (2,-4);

\draw[line width=0.1mm] (2,-4) circle (.6mm);
\filldraw[aqua] (2,-4) circle (.5mm);
\draw[line width=0.1mm] (2,-5) circle (.6mm);
\filldraw[aqua] (2,-5) circle (.5mm);
\draw[line width=0.1mm] (3,-5) circle (.6mm);
\filldraw[aqua] (3,-5) circle (.5mm);

\draw[blue] (1.7,-3.75) node{{\small $4$}};
\draw[blue] (1.8,-5.25) node{{\small $1$}};
\draw[blue] (3.2,-5.25) node{{\small $1$}};

\fill[lblue!5!white] (0,-7) -- (0,-14) -- (7,-14) -- (0,-7);
\draw[line width=0.3mm, blue] (2,-11) -- (2,-12) -- (3,-12) -- (2,-11);

\draw[line width=0.1mm] (2,-11) circle (.6mm);
\filldraw[aqua] (2,-11) circle (.5mm);
\draw[line width=0.1mm] (2,-12) circle (.6mm);
\filldraw[aqua] (2,-12) circle (.5mm);
\draw[line width=0.1mm] (3,-12) circle (.6mm);
\filldraw[aqua] (3,-12) circle (.5mm);
\draw[blue] (1.7,-10.75) node{{\small $1$}};
\draw[blue] (1.7,-12.25) node{{\small $4$}};
\draw[blue] (3.2,-12.25) node{{\small $1$}};

\fill[lblue!5!white] (7,-7) -- (7,-14) -- (14,-14) -- (7,-7);
\draw[line width=0.3mm, blue] (9,-11) -- (9,-12) -- (10,-12) -- (9,-11);
\draw[line width=0.1mm] (9,-11) circle (.6mm);
\filldraw[aqua] (9,-11) circle (.5mm);
\draw[line width=0.1mm] (9,-12) circle (.6mm);
\filldraw[aqua] (9,-12) circle (.5mm);
\draw[line width=0.1mm] (10,-12) circle (.6mm);
\filldraw[aqua] (10,-12) circle (.5mm);
\draw[blue] (8.7,-10.75) node{{\small $1$}};
\draw[blue] (8.7,-12.25) node{{\small $1$}};
\draw[blue] (10.2,-12.25) node{{\small $4$}};

\fill[lblue!5!white] (0,-7) -- (7,-7) -- (7,-14) -- (0,-7);
\draw[line width=0.3mm, blue] (4.5,-9) -- (4.5,-10) -- (5.5,-10) -- (4.5,-9);
\draw[line width=0.1mm] (4.5,-9) circle (.6mm);
\filldraw[aqua] (4.5,-9) circle (.5mm);
\draw[line width=0.1mm] (4.5,-10) circle (.6mm);
\filldraw[aqua] (4.5,-10) circle (.5mm);
\draw[line width=0.1mm] (5.5,-10) circle (.6mm);
\filldraw[aqua] (5.5,-10) circle (.5mm);
\draw[blue] (4.2,-8.75) node{{\small $2$}};
\draw[blue] (4.25,-10.25) node{{\small $2$}};
\draw[blue] (5.7,-10.25) node{{\small $2$}};

\draw[line width=0.3mm, byzant] (-.3,1.3) -- (-.3,.3) -- (.7,.3) -- (-.3,1.3);

\draw[line width=0.1mm] (-.3,1.3) circle (.6mm);
\filldraw[aqua] (-.3,1.3) circle (.5mm);
\draw[line width=0.1mm] (-.3,.3) circle (.6mm);
\filldraw[aqua] (-.3,.3) circle (.5mm);
\draw[line width=0.1mm] (.7,.3) circle (.6mm);
\filldraw[aqua] (.7,.3) circle (.5mm);

\draw[byzant] (-.6,1.5) node{{\small $6$}};
\draw[byzant] (-.6,0.3) node{{\small $0$}};
\draw[byzant] (.95,0.3) node{{\small $0$}};

\draw[line width=0.3mm, byzant] (-.3,-14.3) -- (-.3,-15.3) -- (.7,-15.3) -- (-.3,-14.3);

\draw[line width=0.1mm] (-.3,-14.3) circle (.6mm);
\filldraw[aqua] (-.3,-14.3) circle (.5mm);
\draw[line width=0.1mm] (-.3,-15.3) circle (.6mm);
\filldraw[aqua] (-.3,-15.3) circle (.5mm);
\draw[line width=0.1mm] (.7,-15.3) circle (.6mm);
\filldraw[aqua] (.7,-15.3) circle (.5mm);

\draw[byzant] (-.6,-14.3) node{{\small $0$}};
\draw[byzant] (-.6,-15.3) node{{\small $6$}};
\draw[byzant] (.9,-15.3) node{{\small $0$}};

\draw[line width=0.3mm, byzant] (14.4,-14.3) -- (14.4,-15.3) -- (15.4,-15.3) -- (14.4,-14.3);

\draw[line width=0.1mm] (14.4,-14.3) circle (.6mm);
\filldraw[aqua] (14.4,-14.3) circle (.5mm);
\draw[line width=0.1mm] (14.4,-15.3) circle (.6mm);
\filldraw[aqua] (14.4,-15.3) circle (.5mm);
\draw[line width=0.1mm] (15.4,-15.3) circle (.6mm);
\filldraw[aqua] (15.4,-15.3) circle (.5mm);

\draw[byzant] (14.05,-14.35) node{{\small $0$}};
\draw[byzant] (14.05,-15.35) node{{\small $0$}};
\draw[byzant] (15.7,-15.35) node{{\small $6$}};

\draw[line width=0.3mm, byzant] (7.2,-14.3) -- (7.2,-15.3) -- (8.2,-15.3) -- (7.2,-14.3);

\draw[line width=0.1mm] (7.2,-14.3) circle (.6mm);
\filldraw[aqua] (7.2,-14.3) circle (.5mm);
\draw[line width=0.1mm] (7.2,-15.3) circle (.6mm);
\filldraw[aqua] (7.2,-15.3) circle (.5mm);
\draw[line width=0.1mm] (8.2,-15.3) circle (.6mm);
\filldraw[aqua] (8.2,-15.3) circle (.5mm);

\draw[byzant] (6.9,-14.4) node{{\small $0$}};
\draw[byzant] (6.9,-15.3) node{{\small $3$}};
\draw[byzant] (8.4,-15.3) node{{\small $3$}};

\draw[line width=0.3mm, byzant] (7.4,-5.7) -- (7.4,-6.7) -- (8.4,-6.7) -- (7.4,-5.7);

\draw[line width=0.1mm] (7.4,-5.7) circle (.6mm);
\filldraw[aqua] (7.4,-5.7) circle (.5mm);
\draw[line width=0.1mm] (7.4,-6.7) circle (.6mm);
\filldraw[aqua] (7.4,-6.7) circle (.5mm);
\draw[line width=0.1mm] (8.4,-6.7) circle (.6mm);
\filldraw[aqua] (8.4,-6.7) circle (.5mm);

\draw[byzant] (7.05,-5.5) node{{\small $3$}};
\draw[byzant] (7.05,-6.6) node{{\small $0$}};
\draw[byzant] (8.7,-6.6) node{{\small $3$}};

\draw[line width=0.3mm, byzant] (-1.5,-5.7) -- (-1.5,-6.7) -- (-.5,-6.7) -- (-1.5,-5.7);
\draw[line width=0.1mm] (-1.5,-5.7) circle (.6mm);
\filldraw[aqua] (-1.5,-5.7) circle (.5mm);
\draw[line width=0.1mm] (-1.5,-6.7) circle (.6mm);
\filldraw[aqua] (-1.5,-6.7) circle (.5mm);
\draw[line width=0.1mm] (-.5,-6.7) circle (.6mm);
\filldraw[aqua] (-.5,-6.7) circle (.5mm);
\draw[byzant] (-1.85,-5.5) node{{\small $3$}};
\draw[byzant] (-1.85,-6.7) node{{\small $3$}};
\draw[byzant] (-.25,-6.7) node{{\small $0$}};

\draw[line width=0.7mm,pgreen] (0,-7) -- (7,-7) -- (7,-14) -- (0,-7);

\draw[line width=0.3mm, pgreen]  (3,-7.5) -- (4,-7.5);
\draw[line width=0.3mm, pgreen, dotted] (3,-7.5) -- (3,-6.5) -- (4,-7.5);
\draw[line width=0.1mm] (3,-6.5) circle (.6mm);
\filldraw[aqua] (3,-6.5) circle (.5mm);
\draw[line width=0.1mm] (3,-7.5) circle (.6mm);
\filldraw[aqua] (3,-7.5) circle (.5mm);
\draw[line width=0.1mm] (4,-7.5) circle (.6mm);
\filldraw[aqua] (4,-7.5) circle (.5mm);
\draw[pgreen] (2.75,-7.5) node{{\small $1$}};
\draw[pgreen] (4.25,-7.5) node{{\small $1$}};
\draw[pgreen] (2.75,-6.2) node{{\small $4$}};

\draw[line width=0.3mm, pgreen]  (1,-8.5) -- (2,-9.5);
\draw[line width=0.3mm, pgreen, dotted] (1,-8.5) -- (1,-9.5) -- (2,-9.5);
\draw[line width=0.1mm] (1,-8.5) circle (.6mm);
\filldraw[aqua] (1,-8.5) circle (.5mm);
\draw[line width=0.1mm] (1,-9.5) circle (.6mm);
\filldraw[aqua] (1,-9.5) circle (.5mm);
\draw[line width=0.1mm] (2,-9.5) circle (.6mm);
\filldraw[aqua] (2,-9.5) circle (.5mm);
\draw[pgreen] (0.75,-8.5) node{{\small $1$}};
\draw[pgreen] (2.25,-9.7) node{{\small $1$}};
\draw[pgreen] (0.75,-9.7) node{{\small$4$}};

\draw[line width=0.3mm, pgreen]  (7.5,-8.5) -- (7.5,-9.5);
\draw[line width=0.3mm, pgreen, dotted] (7.5,-8.5) -- (8.5,-9.5) -- (7.5,-9.5);
\draw[line width=0.1mm] (7.5,-8.5) circle (.6mm);
\filldraw[aqua] (7.5,-8.5) circle (.5mm);
\draw[line width=0.1mm] (8.5,-9.5) circle (.6mm);
\filldraw[aqua] (8.5,-9.5) circle (.5mm);
\draw[line width=0.1mm] (7.5,-9.5) circle (.6mm);
\filldraw[aqua] (7.5,-9.5) circle (.5mm);
\draw[pgreen] (7.25,-8.5) node{{\small $1$}};
\draw[pgreen] (7.25,-9.7) node{{\small$1$}};
\draw[pgreen] (8.7,-9.7) node{{\small $4$}};

\draw[line width=0.7mm,brown] (0,0) -- (0,-14) -- (14,-14) -- (0,0);

\draw[line width=0.3mm, brown] (-1.5,-2.7) -- (-1.5,-3.7) -- (-.5,-3.7) -- (-1.5,-2.7);
\draw[line width=0.1mm] (-1.5,-2.7) circle (.6mm);
\filldraw[aqua] (-1.5,-2.7) circle (.5mm);
\draw[line width=0.1mm] (-1.5,-3.7) circle (.6mm);
\filldraw[aqua] (-1.5,-3.7) circle (.5mm);
\draw[line width=0.1mm] (-.5,-3.7) circle (.6mm);
\filldraw[aqua] (-.5,-3.7) circle (.5mm);
\draw[brown] (-1.85,-2.5) node{{\small $5$}};
\draw[brown] (-1.85,-3.7) node{{\small $2$}};
\draw[brown] (-.5,-4) node{{\small $-1$}};

\draw[line width=0.3mm, brown] (4.3,-2.7) -- (4.3,-3.7) -- (5.3,-3.7) -- (4.3,-2.7);
\draw[line width=0.1mm] (4.3,-2.7) circle (.6mm);
\filldraw[aqua] (4.3,-2.7) circle (.5mm);
\draw[line width=0.1mm] (4.3,-3.7) circle (.6mm);
\filldraw[aqua] (4.3,-3.7) circle (.5mm);
\draw[line width=0.1mm] (5.3,-3.7) circle (.6mm);
\filldraw[aqua] (5.3,-3.7) circle (.5mm);
\draw[brown] (3.95,-2.5) node{{\small $5$}};
\draw[brown] (4.4,-4) node{{\small $-1$}};
\draw[brown] (5.6,-3.7) node{{\small $2$}};

\draw[line width=0.3mm, brown] (-1.5,-9.7) -- (-1.5,-10.7) -- (-.5,-10.7) -- (-1.5,-9.7);
\draw[line width=0.1mm] (-1.5,-9.7) circle (.6mm);
\filldraw[aqua] (-1.5,-9.7) circle (.5mm);
\draw[line width=0.1mm] (-1.5,-10.7) circle (.6mm);
\filldraw[aqua] (-1.5,-10.7) circle (.5mm);
\draw[line width=0.1mm] (-.5,-10.7) circle (.6mm);
\filldraw[aqua] (-.5,-10.7) circle (.5mm);
\draw[brown] (-1.85,-9.5) node{{\small $2$}};
\draw[brown] (-1.85,-10.7) node{{\small $5$}};
\draw[brown] (-.5,-11) node{{\small $-1$}};

\draw[line width=0.3mm, brown] (11.5,-9.7) -- (11.5,-10.7) -- (12.5,-10.7) -- (11.5,-9.7);
\draw[line width=0.1mm] (11.5,-9.7)  circle (.6mm);
\filldraw[aqua] (11.5,-9.7)  circle (.5mm);
\draw[line width=0.1mm] (11.5,-10.7) circle (.6mm);
\filldraw[aqua] (11.5,-10.7) circle (.5mm);
\draw[line width=0.1mm] (12.5,-10.7) circle (.6mm);
\filldraw[aqua] (12.5,-10.7) circle (.5mm);
\draw[brown] (11.2,-9.5) node{{\small $2$}};
\draw[brown] (11.65,-11) node{{\small $-1$}};
\draw[brown] (12.8,-10.7) node{{\small $5$}};

\draw[line width=0.3mm, brown] (3.2,-14.3) -- (3.2,-15.3) -- (4.2,-15.3) -- (3.2,-14.3);
\draw[line width=0.1mm] (3.2,-14.3) circle (.6mm);
\filldraw[aqua] (3.2,-14.3) circle (.5mm);
\draw[line width=0.1mm] (3.2,-15.3) circle (.6mm);
\filldraw[aqua] (3.2,-15.3) circle (.5mm);
\draw[line width=0.1mm] (4.2,-15.3) circle (.6mm);
\filldraw[aqua] (4.2,-15.3) circle (.5mm);
\draw[brown] (2.75,-14.3) node{{\small $-1$}};
\draw[brown] (2.9,-15.3) node{{\small $5$}};
\draw[brown] (4.5,-15.3) node{{\small $2$}};

\draw[line width=0.3mm, brown] (10.2,-14.3) -- (10.2,-15.3) -- (11.2,-15.3) -- (10.2,-14.3);
\draw[line width=0.1mm] (10.2,-14.3) circle (.6mm);
\filldraw[aqua] (10.2,-14.3) circle (.5mm);
\draw[line width=0.1mm] (10.2,-15.3) circle (.6mm);
\filldraw[aqua] (10.2,-15.3) circle (.5mm);
\draw[line width=0.1mm] (11.2,-15.3) circle (.6mm);
\filldraw[aqua] (11.2,-15.3) circle (.5mm);
\draw[brown] (9.75,-14.3) node{{\small $-1$}};
\draw[brown] (9.9,-15.3) node{{\small $2$}};
\draw[brown] (11.5,-15.3) node{{\small $5$}};

\draw[line width=0.2mm] (0,0) circle (1mm);
\filldraw[byzant] (0,0) circle (1mm);

\draw[line width=0.2mm] (0,-7) circle (1mm);
\filldraw[byzant] (0,-7) circle (1mm);

\draw[line width=0.2mm] (7,-7) circle (1mm);
\filldraw[byzant] (7,-7) circle (1mm);

\draw[line width=0.2mm] (0,-14) circle (1mm);
\filldraw[byzant] (0,-14) circle (1mm);

\draw[line width=0.2mm] (7,-14) circle (1mm);
\filldraw[byzant] (7,-14) circle (1mm);
\draw[line width=0.2mm] (14,-14) circle (1mm);
\filldraw[byzant] (14,-14) circle (1mm);

\end{tikzpicture} 
\caption{Example of a tiling in Theorem~\ref{thm:reduction-admissible-tiling}. The graph $G$ is the dual graph of a curve with three components. The edge lengths are all equal to one. The divisor $\varD$ on the curve $\varX$ has degree 6 and its tropicalization has support in the vertices. The $G$-admissible divisor corresponding to each stratum in the tiling is given with the same color. As in Figure~\ref{fig:tiling3}, we actually depict the projection of the tiling to $\R^2$. 
} 
\label{fig:tiling-simplex}
\end{figure}

An example of such a tiling is depicted in Figure~\ref{fig:tiling-simplex}.

\smallskip
From the proof  of Theorem~\ref{thm:reduction-admissible-tiling-sccp} and that of Theorem~\ref{thm:semistability-tiling-sccp} we get the following result.

\begin{thm}\label{thm:mixed-tiling-sccp} Let $\alpha,\beta$ be two positive real numbers. Then, the collection of supermodular functions $\alpha\mu!_D+\beta\nu!_D$, for $G$-admissible divisors $D \in \cl_{\valgroup}$, is separated, closed, positive and complete.  
\end{thm}

As a corollary, we deduce the following theorem.

\begin{thm}\label{thm:mixed-tiling}  Let $\alpha,\beta$ be two positive real numbers and $n\coloneqq\alpha d+ \beta(r+1)$. For each $G$-admissible divisor $D$ in $\cl_{\valgroup}$, let $\P_{\alpha, \beta, D}$ be the polytope associated to the supermodular function $\alpha\mu!_D+\beta\nu!_D$. Then, $\P_{\alpha,\beta,D} = \alpha \P_D+\beta \Q_D$, and the collection of polytopes $\P_{\alpha, \beta, D}$ gives a tiling of $H_{n}$. 
\end{thm}

An example of a tiling given by this theorem is depicted in Figure~\ref{fig:mixed-tiling}. Note that the first assertion in Theorem~\ref{thm:mixed-tiling} follows directly from Proposition~\ref{prop:positive-combination}: the polytope associated to $\alpha\mu!_D+\beta\nu!_D$ coincides with $\alpha \P_D + \beta \Q_D$.

\medskip

The rest of this section is devoted to the proof of these theorems. 

\subsection{Induced maps for the $\ssW_h$} 

Let $\ssh_1,\ssh_2$ be two elements of $\adfun_G(\tau(\varD);\valgroup)$. Let 
\[\lambda!_2\coloneqq\min_{v\in V}(\ssh_2(v)-\ssh_1(v)) = \min_{x\in\Gamma}(\ssh_2(x)-\ssh_1(x)),\]
where the last equality comes from the domination property. So we have 
\[\ssh_2(v)-\lambda!_2 \geq \ssh_1(v)\]
for every $v\in V$, with equality on at least one 
vertex $v\in V$.

 Recall Convention~\ref{conv:section}, with which we fix a section $\lambda\mapsto a_{\lambda}$ of the valuation $\val\colon\K^{\times}\to\valgroup$, and use it to define the relative reduction maps $\red_h$. Then, 
 $\ssa_{\lambda!_2}^{-1} \varM_{h_2} \subseteq \varM_{h_1}$, and this induces a map
 \begin{align*}
 \ssphi^{h_2}_{h_1}\colon\ssW_{h_2}&\to\ssW_{h_1}\\
 \red_{h_2}(\varf) &\mapsto \red_{h_1}(\ssa_{\lambda!_2}^{-1}\varf) \quad \forall\,\, \varf\in\varM_{h_2}.
 \end{align*}
 More precisely, this map is the restriction of the endomorphism of $\VS$ given by $\varphi \in \k^V$ with 
 \[
 \varphi_v\coloneqq \red\big(\ssa_{\ssh_1(v)}^{-1} \ssa_{\lambda!_2}^{-1} \ssa_{\ssh_2(v)}\big) =\begin{cases} 1 & \textrm{if } v\in \ssI_2\\
 0 & \textrm{otherwise}
 \end{cases}
 \]
where $\ssI_2 \subset V$ is the set of all $v\in V$ where $\ssh_2-\ssh_1$ takes its minimum, i.e., 
\[
\ssI_2\coloneqq \textrm{argmin}(\ssh_2-\ssh_1)\cap V = \left\{v\in V \st \ssh_2(v)-\ssh_1(v) =\lambda!_2 \right\}.
\] 
Then, we have
\[
\ssphi^{h_2}_{h_1}=\iota!_{I_2}\circ \proj{\ssI_2}\rest{\ssW_{\ssh_2}}\colon\ssW_{h_2}\longrightarrow\ssW_{h_1}.
\]
By symmetry, we get a map
\[
\ssphi_{h_2}^{h_1}\colon \ssW_{h_1}\to \ssW_{h_2},
\]
induced by the element $\ssa_{\lambda_1}\in \K$, where
\[
\lambda!_1 \coloneqq \min_{v\in V}(\ssh_1(v)-\ssh_2(v)) = \min_{x\in\Gamma}(\ssh_1(x)-\ssh_2(x)).
\]
Putting
\[
\ssI_1\coloneqq \textrm{argmin}(\ssh_1-\ssh_1)\cap V = \left\{v\in V \st \ssh_1(v)-\ssh_2(v) =\lambda!_1 \right\},
\]
we have that 
\[
\ssphi^{h_1}_{h_2}=\iota!_{I_1}\circ \proj{\ssI_1}\rest{\ssW_{h_1}}\colon\ssW_{h_1}\longrightarrow\ssW_{h_2}.
\]

If we were to change the section of $\val$, all the $\ssW_h$ would change in tandem to other subspaces of $\VS$ in the same orbit of the original $\ssW_h$ by the action of $\prod_{v\in V}\k^{\times}$, in such a way that $\ssphi^{h_2}_{h_1}$ would be the restriction of the same endomorphism of $\VS$ for each $\ssh_1,\ssh_2\in\adfun_G(\tau(\varD);\valgroup)$.

\begin{prop}\label{prop:comp_zero} Notation as above, 
$\ker(\ssphi_{h_1}^{h_2}) = \ssW_{h_2}^{\ssI_2^c}$ and $\image(\ssphi_{h_1}^{h_2})\subseteq
\ssW_{h_1}^{\ssI_2}$. Also, if $\ssh_2 -\ssh_1$ is not constant, then $\ssphi_{h_2}^{h_1}\circ \ssphi_{h_1}^{h_2} =0$.
\end{prop}

\begin{proof} Both assertions in the first statement follow from the fact, observed above, that $\ssphi^{h_2}_{h_1}$ is induced by an element $\varphi\in\k^V$ whose value is zero at $v\in V$ if and only if $v\in \ssI_2^c$. 

As for the second statement, since $\ssh_2-\ssh_1$ is not constant, we have $\ssI_2 \cap \ssI_1 = \emptyset,$ and so $\ssI_2\subseteq \ssI_1^c$. By symmetry, $\ker(\ssphi_{h_2}^{h_1}) = \ssW_{h_1}^{\ssI_1^c}$. Therefore,
\[
\image(\ssphi_{h_1}^{h_2})\subseteq
\ssW_{h_1}^{\ssI_2}\subseteq \ssW_{h_1}^{\ssI_1^c}=
\ker(\ssphi_{h_2}^{h_1}),
\]
and so $\ssphi_{h_2}^{h_1}\circ\ssphi_{h_1}^{h_2}=0$.
\end{proof}

Let $\ssD_1 = \tau(\varD) + \div(\ssh_1)$ and $\ssD_2 = \tau(\varD) + \div(\ssh_2)$.

\begin{prop}\label{cor:splitting} Notation as above, assume $\ssh_2 -\ssh_1$ is not constant. Then, the pair $(\ssI_2, \ssI_2^c)$ is a separation for $\nu!_{D_1}$ and $ \nu!_{D_2}$, which is strict if and only if $\ssW_{h_1}^{\ssI_2}\neq\image(\ssphi_{h_1}^{h_2})$. Similarly, the pair $(\ssI_1, \ssI_1^c)$ is a separation for $\nu!_{D_2}$ and $\nu!_{D_1}$, which is strict if and only if $\ssW_{h_2}^{\ssI_1}\neq\image(\ssphi_{h_2}^{h_1})$.
\end{prop}

\begin{proof} Applying the previous proposition, we infer that 
\begin{align*}
\nu!_{D_2}(\ssI_2^c) + \nu!_{D_1}(\ssI_2) &=\dim_{\k}\left(\ssW_{h_2}^{\ssI_2^c}\right) + \dim_{\k}\left(\ssW_{h_1}^{\ssI_2}\right)\\
&
 \geq\dim_{\k}\left(\ker\left(\ssphi_{h_1}^{h_2}\right)\right)+\dim_{\k}\left(\image\left(\ssphi_{h_1}^{h_2}\right)\right) =\dim_{\k}(\ssW_{h_2}) = r+1.
\end{align*}
This proves the first claim. The second follows by symmetry.
\end{proof}

\begin{prop}\label{prop:compoh} Let $h_1,h_2,h_3\in\adfun_G(\tau(\varD);\valgroup)$. If 
$\ssphi_{h_3}^{h_2}\circ\ssphi_{h_2}^{h_1}\neq 0$, then 
$\ssphi_{h_3}^{h_2}\circ\ssphi_{h_2}^{h_1}=\ssphi_{h_3}^{h_1}$.
\end{prop}

\begin{proof} Let 
\[\lambda!_{i,j}\coloneqq\min_{v\in V}(\ssh_i(v)-\ssh_j(v))
\quad\text{and}\quad
\ssI_{i,j}\coloneqq \textrm{argmin}(\ssh_i-\ssh_j)
\]
for $i,j\in\{1,2,3\}$. Then, $\lambda!_{1,3}\geq\lambda!_{1,2}+\lambda!_{2,3}$, with equality if and only if $\ssI_{1,2}\cap\ssI_{2,3}\neq\emptyset$, in which case $\ssI_{1,2}\cap\ssI_{2,3}=\ssI_{1,3}$. Now, if $\ssI_{1,2}\cap\ssI_{2,3}=\emptyset$, then $\ssphi_{h_3}^{h_2}\circ \ssphi_{h_2}^{h_1}= 0$. Otherwise, $\ssI_{1,2}\cap\ssI_{2,3}=\ssI_{1,3}$, and hence 
$\ssphi_{h_3}^{h_2}\circ \ssphi_{h_2}^{h_1}=\ssphi_{h_3}^{h_1}$.
\end{proof}

\subsection{$\C!_{(\varX,G)}(\varD, \varH)$ is separated}

\begin{prop}\label{prop:separation_simplex} The collection $\C!_{\varX,G}(\varD, \varH)$ is separated.
\end{prop} 

\begin{proof} Let $\ssh_1,\ssh_2\in\adfun_G(\tau(\varD);\valgroup)$ and put $\ssD_i\coloneqq \tau(\varD)+\div(\ssh_i)$ and $\nu!_i\coloneqq \nu!_{D_1}$ for $i=1,2$. We have
$\ssD_2=\ssD_1+\div(\ssh_2-\ssh_1)$. Let $\epsilon!_1 < \epsilon!_2< \dots< \epsilon!_r$ be the increasing sequence of values of $\ssh_2-\ssh_1$ on $V$ taken on subsets $\sspi_1, \dots, \sspi_r$, respectively, so that $\pi = (\sspi_1, \dots, \sspi_r)$ is a nontrivial ordered partition of $V$. 
Assume $\ssD_1\ne\ssD_2$, so $r>1$.

Let $I=\sspi_1$ and $J = \sspi_1^c$, and consider the ordered bipartition $\varpi=(I,J)$. Applying Proposition~\ref{cor:splitting}, we deduce that 
$\varpi$ is a separation for $\nu!_{1}$ and $\nu!_{2}$. 

If this separation is trivial, then we must have, using Proposition~\ref{cor:splitting} again, that
 \begin{itemize}
  \item $\nu!_{1}$ and $\nu!_{2}$ are $\varpi$-split, and
  \item $\image(\ssphi_{h_1}^{h_2})=\ssW_{h_1}^I$.
 \end{itemize}
The first property is equivalent to $\iota!_I\circ \proj{I}$ inducing isomorphisms $\ssW_{h_i}^I \to \ssW_{h_i,I}$ for $i=1,2$, which lead to decompositions: 
 \[
 \ssW_{h_1} = \ssW_{h_1}^{I} \oplus \ssW_{h_1}^{J}
 \quad \textrm{and} \quad
 \ssW_{h_2} = \ssW_{h_2}^{I} \oplus \ssW_{h_2}^{J}.
 \]
As for the second, since $\ker(\ssphi_{h_1}^{h_2}) = \ssW_{h_2}^{J}$ and $\ssW_{h_1}^I=\ssW_{h_1,I}$, it is equivalent to $\ssphi_{h_1}^{h_2}$ inducing an isomorphism from $\ssW_{h_2,I}$ to $\ssW_{h_1,I}$, or equivalently, an isomorphism $\ssW_{h_2}^I\to\ssW_{h_1}^I$.

 Let $\hat f\in\Rat(\Gamma)$ be the $G$-admissible extension of
 $f\coloneqq \epsilon!_2 \one!_{J}+\epsilon!_1 \one!_{I}$ with respect to $\ssD_1$, and put $D' \coloneqq  \ssD_1+\div(\hat f)$. We have
$D'=\tau(\varD)+\div(\ssh')$ where $\ssh'\coloneqq \ssh_1+\hat f$. Since the $\epsilon!_j$ are in $\valgroup$, we get that also $\ssh'\in\adfun_G(\tau(\varD);\valgroup)$. 
 Notice that the set of 
 vertices of $V$ where $\ssh'-\ssh_1$ (resp.~$\ssh_2-\ssh'$) takes minimum value is $\sspi_1$ (resp.~$\sspi_1 \cup \sspi_2$). 

It follows that the maps $\ssphi_{\ssh'}^{h_2}$ and $\ssphi_{h_1}^{\ssh'}$ induce injections $\ssW_{h_2}^{I}\to\ssW_{\ssh'}^{I}$ and $\ssW_{\ssh'}^{I}\to\ssW_{h_1}^{I}$. Also, by Proposition~\ref{prop:compoh}, the composition of these injections is the map $\ssW_{h_2}^{I} \to\ssW_{h_1}^I$ induced by $\ssphi_{h_1}^{h_2}$, which is an isomorphism. We infer that both injections are isomorphisms. Now, since $\ssW_{h_1}^I=\ssW_{h_1,I}$, the map $\ssphi_{h_1}^{\ssh'}$ induces as well an isomorphism $\ssW_{\ssh',I} \to\ssW_{h_1,I}$, and  we have $\ssW_{\ssh'}^I=\ssW_{\ssh',I}$. But then the spaces $\ssW_{\ssh'}^J$ and $\ssW_{h_1}^J$ have the same dimension, and since $\ssphi_{\ssh'}^{h_1}$ induces an injection 
 $\ssW_{h_1}^J\to\ssW_{\ssh'}^J$, that injection is an isomorphism too. We conclude that 
\[
 \ssW_{\ssh'} = \ssW_{\ssh'}^{I} \oplus \ssW_{\ssh'}^{J}
 = \ssW_{h_1}^{I} \oplus \ssW_{h_1}^{J}.
 \]
It follows from this equality that $\nu!_{D'} =\nu!_{1}$.

On the other hand, $\ssh_2-\ssh'$ takes only $r-1$ values. If $r=2$, then $\ssD_2=\ssD'$, whence $\nu!_{2}=\nu!_{\ssD'}=\nu!_{1}$. Thus, if $\nu!_1\neq\nu!_2$ and $r=2$,  the separation $\varpi$ is nontrivial. Arguing now by induction on $r$, if $r>2$ and $\nu!_1\neq\nu!_2$, then we have $\nu!_2\neq\nu!_{\ssD'}$, and thus by induction there is a nontrivial separation for 
$\nu!_{\ssD'}$ and $\nu!_2$. Since $\nu!_{\ssD'}=\nu!_1$, we conclude.
\end{proof}

\subsection{Exact pairs} Consider two functions $\ssh_1, \ssh_2 \in \adfun_G(\tau(\varD);\valgroup)$. In this section we prove Proposition~\ref{prop:vanishing}, which roughly speaking shows that if $\ssh_1$ and $\ssh_2$ are distinct but close to each other, then the two spaces $\ssW_{h_1}$ and $\ssW_{h_2}$ are in close relation.

More precisely, we say that $\ssW_{h_1}$ and $\ssW_{h_2}$ form an \emph{exact pair} if we have 
\[\ker(\ssphi_{h_1}^{h_2}) = \image(\ssphi_{h_2}^{h_1}) \quad \textrm{and}  \quad \ker(\ssphi_{h_2}^{h_1}) = \image(\ssphi_{h_1}^{h_2}),\]
that is, if the sequence 
\[\ssW_{h_1} \to \ssW_{h_2} \to \ssW_{h_1} \to \ssW_{h_2}\]
induced by the maps $\ssphi_{h_2}^{h_1}$ and $\ssphi_{h_1}^{h_2}$ is exact. In this case, we say $(\ssW_{h_1}, \ssW_{h_2})$ is exact. 

\begin{prop} \label{prop:exact-one-way} A pair $(\ssW_{h_1}, \ssW_{h_2})$ is exact if and only if $\ker(\ssphi_{h_1}^{h_2}) = \image(\ssphi_{h_2}^{h_1})$.
\end{prop}
\begin{proof} One direction is trivial. Assume $\ker(\ssphi_{h_1}^{h_2}) = \image(\ssphi_{h_2}^{h_1})$. In any case, Proposition~\ref{prop:comp_zero} yields 
$\image(\ssphi_{h_1}^{h_2}) \subseteq  \ker(\ssphi_{h_2}^{h_1})$. But now
\[
\dim\left( \image(\ssphi_{h_1}^{h_2}))\right) = r+1 - \dim\left(\ker(\ssphi_{h_1}^{h_2})\right) =  r+1 - \dim\left( \image(\ssphi_{h_2}^{h_1}) \right) = \dim\left(\ker(\ssphi_{h_2}^{h_1}) \right),
\]
from which the proposition follows.
\end{proof}

Notice that if $(\ssW_{h_1}, \ssW_{h_2})$ is exact, then $\ssh_1-\ssh_2$ is not constant. Moreover, denoting by $\ssI_1$ (resp.~$\ssI_2$) the subset of $V$ where $\ssh_1-\ssh_2$ (resp.~$\ssh_2-\ssh_1$) takes minimum value, it follows from Proposition~\ref{prop:comp_zero} that $\ssphi_{h_2}^{h_1}$ (resp.~$\ssphi_{h_1}^{h_2}$) factors through an isomorphism 
$\ssW_{h_1,I_1}\cong\ssW_{h_2}^{I_1}$ (resp.~
$\ssW_{h_2,I_2}\cong\ssW_{h_1}^{I_2}$). 
Clearly, $I_1\cap I_2=\emptyset$.

\smallskip

The following proposition gives sufficient conditions for exactness.

\begin{prop}\label{prop:vanishing} Let $\ssh_1,\ssh_2 \in \adfun_G(\tau(\varD);\valgroup)$. Assume that the difference $\ssh_2-\ssh_1$ takes exactly two values on $V$. Let $d$ be the positive difference between these two values, and let $F\subseteq V$ be the subset where $\ssh_2-\ssh_1$ takes the maximum value. Then, if $d$ is small enough, we have 
\[\ssW_{h_2} = \ssW_{h_2}^{F} \oplus \ssW_{h_2}^{F^c} \quad \textrm{and the pair }
(\ssW_{h_1}, \ssW_{h_2})\textrm{ is exact}.\]
\end{prop}

The proof gives an explicit description in terms of a constant $c>0$, defined in terms of $\ssh_1$ and the reduction map $\red_{h_1}$, on how small $d$ should be in the proposition. 

We introduce the following notation before going through the proof. Given an element $h\in \adfun_G(\tau(\varD);\valgroup)$ and a subset $F\subseteq V$, we define 
\[
\varM_{h}^{F} \coloneqq \left\{\varf \in \varM_{h} \st \trop(\varf)(v) > \ssh(v) \,\, \textrm{for all} \,\, v\in \ssF^c \right\}.
\]
This is an $\varR$-submodule of $\varM_{h}$. Obviously, we have $\varm \varM_{h} \subseteq \varM_{h}^{F}$. Moreover, $\rquot{\varM_{h}^{F}}{\varm \varM_{h}}$ is naturally isomorphic to $\ssW_{h}^{F}$.

\begin{proof} Applying Proposition~\ref{prop:comp_zero}, we see that  $\ssphi_{h_2}^{h_1}$ factors through an 
injection $\ssW_{h_1,F}\to\ssW_{h_2}^{F}$, and $\ssphi_{h_1}^{h_2}$ factors through an injection 
$\ssW_{h_2,F^c}\to\ssW_{h_1}^{F^c}$.  By Proposition~\ref{prop:exact-one-way}, the former is an isomorphism if and only if the latter is. Also, the composition of the two maps
\[\ssW_{h_2}^{F^c}\hookrightarrow\ssW_{h_2,F^c} \to\ssW_{h_1}^{F^c}\] 
is the map $\ssW_{h_2}^{F^c}\to\ssW_{h_1}^{F^c}$ induced by $\ssphi_{h_1}^{h_2}$. It is thus enough to 
show that this induced map is an isomorphism, provided that $d$ is small enough.

As observed before the proof, the reduction map $\red_{h_1}$ takes $\varM_{h_1}^{F^c}$ onto $\ssW_{h_1}^{F^c}$. Choose elements $\varf_1, \dots, \varf_p \in \varM_{h_1}^{F^c}$ whose reductions $\red_{h_1}(\varf_j)$ generate $\ssW_{h_1}^{F^c}$. Consider the $\varR$-submodule $\varN$ of $\varM_{h_1}$ generated by $\varf_1, \dots, \varf_p$.  

 Define 
\begin{equation}\label{eq:exactness}
c\coloneqq \min_{j\in [p]} \left(\min_{v\in F}\Bigl(\trop(\varf_j)(v) -\ssh_1(v)\Bigr)\right).
\end{equation}
Then, $c$ is positive and belongs to the value group $\valgroup$. Also, each $\varf  \in \varN$ can be written as 
$\varf = \ssa_1\varf_1+\dots+ \ssa_p\varf_p$ for certain $\ssa_j\in \varR$, whence the ultrametric triangle inequality yields $\trop(\varf)(v) -\ssh_1(v)\geq c$ for every $v\in F$.

After summing a constant, if necessary, we may assume the difference $\ssh_1-\ssh_2$ takes maximum value zero. Then, for each $\varf \in \varM_{h_2}$, we have  $\trop(\varf) \geq \ssh_1$, and so $\varM_{h_2} \subseteq \varM_{h_1}$.  Moreover, $\trop(\varf)(v) \geq \ssh_2(v) > \ssh_1(v)$ for all $v\in F$. That is, $\varM_{h_2} \subseteq \varM_{h_1}^{F^c}$. 

Assume $d<c$. Then $\varN \subseteq \varM_{h_2}^{F^c}$.  
Since $\ssh_1$ and $\ssh_2$ agree on $F$, we have 
$\red_{h_1}(\varf_j)=\red_{h_2}(\varf_j)$ for each $j\in[p]$. Since $\ssW_{h_1}^{F^c}$ is generated by the reductions
$\red_{h_1}(\varf_j)$ and $\ssW_{h_2}^{F^c}$ contains the 
$\red_{h_2}(\varf_j)$, we have that $\ssphi_{h_1}^{h_2}$ induces an isomorphism 
$\ssW_{h_2}^{F^c}\to\ssW_{h_1}^{F^c}$, as required. 
 \end{proof}

\subsection{$\C!_{(\varX,G)}(\varD, \varH)$ is closed} \label{sec:closed-admissible-lls} 

\begin{prop} \label{prop:closed-lls}The collection  $\C!_{(\varX,G)}(\varD, \varH)$ is closed.
\end{prop}

\begin{proof} Consider a $G$-admissible function $h\in \adfun_G(\tau(\varD);\valgroup)$, and let $D\coloneqq\tau(\varD)+\div(h)$. Let $\nu=\nu!_D$. Let $\pi=(\ssF^c,F)$ be an ordered bipartition of $V$. By Proposition~\ref{refineV}, it is enough to show that $\nu!_{\pi}$ belongs to $\C!_{(\varX,G)}(\varD, \varH)$.

Let $\hat f\in\Rat(\Gamma)$ be the $G$-admissible extension of $f\coloneqq \epsilon \one!_{F}$ with respect to $D$, and put 
$D'\coloneqq D+\div(\hat f)$. Then, $D'=\tau(\varD)+\div(h')$ for $h'\coloneqq h+\hat f$. Since $\epsilon\in\valgroup$, we have that also 
$h'\in\adfun_G(\tau(\varD);\valgroup)$. Let $\nu'\coloneqq \nu!_{D'}$. We claim that $\nu!_{\pi} = \nu'$.

Indeed, by Proposition~\ref{prop:vanishing}, we know that $\nu'$ is $\pi$-split. In addition, and because of that, $\phi_{h'}^{h}$ and $\ssphi_{h}^{h'}$ induce isomorphisms $\ssW_{h,F}\to\ssW_{h',F}$ and 
$\ssW_{h'}^{F^c}\to\ssW_{h}^{F^c}$. It follows from the second isomorphism that
$\nu'(I)=\nu(I)$ for each $I\subseteq \ssF^c$, and from the first that $\ssub{(\nu')}!^*(I)=\nu!^*(I)$ for each 
$I\subseteq F$. By adjunction, this latter property is equivalent to $\nu'(J)=\nu(J)$ for each $J\subseteq V$ with $\ssF^c\subseteq J$. We infer that
\[
\nu!_{\pi}(I)=\nu(I\cup \ssF^c)+\nu(I\cap \ssF^c)-\nu(\ssF^c)=\nu'(I\cup \ssF^c)+\nu'(I\cap \ssF^c)-\nu'(\ssF^c)=\nu!_{\pi}'(I)=\nu'(I)
\]
for each $I\subseteq V$, and the proposition follows.
\end{proof}

\subsection{Threshold of $\pi$-splitting} Let  $\nu\coloneqq \nu!_D$ be an element of $\C!_{(\varX,G)}(\varD, \varH)$ 
such that $\Q_{D}$ intersects $\ring\Delta_{r+1}$. Let 
$\pi=(I,J)$ be a bipartition of $V$ such that $\nu!_\pi = \nu$. For each $t\in\valgroup$, let 
$\ssD_t\coloneqq D+\div(t\one!_J,D)$.

In this section we prove the following result.

\begin{prop}\label{prop:threshold} There exists a largest positive real $c \in (0,+\infty)$ such that $\nu!_{D_t} =\nu$ for each $t\in [0,c)\cap \valgroup$. Moreover, $c\in\valgroup$.
\end{prop}

Notice that, for $t$ small, $\ssD_t$ is the $G$-admissible chip-firing induced by $I$. In analogy with Definition~\ref{defi:minimum-length}, we define the following quantity associated to $D$ and $\pi=(I, J)$; see Remark~\ref{rmk:mlD}.

\begin{defi}\label{defi:threshold} Notation as above, we denote the element of $\valgroup$ obtained from Proposition~\ref{prop:threshold} by $\ts_{D}(I)$ and call it the \emph{threshold of $\pi$-splitting of $\nu!_D$.}
\end{defi}

The rest of this section is devoted to the proof of Proposition~\ref{prop:threshold}. Let $h\in \adfun_G(\tau(\varD);\valgroup)$ be such that 
$D=\tau(\varD)+\div(h)$. For each positive real number $t$ in the value group $\valgroup$, let $\ssfhat_{t}$ be the $G$-admissible extension of $t\one!_J$ with respect to $D$. Let $\ssh_t\coloneqq h+\ssub{\hat f}!_{t}$. Note that we have $\ssh_t\in\adfun_G(\tau(\varD);\valgroup)$ and  $\ssD_t= \tau(\varD)+\div(\ssh_t)$.

\begin{lemma}\label{lem:open-S}
Let $S\subset\valgroup \subset \R$ be the subset consisting of all real numbers $t\in \valgroup$ such that $\nu!_{D_t} =\nu$.  There exists an open subset $U$ of $\R$ such that $S = U \cap \valgroup$.
\end{lemma}
\begin{proof}
Since $\nu$ is $\pi$-split, we have 
\[
\ssW_{h} = \ssW_{h}^{I} \oplus \ssW_{h}^{J}.
\]
Since $\Q_{D}$ intersects $\ring\Delta_{r+1}$, we must have 
$\ssW_{h,I}\neq 0$. 

We apply Proposition~\ref{prop:vanishing} to the pair of elements $h$ and $\ssh_t$ in $\adfun_G(\tau(\varD);\valgroup)$. It yields $c\in\R$ such that for each $t\in(0,c)\cap \valgroup$, we have
\[
\ssW_{h_t}=\ssW_{h_t}^I\oplus \ssW_{h_t}^{J}=\ssW_{h,I}\oplus \ssW_{h_t,J}=
\ssW_{h}^I\oplus \ssW_{h}^{J}=\ssW_{h}.
\]
Thus, $\nu!_{D_t} =\nu$ for every $t\in [0, c)\cap \valgroup$.  Applying the same argument to $\nu$ and the bipartition $\sspi^c=(J,I)$ instead of $\pi$ yields a negative number $c'\in\R$ such $\nu!_{D_t} =\nu$ for every $t\in(c',0]\cap \valgroup$.

Notice that, since $\nu!_{D_t}=\nu$, we have that 
$\nu!_{D_t}$ is $\pi$-split and $\Q_{D_t}$ intersects $\ring\Delta_{r+1}$ for each $t\in S$. The above reasoning can be thus applied to each point of $S$, showing the existence of an open $U \subset \R$ such that $S = U \cap \valgroup$.   
\end{proof}

  We will need the following lemma.

\begin{lemma}\label{lem:vanishing-high-value} Let $\pi=(I,J)$ be a bipartition of $V$. Let $D\in\stable!_G(\cl_{\valgroup})$ and 
$\ssD_t\coloneqq D+\div(t\one!_J;D)$ for each $t\in\valgroup$. Then, $\ssW_{D_t,I}=0$ for large enough $t$.
\end{lemma}

\begin{proof} Let $h\in \adfun_G(\tau(\varD);\valgroup)$ such that $D=\tau(\varD)+\div(h)$. For each $t\in\valgroup$, let $\ssh_t\in \adfun_G(\tau(\varD);\valgroup)$ such that $\ssh_t-h$ is the $G$-admissible extension of $t\one!_J$ with respect to $D$. By definition, $\ssW_{D_t}=\ssW_{h_t}$ is the reduction of $\varM_{h_t}$. We claim that, if $t$ is large enough, for every $\varf\in \varM_{h_t}$, we have $\ssub{(\red_{h_t}(\varf))}!_v=0$ for all $v\in I$.

 By definition, $\varM_{h_t}$ consists of those $\varf \in \varH$ that verify $\trop(\varf) \geq \ssh_t$. By the definition of the reduction map $\red_{h_t}$, we need to prove that $\trop(\varf)(v) >\ssh_t(v) = h(v)$ for every $\varf \in \varM_{h_t}$ and $v\in I$, provided that $t$ is large enough.
 
  Let $\varf\in \varH$. Using the equality $\tau(\varD) + \div(\trop(\varf)) = \tau(\varD + \div(\varf))$, we get 
 \[\div(\trop(\varf)) + \tau(\varD) \geq 0.\]
 By~\cite[Lemma~1.8]{GK08}, this implies that the slopes taken by $\trop(\varf)$ are all bounded in absolute value. Using connectivity and compactness of $\Gamma$, this implies that there exists a constant $C>0$ such that 
 $\trop(\varf)(u)-\trop(\varf)(v)<C$ for each $u,v\in V$. Choose $t$ larger than $C+h(v)-h(u)$ for every $v\in I$ and $u\in J$. Then, for each $\varf\in\varH$ and $v\in I$, if $\trop(\varf)(v)=h(v)$, we get 
 \[
 \trop(\varf)(u)< \trop(\varf)(v)+C=h(v)+C<t+h(u)=\ssh_t(u)
 \]
  for each $u\in J$, that is, $\varf \notin \varM_{h_t}$.  This proves the claim and the lemma follows.
\end{proof}

\begin{proof}[Proof of Proposition~\ref{prop:threshold}] 
We start by proving the existence part of the proposition. Since $\Q_D\cap\ring\Delta_{r+1}\neq\emptyset$, we have $\ssW_{\ssD,I}\neq 0$, and thus also $\ssW_{D_t,I}\neq 0$ for each $t\in S$, as $\nu!_{D_t} =\nu$. 
We apply now Lemma~\ref{lem:vanishing-high-value} and infer that $S$ does not contain the full set  $[0, +\infty)\cap \valgroup$. In particular, there exists a largest $c\in (0,+\infty)$ such that for each $t\in [0,c)\cap \valgroup$, we have $\nu!_{D_t} =\nu$. This proves that $\ts_D(I)$ exists.

\smallskip

 We show now that $\ts_D(I)$ belongs to $\valgroup$. Consider the set $\trop(\varH)$ consisting of the tropicalizations of all functions in $\varH$, i.e.,
\[\trop(\varH)\coloneq \left\{\trop(\varf)\, \st \,\varf\in \varH\right\}.\]
The set $\trop(\varH)$ is closed in the sense that any pointwise limit $\ssf\colon \Gamma \to \R$ of functions $\ssf_1, \ssf_2, \dots$ in $\trop(\varH)$ belongs to $\trop(\varH)$, see~\cite[Theorem 8.3]{AG22}.

 Consider a sequence $(\sst_n)_{n\in\N}$ of positive real numbers in $\valgroup$ converging to 
$\ts_D(I)$ from below. We have $t_n\in S$ by the definition of $\ts_D(I)$. Applying Proposition~\ref{prop:vanishing}, there is $\ssc_n\in\R_{>0}$ such that $[\sst_n, \sst_n+\ssc_n)\cap \valgroup \subset S$ for each $n\in\N$ (see Lemma~\ref{lem:open-S}). In particular, $\sst_n+\ssc_n\leq \ts_D(I)$. 
Moreover, as the proof of Proposition~\ref{prop:vanishing} shows, we may choose $c_n$ such that Equation~\eqref{eq:exactness} holds. In particular, there is a rational function $\varf_n \in \varM_{h_{\sst_n}}^{I}$ verifying $\trop(\varf_n)\rest{I} = \ssh_{t_n}\rest{I}$ such that
\[
\ssc_n \coloneqq \min_{v\in J}\bigl( \trop(\varf_n)(v) - \ssh_{t_n}(v)\bigr)= \min_{v\in J} \bigl(\trop(\varf_n)(v) - \ssh(v)-\sst_n\bigr).
\]
Since $\trop(\varf_n)$ coincides with $\ssh$ on $I$, and the slopes of $\trop(\varf_n)$ are bounded (see proof of Lemma~\ref{lem:vanishing-high-value}), we infer, by applying Arzel\`a--Ascoli Theorem, that passing to a subsequence if necessary, the sequence $\bigr(\trop(\varf_n)\bigl)_{n\in \N}$ converges uniformly. Since $\trop(\varH)$ is closed,  the sequence converges to $\trop(\varf)$ for some element $\varf\in \varH$. To conclude, note that
\begin{align*}
\ts_D(I)&=\lim_{n\to \infty} \sst_n = \lim_{n\to \infty} (\sst_n+\ssc_n)= \lim_{n\to \infty}\left(\min_{v\in J} \bigl(\trop(\varf_n)(v) - h(v) \bigr)\right)\\
&=\min_{v\in J} \bigl(\trop(\varf)(v) - h(v)\bigr) \in \valgroup.\qedhere
\end{align*}
\end{proof}

\subsection{$\C!_{(\varX,G)}(\varD, \varH)$ is relatively complete}\label{sec:complete-admissible-lls}  

\begin{prop} \label{prop:complete-admissible-lls} The collection  $\C!_{(\varX,G)}(\varD, \varH)$ is complete relative to $\ring\Delta_{r+1}$, the interior of $\Delta_{r+1}$.
\end{prop}

\begin{proof} Since $\C!_{(\varX,G)}(\varD, \varH)$ is closed by Proposition~\ref{prop:closed-lls}, we only need to consider an element $\nu\coloneqq \nu!_D$ of $\C!_{(\varX,G)}(\varD, \varH)$ and 
a bipartition $\pi=(F,\ssF^c)$ of $V$ such that $\nu!_\pi = \nu$ and $\Q_{D}$ intersects $\ring\Delta_{r+1}$, and to show the 
existence of $\nu'$ in $\C!_{(\varX,G)}(\varD, \varH)$ satisfying $\nu' \neq \nu$ with $\pi$-splitting $\nu!_{\pi}' =\nu$. 

Let $h\in \adfun_G(\tau(\varD);\valgroup)$ be such that 
$D=\tau(\varD)+\div(h)$. For each positive real number $t$ in the value group $\valgroup$, let $\ssub{\hat f}!_{t}$ be the $G$-admissible extension of $\ssub{f}!_{t}\coloneqq t\one!_{\ssF}$ with respect to $D$. Let
$\ssh_t\coloneqq h+\ssub{\hat f}!_{t}$ and $\ssD_t\coloneqq\tau(\varD)+\div(\ssh_t)$. 

We apply Proposition~\ref{prop:threshold} to $\nu=\nu!_D$ and $\sspi^c=(\ssF^c, F)$. Let $b\coloneq \ts_D(\ssF^c)$ be the threshold of $\sspi^c$-splitting of $\nu$. We have $b\in\valgroup$. Let  $\ssh'\coloneqq \ssh_b$ and $\nu'\coloneqq \nu!_{D_b}$. Then, we must have $\nu'\neq\nu$, as otherwise, by Lemma~\ref{lem:open-S}, there would exist an $\epsilon>0$ such that for every $t\in [b,b+\epsilon]\cap \valgroup$, we would get $\nu!_{D_t}=\nu!_{D_b}=\nu$, contradicting the definition of $\ts_D(\ssF^c)$. 

Apply now Proposition~\ref{prop:vanishing} to the pair of elements $h', \ssh_{b-t}$ for a small enough positive number $t \in \valgroup$. Since $t$ is small, $\ssW_{h_{b-t}} = \ssW_{h}$. By that proposition, 
\[
\ssW_{h}=\ssW_{h,F} \oplus \ssW_{h}^{F^c}=
\ssW_{h'}^F\oplus\ssW_{h',F^c},
\]
whence $\nu!_{\pi}' = \nu$. We have proved the existence of $\nu' \neq \nu$ with $\nu!_\pi' = \nu$, as required.
\end{proof}

\subsection{Proof of Theorems~\ref{thm:reduction-admissible-tiling-sccp}~and~\ref{thm:reduction-admissible-tiling}} We have verified that the collection $\C=\C!_{(\varX,G)}(\varD, \varH)$ is separated, closed, and complete for $\ring\Delta_{r+1}$, our Propositions~\ref{prop:separation_simplex},~\ref{prop:closed-lls}~and~\ref{prop:complete-admissible-lls}. Since the supermodular functions in $\C!_{(\varX,G)}(\varD, \varH)$ take integer values, the spread of the simple ones is at least $1$. This proves that 
$\C!_{(\varX,G)}(\varD, \varH)$ is positive, finishing the proof of Theorem~\ref{thm:reduction-admissible-tiling-sccp}. We thus infer from Theorems~\ref{thm:tiling}~and~\ref{thm:tiling2} that the corresponding polytopes provide a tiling of $\P_{\C}$ containing the interior of $\Delta_{r+1}$. Since $\P_{\C}$ is closed and contained in $\Delta_{r+1}$, we conclude. \qed


 \subsection{Proof of Theorems~\ref{thm:mixed-tiling-sccp}~and~\ref{thm:mixed-tiling}} Consider the collection $\C_{\alpha, \beta}$ given by the supermodular functions $\alpha\mu!_{D}+\beta\nu!_{D}$ for 
 $D\in\stable!_G(\cl_{\valgroup})$, where $\cl$ is the class of $\tau(\varD)$ in $\Pic^d(\Gamma)$. 

To show that $\C_{\alpha, \beta}$ is separated, 
let $\ssD_1$ and $\ssD_2$ be distinct divisors in $\stable!_G(\cl_{\valgroup})$. Put $\mu!_i\coloneqq\mu!_{D_i}$ and $\nu!_i\coloneqq\nu!_{D_i}$, as well as $\ssrho_i\coloneqq\alpha\mu!_i+\beta\nu!_i$ for $i=1,2$. 
Observe that the proofs of Proposition~\ref{prop:Cseparated} and Proposition~\ref{prop:separation_simplex} follow the same argument. They start by producing a bipartition $\varpi$ of $V$ that is a separation both for $\mu!_1$ and $\mu!_2$, and for $\nu!_1$ and $\nu!_2$, whence for $\ssrho_1$ and $\ssrho_2$. If the last separation is trivial, so are the first two, and both proofs show that a particular $D'\in\stable!_G(\cl_{\valgroup})$ satisfies $\mu!_{D'}=\mu!_1$ and $\nu!_{D'}=\nu!_1$, and hence $\ssrho'=\ssrho_1$ where 
$\ssrho'\coloneqq\alpha\mu!_{D'}+\beta\nu!_{D'}$.
The proof concludes now as in the propositions, by observing that if 
$\ssrho_1\neq\ssrho_2$, then either $\ssrho'\neq\ssrho_1$, and hence 
$\varpi$ is a nontrivial separation for $\ssrho_1$ and $\ssrho_2$, or $\ssrho'=\ssrho_1$, whence $\ssrho'\neq\ssrho_2$, but 
$D'$ is closer to $\ssD_2$ than $\ssD_1$, and hence there is by induction a nontrivial separation for 
$\ssrho'$ and $\ssrho_2$, and thus for $\ssrho_1$ and $\ssrho_2$.

To show that the collection $\C_{\alpha, \beta}$ is closed, let 
$h\in\adfun_G(\tau(\varD);\valgroup)$, and let $\pi=(\ssF^c,F)$ be a bipartition of $V$. Let $D\coloneqq \tau(\varD)+\div(h)$, and $\mu\coloneqq\mu!_{D}$ and $\nu\coloneqq\nu!_{D}$. 
Then, the proofs of  Proposition~\ref{prop:Cclosed} and Proposition~\ref{prop:closed-lls} show that for a small positive number $\epsilon$ in the value group $\valgroup$, we have that $h'\in\adfun_G(\tau(\varD);\valgroup)$, where $h'\coloneqq h+\hat f$, for $\hat f$ the $G$-admissible extension of $f\coloneqq \epsilon\one!_F$ with respect to $D$, and that $\mu!_{\pi} = \mu!_{D'}$ and $\nu!_{\pi} = \nu!_{D'}$, where $D'\coloneqq\tau(\varD)+\div(h')$. Thus,
\[
\ssub{(\alpha \mu+\beta\nu)}!_{\pi} = \alpha \mu!_{\pi}+\beta\nu!_{\pi} = \alpha \mu!_{D'}+\beta\nu!_{D'}
\]

To show that $\C_{\alpha, \beta}$ is complete, let $h\in\adfun_G(\tau(\varD);\valgroup)$ and $D\coloneqq \tau(\varD)+\div(h)$. Put $\mu\coloneqq\mu!_D$ and $\nu\coloneqq\nu!_D$, and let 
$\ssrho\coloneqq \alpha\mu+\beta\nu$. Let $\pi=(F,\ssF^c)$ be a bipartition of $V$ such that $\ssrho_{\pi}=\ssrho$. Then, 
$\mu!_{\pi}=\mu$ and $\nu!_{\pi}=\nu$. For each positive real number $t$ in the value group $\valgroup$, let $\ssub{\hat f}!_{t}$ be the $G$-admissible extension of $\ssub{f}!_{t}\coloneqq t\one!_F$ with respect to $D$. Put $\ssh_t\coloneqq\ssh+\ssub{\hat f}!_{t}$ and $\ssD_t\coloneqq\tau(\varD)+\div(\ssh_t)$. Then, 
$\ssh_t\in\adfun_G(\tau(\varD);\valgroup)$. For each $t$, let $\mu!_t\coloneqq\mu!_{D_t}$, $\nu!_t\coloneqq\nu!_{D_t}$, and 
$\ssrho_{t}\coloneqq\alpha\mu!_t+\beta\mu!_t$.

Let $b\coloneqq \min\bigl(\ml_D(F^c), \ts_D(\ssF^c)\bigr)$, with $\ml_D(\ssF^c)$ and $\ts_D(\ssF^c)$ introduced in Definitions~\ref{defi:minimum-length} and~\ref{defi:threshold}, respectively. Since $D$ is $\valgroup$-rational, Propositions~\ref{prop:minimum-length} and~\ref{prop:threshold} yield that $b\in \valgroup$. By definition of these numbers, we have $\mu!_t=\mu$ and $\nu!_t=\nu$ for each $t\in [0,b)\cap\valgroup$, and that either $\mu!_b$ or $\nu!_b$ is not $\pi$-split.  In any case, the results show as well that $\mu!_{b,\pi}=\mu$ and $\nu!_{b,\pi}=\nu$. Thus, $\ssrho_{b,\pi}=\ssrho$, and $\ssrho_{b}$ is not $\pi$-split. We infer that $\ssrho_{b}\neq\ssrho$ and the claim follows. 

Finally, the collection $\C_{\alpha, \beta}$ is obviously positive. This proves Theorem~\ref{thm:mixed-tiling-sccp}. As before, Theorem~\ref{thm:mixed-tiling} follows from 
Theorems~\ref{thm:tiling}~and~\ref{thm:tiling2}.\qed


\section{Geometric features of the tilings} \label{sec:geometric-features}
In this section we discuss geometric properties of the tilings in Theorems~\ref{thm:semistability-tiling},~\ref{thm:reduction-admissible-tiling}, and~\ref{thm:mixed-tiling}.

\subsection{Regularity of the tilings} \label{sec:regularity}

Our first result in this section is the following.

\begin{thm}\label{thm:regularity-tiling} The tilings in Theorem~\ref{thm:semistability-tiling}, in Theorem~\ref{thm:reduction-admissible-tiling}, and in Theorem~\ref{thm:mixed-tiling} are all regular. 
\end{thm}

\begin{proof} In each case, we let 
$\stable!_G^{^\circ}(\cl_{\valgroup})\subset \stable!_G(\cl_{\valgroup})$ be the collection of $G$-admissible $\valgroup$-rational divisors $D$ whose associated polytope, be it $\P_D$, $\Q_D$ or $\P_{\alpha,\beta,D}$, is full-dimensional. To simplify the presentation, we let $\ssrho_D$ be the supermodular function and $\P_D'$ the polytope associated to $D$ in each case. Also, let $H\subset\R^V$ be the hyperplane of $\R^V$ containing all the $\P_D'$.

Each function $f\colon V\to\R$ gives rise to an $\R$-linear function on $\R^V$, denoted $(f,\cdot)$, that sends $q$ to 
\[(f,q)\coloneqq \sum_{v\in V} f(v)q(v).\] 
Fix $\ssD_0\in\stable!_G^{^\circ}(\cl_{\valgroup})$, as well as $\ssf_{0}\in\R^V$ and 
$\lambda!_{0}\in\R$. We claim that  there are $\ssf_D\in\R^V$ and $\lambda!_D\in\R$, for each $D\in\stable!_G^{\circ}(\cl_{\valgroup})$, such that the following four properties hold:
\begin{enumerate}[label=(\arabic*)]
    \item $\ssf_{D_0} =\ssf_0$ and $\lambda!_{D_0} = \lambda!_{0}$.
    \item \label{seq:2} $\ssD=\ssD_0+\div(\ssf_D-\ssf_{D_0};\ssD_0)$ for each $D\in\stable!_G^{^\circ}(\cl_{\valgroup})$.
    \item \label{seq:3} The restrictions of the affine functions 
    $(\ssf_D,\cdot)-\lambda!_D$ to $\P_D'$  for $D\in\stable!_G^{^\circ}(\cl_{\valgroup})$ glue to a global function $\Phi$ defined on the union of all the $\P_D'$ in $H$.
    \item \label{seq:4} The function $\Phi$ is convex with domains of linearity given by the $\P_{D}'$, $D\in \stable!_G^{^\circ}(\cl_{\valgroup})$. 
\end{enumerate}

To prove the claim, put $\ssf_{D_0} \coloneqq \ssf_0$ and $\lambda!_{D_0} \coloneqq \lambda!_{0}$. For each 
$D\in\stable!_G^{^\circ}(\cl_{\valgroup})$, we choose a sequence $\ssD_0,\ssD_1,\dots,\ssD_p=D$ in $\stable!_G^{^\circ}(\cl_{\valgroup})$ such that, for each $i=1,\dots,p$, the polytopes $\P_{D_{i-1}}'$ and $\P_{D_{i}}'$ intersect in a common facet. The proofs of Propositions~\ref{prop:Cclosed}~and~\ref{prop:Csimple} in one case, and those of Propositions~\ref{prop:closed-lls}~and~\ref{prop:complete-admissible-lls} in another, and finally that of Theorem~\ref{thm:mixed-tiling-sccp} in the remaining ``mixed" case, show that for each $i=1,\dots,p$, there exist a bipartition $\sspi_i=(\ssF_i^c,\ssF_i)$ of $V$ and a $b_i\in\valgroup \cap \R_{>0}$ satisfying the following three properties:
\begin{enumerate}
    \item[(a)] $\ssD_i=\ssD_{i-1}+\div(b_i\one!_{F_i};\ssD_{i-1})$.
    \item[(b)] The base polytope of $\ssub{(\ssrho_{D_{i-1}})}!_{\pi_i}$ is a facet of $\P_{D_{i-1}}'$.
    \item[(c)] The $\pi_i$-splitting of $\ssrho_{D_{i-1}}$ is the $\pi_i^c$-splitting of $\ssrho_{D_{i}}$, that is, $\ssub{(\ssrho_{D_{i-1}})}!_{\pi_i} = \ssub{(\ssrho_{D_{i}})}!_{\pi_i^c}$. 
\end{enumerate}

Property (b) implies the uniqueness of $\sspi_i$ for each $i$. Also, $b_i = \ml_{D_{i-1}}(\ssF_{i}^c)$, or $\ts_{D_{i-1}}(\ssF_i^c)$, or $\min\bigl(\ml_{D_{i-1}}(\ssF_i^c), \ts_{D_{i-1}}(\ssF^c_i)\bigr)$, depending on whether we consider $\P_D$, $\Q_D$ or $\P_{\alpha,\beta,D}$, respectively.

\smallskip
It follows from (c) above that the common facet of $\P_{D_{i-1}}'$ and $\P_{D_{i}}'$ coincides with the base polytope of $\ssub{(\ssrho_{D_{i-1}})}!_{\pi_i} = \ssub{(\ssrho_{D_{i}})}!_{\pi_i^c}$.

We define 
\begin{align} \label{eq:regularity}
\ssf_D\coloneqq \ssf_{D_0}+\sum_{i=1}^p \ssub{b}!_{i}\ssub{\one}!_{F_i}
\quad\text{and}\quad
\lambda!_D:=\lambda!_{D_0}+\sum_{i=1}^p \ssub{b}!_{i}\,\ssrho_{D_i}(\ssF_i).
\end{align}
It is clear from Property~(a) above, and Proposition~\ref{cor:sumcan}, that Property~\ref{seq:2} holds for every $i$. 

We will prove below that, though $\ssf_D$ and $\lambda!_D$ might be different for a different choice of the sequence $\ssD_i$, the restriction of $(\ssf_D,\cdot)-\lambda!_D$ to the hyperplane $H\subset\R^V$ does not depend on the choice. Assuming this for the moment, we show~\ref{seq:3} and~\ref{seq:4}. 

To prove \ref{seq:3}, consider $D, D' \in\stable!_G^{^\circ}(\cl_{\valgroup})$ such that $\P_{D}'$ and $\P_{D'}'$ share a facet. We may assume the sequence $\ssD_0, \dots, \ssD_{p-1}, \ssD_p=D$ chosen for $D$ satisfies $D'=\ssD_{p-1}$. In this case, for every point $q$ in the intersection of $\P_{D'}'$ and $\P_D'$, we have
\[
(\ssf_D-\ssf_{D'},q)=\ssub{b}!_{p}\ssq(\ssF_p)=\ssub{b}!_{p}
\ssrho_{D}(\ssF_p)=\lambda!_D-\lambda!_{D'}.
\]
We infer that the restrictions of $(\ssf_D,\cdot)-\lambda!_D$ and $(\ssf_{D'},\cdot)-\lambda!_{D'}$ to the intersection $\P_{D}'\cap \P_{D'}'$ coincide. This implies that the restrictions of the $(\ssf_D,\cdot)-\lambda!_D$ to $\P'_D$ glue together to a function $\Phi$ on the union of the $\P_D'$. So~\ref{seq:3} follows.

We now prove~\ref{seq:4}. Obviously, by definition, $\Phi\rest{\P_D'}$ is the restriction of an affine function. Since the union of the $\P_D'$ for $D\in \stable!_G^{^\circ}(\cl_{\valgroup})$ is a convex domain in the hyperplane $H\subset\R^V$, it will be enough to show that the function $\Phi$ is strictly convex around $\P_D\cap\P_D'$ for each $D,D'\in \stable!_G^{^\circ}(\cl_{\valgroup})$ whose associated polytopes share a common facet. Here, strict convexity is understood in the polyhedral geometric context, in the sense that there exists an affine function $\phi$ on $\R^V$ which coincides with $\Phi$ on the intersection $\P_D'\cap \P_{D'}'$, but $\bigl(\Phi-\phi\bigr) \rest{\P_D'\setminus \P_{D'}'}>0$ and $\bigl(\Phi-\phi\bigr) \rest{\P_{D'}'\setminus \P_{D}'}>0$. 

As in the proof of \ref{seq:3}, we may assume the sequence $\ssD_0, \dots, \ssD_{p-1}, \ssD_p=D$ chosen for $D$ satisfies $D'=\ssD_{p-1}$. Then, the base polytope of 
$\ssub{(\ssrho_{D'})}!_{\pi_p}$ is the intersection $\P_D'\cap \P_{D'}'$, and 
$\ssf_D-\ssf_{D'}=b_p\one!_{F_p}$ and $\lambda!_{D}-\lambda!_{D'}=b_p\ssrho_{D}(\ssF_p)$. We claim
\[\phi\coloneq \frac12\big((\ssf_{D},\cdot) -\lambda!_D+(\ssf_{D'},\cdot)-\lambda!_{D'}\big)\]
works as desired. Indeed,
\[\phi\rest{\P_D'\cap \P_{D'}'} = \Phi\rest{\P_D'\cap\P_{D'}'},\]
given that $(\ssf_D,\cdot) -\lambda!_D$ and $(\ssf_{D'},\cdot) -\lambda!_{D'}$ coincide on the facet. 
In addition,
\begin{align*}
\bigl(\Phi-\phi\bigr) \rest{\P_{D}'\setminus \P_{D'}'} =& \frac12
\big((\ssf_{D},\cdot)-(\ssf_{D'},\cdot)-\lambda!_D+\lambda!_{D'}\big)
\rest{\P_D'\setminus \P_{D'}'}\\
=&\frac 12\big((b_p\one!_{F_p},\cdot)-
b_p\ssrho_{D}(\ssF_p)\big)\rest{\P_D'\setminus \P_{D'}'},
\end{align*}
which is positive because 
$q(\ssF_p)>\ssrho_D(\ssF_p)$ for each $q\in \P_D'\setminus \P_{D'}'$. An analogous argument shows that 
$\bigl(\Phi-\phi\bigr) \rest{\P_{D'}'\setminus \P_{D}'}>0$
as well.

It remains only to prove that the restriction of $(\ssf_D,\cdot)-\lambda!_D$ to the hyperplane $H\subset\R^V$ does not depend on the choice of the sequence $\ssD_i$. Given two sequences of divisors $\ssD_0, \dots, \ssD_r=D$ and $\ssD_{r+q}=\ssD_0, \ssD_{r+q-1}, \dots, \ssD_{r+1}, \ssD_{r}=D$ from $\ssD_0$ to $\ssD$, we merge them together to get a sequence $\ssD_0, \dots, \ssD_r=D, \ssD_{r+1}, \dots, \ssD_{r+q} = \ssD_0$. This way, it will be enough to prove that along any \emph{cycle} $\ssD_0, \ssD_1,  \dots, \ssD_p=\ssD_0$ with consecutive associated polytopes sharing a facet, if we pursue as in Definition~ \eqref{eq:regularity}, we get that $(\ssf_{D_0},\cdot)-\lambda!_{D_0}$ 
coincides with $(\ssf_{D_p},\cdot)-\lambda!_{D_p}$ on $H$. Note that this cycle represents an actual cycle in the \emph{dual graph} of the tiling. This dual graph has vertices in bijection with the $\P_D'$ for $D\in \stable!_G^{^\circ}(\cl_{\valgroup})$, with an edge connecting a pair of vertices if the corresponding tiles share a facet. 

Cutting the cycle into smaller cycles in the dual graph, we can reduce to the case of a cycle $\ssD_0, \dots, \ssD_{p-1}, \ssD_p=\ssD_0$ as above where the intersection of all the polytopes $\P_{D_i}'$, $i=0, \dots, p$, is nonempty. Choose a point $q \in \bigcap_{i=0}^p\P_{D_i}'$. We show that $(\ssf_{D_0},\cdot)-\lambda!_{D_0}$ 
coincides with $(\ssf_{D_p},\cdot)-\lambda!_{D_p}$ on $H$.

Since $\ssD_p=\ssD_0$ and $\ssD_p = \ssD_0 + \div(\ssf_{D_p}-\ssf_{D_0};\ssD_0)$, we have that $\ssf_{D_p}=\ssf_{D_0}+c\one!_{V}$ for some constant $c\in\R$. Thus, $\ssf_{D_p}-\ssf_{D_0}$ is constant on $H$. Since $q\in\P_{D_i}'$ for every $i$, that constant can be determined as
\[
(\ssf_{D_p}-\ssf_{D_0},q)
=\sum_{i=1}^p \ssub{b}!_{i}(\ssub{\one}!_{F_i},q)=
\sum_{i=1}^p \ssub{b}!_{i}\,q(\ssF_i)=\sum_{i=1}^p \ssub{b}!_{i}\ssrho_{D_i}(\ssF_i)=\lambda!_{D_p}-\lambda!_{D_0}.
\]
We infer that $(\ssf_{D_p},\cdot)-\lambda!_{_p}$ coincides with $(\ssf_{D_0},\cdot)-\lambda!_{D_0}$ on $H$. The proof of the theorem is complete.
\end{proof}

\begin{remark}
 We note that in equal characteristic, the regularity of the tilings in Theorem~\ref{thm:reduction-admissible-tiling} follows as well from Theorem~\ref{thm:comparison-chow-quotient}, proved in the next section, and the regularity of the tilings associated to Chow quotients, proved by Lafforgue in~\cite{Laf99}. The above given proof handles all the tilings at the same time, and is needed as well for treating those coming from Theorem~\ref{thm:mixed-tiling}.  
\end{remark}

\subsection{Chow quotients} \label{sec:chow-quotients} We formulate a connection between the tiling appearing in Section~\ref{sec:tilings-simplex-mixed} and the works by Kapranov~\cite{Kap-chow}, Lafforgue~\cite{Laf99} and Giansiracusa and Wu~\cite{GW22}. 

We assume in this section that the valuation ring $\varR$ contains a copy of the residue field $\k$, so $\k$ and the valued field $\K$ have the same characteristic.

Recall that $\varH$ is a $\K$-vector space of dimension $r+1$ of the space of global sections of the line bundle $\varL=\mathcal O(\varD)$, that we view in $\K(\varX)$. For each 
$G$-admissible rational function $\ssh \in \adfun_{G}(\tau(\varD);\valgroup) \subset \Rat(\Gamma)$, we have an associated $\varR$-submodule 
$\varM_{\ssh}\subseteq\varH$, which is the collection of functions $\varf\in\varH$ such that 
$\trop(\varf)(v)\geq\ssh(v)$ for each $v\in V$. The inclusion $\varM_{\ssh}\hookrightarrow\varH$ induces a natural isomorphism $\varM_{\ssh}\otimes_{\varR}\K\to\varH$.

Let as before $t \mapsto \ssa_t$ be a section of the valuation map 
$\val\colon\K\setminus \{0\}\to\Lambda$. For each $\ssh_1,\ssh_2\in\adfun_{G}(\tau(\varD);\valgroup)$, put
\[
\lambda=\lambda!^{h_1}_{h_2}\coloneqq\min_{v\in V}(\ssh_1(v)-\ssh_2(v)).
\]
We have an inclusion $\ssa_{\lambda}^{-1}\varM_{h_1}\subseteq\varM_{h_2}$, and the multiplication by $\ssa_{\lambda}^{-1}$ induces a 
$\K$-linear automorphism of $\varH$ that restricts to 
an $\varR$-module homomorphism $\sspsi^{h_1}_{h_2}\colon\varM_{h_1}\to\varM_{h_2}$. Also, it induces the $\k$-linear map 
$\ssphi^{h_1}_{h_2}\colon\ssW_{h_1}\to\ssW_{h_2}$
considered in Section~\ref{sec:tilings-simplex-mixed}. Since $\epsilon\mapsto a_{\epsilon}$ is a group homomorphism, $\ssphi^{h_1}_{h_2}$ is given by multiplication by $\varphi\in\k^V$ with $\varphi_v$ equal to $1$ if $h_1(v)-h_2(v)=\lambda$, and to $0$ otherwise.

Recall that for each $v\in V$, there is $\ssh_v\in\Rat(\Gamma)$ such that the divisor $\tau(\varD)+\div(\ssh_v)$ is $v$-reduced. Furthermore, $\ssh_v\in\adfun_{G}(\tau(\varD);\valgroup)$ by Proposition~\ref{prop:v-reduced-adm}, and the projection map $\ssW_{h_v}\to\VS_v$ is an isomorphism by Proposition~\ref{prop:v-reduced-iso}. Now, the projection to $\VS_v$ is compatible with $\ssphi^{h}_{h_v}$ for each $\ssh\in\adfun_{G}(\tau(\varD);\valgroup)$, that is, 
$\theta_v\ssphi^{h}_{h_v}=\theta_v$ for each $v\in V$. 
The $\ssphi^{\ssh}_{h_v}$ for $v\in V$ induce an injection
\[
\ssW_{\ssh}\longrightarrow\bigoplus_{v\in V}
\ssW_{h_v}.
\]
This injection is induced by tensoring by $\k$ the $\varR$-module homomorphisms $\psi^{h}_{h_v}$. Hence, the $\psi^{h}_{h_v}$ induce an injection
\[
\varM_{\ssh}\longrightarrow\bigoplus_{v\in V}
\varM_{h_v}
\]
with $\varR$-flat cokernel. Tensoring by $\K$, and using the natural isomorphisms $\varM_{\ssh}\otimes_{\varR}\K\to\varH$ and $\varM_{h_v}\otimes_{\varR}\K\to\varH$, we obtain an embedding $\varH\to\oplus_{v\in V}\varH$ whose image is in the orbit of the diagonal by the action of the torus $\varG^V_{\varm}$.

Assume $\varR$ contains a copy of the field $\k$. Choosing an ordered basis for $\varM_{h_v}$ for each $v\in V$, we obtain a morphism 
\[
\beta_h\colon\text{Spec}(\varR)\longrightarrow \varGr,\quad\text{where}\quad \varGr \coloneqq\mathrm{Grass}\big(r+1,\underbrace{\k^{r+1}\oplus \dots \oplus \k^{r+1}}_{|V| \textrm{ copies }}\big),
\]
the Grassmannian parametrizing subspaces of $\k^{(r+1)|V|}$ of dimension $r+1$. The torus 
$\varG^V_{\varm}$ acts on the Grassmannian, each factor acting on the corresponding copy of $\k^{r+1}$. The action is a restriction of the action by the (maximum) torus of dimension $(r+1)|V|$. The morphism $\beta_h$ sends the generic point to a point on the orbit of the diagonal subspace, the image under the diagonal map $\delta\:\varH \hookrightarrow \varH\oplus\dots\oplus\varH$. 

Kapranov~\cite{Kap-chow} considered only the action by the maximum torus on the Grassmannian in his seminal work. Lafforgue~\cite{Laf99} studied exactly the case of interest to us, however, we find it more natural to formulate the connection using the more recent work by Giansiracusa and Wu~\cite{GW22}. This work extends Kapranov's study \cite{Kap-chow} to actions by subtori, and as in \emph{loc.~cit.}, describes the Chow quotient of the action. It covers as well Lafforgue's work~\cite{Laf99}. 

The action of the diagonal $\varG_{\varm}\to\varG^V_{\varm}$ on the Grassmannian is trivial, so we may consider the action by the quotient $\varT\coloneqq\varG^V_{\varm}/\varG_{\varm}$. The diagonal $\delta(\varH)$ is parameterized by a $\K$-point on the open subscheme $\varY$ of the Grassmannian where the action is free, and there is an inclusion $\rquot{\varY}{\varT}\hookrightarrow \mathrm{Chow}(\varGr)$ of the quotient $\rquot{\varY}{\varT}$ to the Chow variety of $\varGr$, taking an orbit to its closure, viewed as a prime cycle of 
$\varGr$. The \emph{Chow quotient} is simply the closure of $\rquot{\varY}{\varT}$ in $\mathrm{Chow}(\varGr)$. It parameterizes algebraic cycles, that Giansiracusa and Wu claim to be reduced in \cite[Prop.~2.2]{GW22}. 

Moreover, they show that a cycle $\sum_im_i\varZ_i$ is parameterized by the Chow quotient only if $m_i=1$ for every $i$, and the prime cycles $\varZ_i$ are closures of $\varT$-orbits of $(r+1)$-dimensional subspaces $L_i\subseteq\k^{(r+1)|V|}$ whose associated polytopes $\Q_{L_i}$ (see Section~\ref{sec:modular-subspace}) are full-dimensional and form a polyhedral decomposition of $\Delta_{r+1}$. 

In our case, the diagonal $\delta(\varH)$ gives a $\K$-point of $\varGr$ whose orbit closure corresponds to a subscheme of $\varGr\times_{\k}\K$. Its closure $\mathcal Z\subseteq \varGr\times_{\k}\mathrm{Spec}(\varR)$ has closed fiber $\varZ$ over $\k$ whose associated cycle is of the form $\sum_i\varZ_i$ as described above.

Let $\cl$ be the class of $\tau(\varD)$ in $\Pic^d(\Gamma)$ and $\cl_{\valgroup}\subseteq\cl$ the subset of $\Lambda$-rational divisors. Denote by $\stable!_G(\cl_{\valgroup})$ the set of $G$-admissible divisors $D$ in $\cl_{\valgroup}$. 

\begin{thm}\label{thm:comparison-chow-quotient} Notation as above, the $\ssW_{D}$ for $D\in \stable!_G(\cl_{\valgroup})$ are in the orbit closures of the $L_i$. Inversely, each $L_i$ is in the orbit of $\ssW_D$ for some $D\in \stable!_G(\cl_{\valgroup})$. In particular, the tiling by the $\Q_D$ in Theorem~\ref{thm:reduction-admissible-tiling} coincides with the tiling by the $\Q_{L_i}$ and their faces.
\end{thm}
\begin{proof} Each map $\beta_h$ sends its generic point to a $\K$-point on the orbit of the diagonal, whence $\beta_h$ factors through $\mathcal Z$. It follows that each $\ssW_D$ for $D\in \stable!_G(\cl_{\valgroup})$ lies in the orbit closure of some $L_i$. This proves the first statement.

Now, for the second statement, Theorem~\ref{thm:reduction-admissible-tiling} implies that for each $i$, there is $D\in \stable!_G(\cl_{\valgroup})$ such that $\Q_{D}$ is full-dimensional and its interior intersects the interior of $\Q_{L_i}$. On the other hand, as pointed out, $\ssW_D$ is in the orbit closure of $L_j$ for some $j$. Since $\Q_{D}$ is full-dimensional, $\ssW_D$ must be in the orbit of $L_j$. But then $\Q_{D}=\Q_{L_j}$. Since $\Q_{D}$ and $\Q_{L_i}$ have an interior point in common, we must have $j=i$. Thus, $L_i$ is in the orbit of $\ssW_D$. 

The last statement follows.
\end{proof}

\subsection{Reduced divisors in the tropicalization of linear series} \label{sec:reduced-divisors} Notation as in Sections~\ref{sec:tropicalization}~and~\ref{sec:tilings-simplex-mixed}, let $\mathrm{M}$ be the tropicalization of $\varH$, that is the set of functions $\trop(\varf)$, $\varf \in \varH$. Endowed with the operation of pointwise minimum and addition by constants, $\mathrm{M}$ is a subsemimodule of $\Rat(\Gamma)$.

For each vertex $v\in V$,  define
\[\ssf_v^{\mathrm M} \coloneqq \inf_{\substack{h\in \mathrm{M}\\h(v)=0} } h.\]
Since $\mathrm{M}$ is closed in $\Rat(\Gamma)$ for the topology of pointwise convergence, we 
have $\ssf_v^{\mathrm{M}}=\trop(\varf_v)$ for an element $\varf_v$ in $\varH$. 
Consider the restriction $\ssf_v^{\mathrm{M}}\rest{V}$ and denote by $\ssh_v^{\mathrm{M}}$ its $G$-admissible extension with respect to $\tau(\varD)$. Finally, set $\ssD_v^{\mathrm{M}} \coloneqq  \div(\ssh_v^{\mathrm{M}}) +\tau(\varD)$, the corresponding $G$-admissible divisor. Since $\ssf_v^{\mathrm{M}}=\trop(\varf_v)$, we have that 
$\ssD_v^{\mathrm{M}}\in\stable!_G(\cl_{\valgroup})$.

For each $v\in V$, let $q_v$ be the vertex of $\Delta_{r+1}$ with $q_v(v)=r+1$, and hence $q_v(u)=0$ for each $u\in V$ distinct from $v$. 
Let $\Q_v \coloneqq \Q_{D_v^{\mathrm{M}}}$.

\begin{thm}\label{thm:reduced-divisors} Notation as above, the supermodular function $\nu!_{D_v^{\mathrm{M}}}$ is simple. The polytope $\Q_v$ is the unique full-dimensional polytope in the tiling of the simplex $\Delta_{r+1}$ that contains the vertex $q_v$.
\end{thm}

\begin{proof} Every function $\varf \in \varH$ with $\trop(\varf)(v)=0$ verifies $\trop(\varf) \geq \ssh_v^{\mathrm{M}}$. This implies that the projection 
\[\ssW_{h_v^{\mathrm{M}}} \to \VS_v\]
is an isomorphism. Thus, $\ssW_{h_v^{\mathrm{M}}} \simeq \ssW_{h^{\mathrm{M}}_v,I}$,  and so $\ssW_{h_V}^{I^c} =0$, or equivalently, $\nu!_{D_v^{\mathrm{M}}}(\ssI^c)=0$ and $\nu!_{D_v^{\mathrm{M}}}^*(I) =r+1$ for each subset $I\subseteq V$ containing $v$.

We first show that $\Q_v$ contains $q_v$. We need to show the inequalities 
\[\nu!_{D_v^{\mathrm{M}}}(I) \leq q_v(I) \leq \nu!_{D_v^{\mathrm{M}}}^*(I) \quad \forall I \subseteq V.\]
But, if $I$ contains $v$, then $q_v(I)=r+1 = \nu!_{D_v^{\mathrm{M}}}^*(I)$. And if $I$ does not contain $v$, then $q_v(I)=0=\nu!_{D_v^{\mathrm{M}}}(I)$. In any case, the inequalities hold. This proves $q_v \in \Q_v$.

We show that $\Q_v$ is full-dimensional. We need to prove that $\nu!_{D_v^{\mathrm{M}}}$ is simple. Assume for the sake of a contradiction that $(\nu!_{D_v^{\mathrm{M}}})_\pi=\nu!_{D_v^{\mathrm{M}}}$ for a bipartition $\pi=(I, J)$ of $V$. Without loss of generality, assume that $v\in I$. 
Then, 
\[
\nu!_{D_v^{\mathrm{M}}}^*(J)=r+1-\nu!_{D_v^{\mathrm{M}}}(I)=\nu!_{D_v^{\mathrm{M}}}(J)=0,
\]
that is, $\ssW_{h_v^{\mathrm{M}},J}=(0)$. On the other hand, since $\ssh_v^{\mathrm{M}}\rest{V} = \trop(\varf_v)\rest{V}$, the reduction of $\varf_v$ at a vertex $u$ of $J$ gives a nonzero element of $\VS_u$, leading to a contradiction.
\end{proof}

\begin{remark}
Note that the divisor $\div(\ssf_v^{\mathrm{M}})+\tau(\varD)$ is $v$-reduced with respect to $\mathrm{M}$ in the terminology of~\cite{AG22}, which defines reduced divisors with respect to semimodules. (This is a relative notion of reducedness, and in general differs from that given in Section~\ref{sec:v-reduced-1}.) We do not know whether the equality $\ssh_v^{\mathrm{M}}=\ssf_v^{\mathrm{M}}$ holds in general.
\end{remark}

\subsection{Fundamental collection, effectivity and finite generation property} \label{sec:finiteness-lls} Consider the setup of Sections~\ref{sec:tropicalization}~and~\ref{sec:tilings-simplex-mixed}.  
The \emph{fundamental collection}  associated to the pair $(\varD, \varH)$ is the collection $\FC_{(\varD, \varH)}$ of spaces $\ssW_D$ for $D \in \stable_G(\cl_{\valgroup})$ with full-dimensional $\Q_D$, equivalently, with $\nu!_D$ simple: 
\[
\FC_{(\varD, \varH)} \coloneqq \left\{\ssW_D \st  D \in \stable_G(\cl_{\valgroup}) \, \textrm{  with }\Q_D \textrm{ full-dimensional}\right\}.
\]

Let $\cl$ be the linear equivalence class of $\tau(\varD)$, and $D \in \stable_G(\cl_{\valgroup})$ be a $G$-admissible $\valgroup$-rational divisor such that $\Q_D$ is full-dimensional, or equivalently, $\nu!_D$ is simple. We can write $D=\div(h)+\tau(\varD)$ for $h\in \adfun_G(\tau(\varD);\valgroup)$.   Since $\Q_D$ is full-dimensional, $\nu!^*_D(v) \neq 0$ for all vertices $v$ of $G$. This means that there exists for each $v\in V$ an element $\varf_v \in \varM_h$ with $\trop(\varf_v)(v)=h(v)$. Taking a generic linear combination of the $\varf_v$ gives an element $\varf \in \varM_h$ with $\trop(\varf)(v)=h(v)$ for all $v$.  The Specialization Lemma \cite[Thm.~4.5]{AB15} yields that 
$\div(\trop(\varf))+\tau(\varD)$ is effective.  The domination property for $G$-admissible divisors, Lemma~\ref{lem:varH}, implies that $\trop(\varf)\geq h$. From this inequality, since $\trop(\varf)\rest{V}=h\rest{V}$, we deduce that $\div(\trop(\varf))(v)\leq \div(h)(v)$ for each $v\in V$, leading to inequalities 
\[D(v) \geq \div(\trop(\varf))(v)+\tau(\varD)(v) \geq 0\]
for $v\in V$. Since $D$ is $G$-admissible, we deduce that it is effective.

\begin{thm}[Finite generation property]\label{thm:finiteness-lls} Each $G$-admissible divisor $D \in \stable_G(\cl_{\valgroup})$ whose polytope $\Q_D$ is full-dimensional is necessarily effective.  Furthermore, the fundamental collection $\FC_{(\varD, \varH)}$ of spaces $\ssW_D$ with $\Q_D$ full-dimensional generates all the other spaces in the following sense: For each $G$-admissible divisor $D'\in \cl_{\valgroup}$, the space $\ssW_{D'}$ is generated by the images of the maps $\ssW_D \to \ssW_{D'}$ for $D\in \mathcal S$.
\end{thm}

\begin{proof} Let $f$ be a nonzero element of $\ssW_D$. It will be enough to show that $f$ is generated by the images of the maps $\ssW_{D'} \to \ssW_D$, for $G$-admissible divisors $D' \in \cl_{\valgroup}$  with $\nu!_{D'}$ simple. We prove this statement proceeding by induction on the codimension of the polytope $\Q_D$. 

If the codimension is zero, then the polytope $\Q_D$ is full-dimensional, that is, $\nu!_D$ is simple. We take $D'=D$, and there is nothing to prove. 

Assume the codimension is positive. Then, $\ssW_D$ is $\pi$-split for some bipartition $\pi = (I,J)$, and we can write $f = f_1+f_2$ with $f_1 \in \ssW_D^I$ and $\ssW_D^J$. It will be enough to treat the case where $f\in \ssW_D^I$. By Proposition~\ref{prop:complete-admissible-lls}, the function $\nu!_D$ is the $\pi$-splitting of $\nu'=\nu!_{D'}$ for a $G$-admissible divisor $D'\in \cl_{\valgroup}$ such that $\nu!_{D'} \neq \nu!_D$. Moreover, using the notation of Section~\ref{sec:complete-admissible-lls}, we have $D' = \ssD_b$ for $b=\ts_D(I)$ the $\pi$-splitting threshold of $\nu!_D$, as in Definition~\ref{defi:threshold}. By  Proposition~\ref{prop:vanishing}, the pair $(\ssW_D, \ssW_{D'})$ is exact. We infer that $\ssW_{D'}\to\ssW_D$ factors through an isomorphism $\ssW_{D',I}\to\ssW_D^I$. Then, there is $f'\in\ssW_{D'}$ whose image in $\ssW_D$ is $f$. The polytope $\Q_{D'}$ has codimension strictly smaller than $\Q_D$, and so  by the induction hypothesis, $f'$ is generated by the images of the maps $\ssW_{D''} \to \ssW_{D'}$, for  $G$-admissible divisors $D'' \in \cl_{\valgroup}$ with $\nu!_{D''}$ simple.  Then, $\ssf$ is generated by the images of the compositions $\ssW_{D''} \to \ssW_{D'} \to \ssW_D$. The compositions which are zero may be discarded. Those which are nonzero are multiples of the maps  $\ssW_{D''} \to \ssW_{D}$ by Proposition~\ref{prop:compoh}. We infer that $f$ is generated by the images of the maps $\ssW_{D''} \to \ssW_{D}$, for $G$-admissible divisors $D'' \in \cl_{\valgroup}$ with $\nu!_{D''}$ simple. The theorem follows.
\end{proof}


\section{Complementary results} \label{sec:discussion}
We discuss a few complementary results related to the content of the paper.

\subsection{Admissible sum and periodicity of the tiling by semistability polytopes} \label{sec:admissible-sum} Let $G$ be a connected graph, $\ell\colon E \to \R_{+}$ an edge length function, and $\Gamma$ the corresponding metric graph. 
Let $\valgroup\subseteq\R$ be an additive subgroup containing the $\ell_e$.

Given two $G$-admissible divisors $\ssD_1$ and $\ssD_2$ on $\Gamma$, we define their \emph{$G$-admissible sum}, denoted $\ssD_1+_\ell \ssD_2$, by
\[
\ssD_1+_\ell \ssD_2\coloneqq \ssD_1+\ssD_2+\div_\ell(0;\ssD_1+\ssD_2).
\]
Then, $\ssD_1+_\ell \ssD_2$ is a $G$-admissible divisor linearly equivalent to $\ssD_1+\ssD_2$. If $D_1$ and $D_2$ are $\valgroup$-rational, then so is $\ssD_1+_\ell \ssD_2$.

We have the following theorem.

\begin{thm} The following statements hold:
\begin{enumerate}
    \item The $G$-admissible sum gives the set of $G$-admissible $\valgroup$-rational divisors on $\Gamma$ an Abelian group structure.
    \item The $G$-admissible $\valgroup$-rational divisors on $\Gamma$ linearly equivalent to $0$ form a subgroup.
    \item For each class $\cl\in\Pic^d(X)$ such that $\cl_{\valgroup}\neq\emptyset$, the set $\stable!_G(\cl_{\valgroup})$ is a coset for that subgroup.
\end{enumerate}
\end{thm}

\begin{proof} The statements follow by direct verification using Proposition~\ref{prop:adm} and Theorem~\ref{thm:admissible2}. We omit the details.
\end{proof}

Note that if $\ssD_1$ is $G$-admissible and $\ssD_2$ has support on the set of vertices of $G$, then $\ssD_1+\ssD_2$ is $G$-admissible and we have $\div_\ell(0, \ssD_1+\ssD_2)=0$. This implies that 
\[\ssD_1+_\ell\ssD_2 = \ssD_1+\ssD_2.\]

Using the above theorem, we obtain the following periodicity result for the tilings given by Theorem~\ref{thm:semistability-tiling}. 

Let $L$ be the set of all $G$-admissible divisors linearly equivalent to the zero divisor which are supported on the set of vertices of $G$. It is in fact a group whose elements are of the form $\div_\ell(f;0)$ for $f\colon V\to \R$  that has the property that for each $e=uv\in\E$, $\frac{f(u)-f(v)}{\ell_e}$ is an integer. Also, notice that $\P_D\subset H_0$ for each $D\in L$, and that the centers $\ssq_D$ of $\P_D$ for $D\in L$ form a lattice in $H_0$. Moreover, the assignment $D\mapsto q_D$ is a group isomorphism, giving an embedding of $L$ in $H_0$.

\begin{thm} Let $\cl \in \Pic^d(\Gamma)$ be a divisor class of degree $d$ with $\cl_{\valgroup}\neq\emptyset$. Then, the tiling by $\P_D$ for $D\in \stable_G(\cl_{\valgroup})$ is $L$-periodic, i.e., it is invariant under addition of the center $\ssq_D$ of $\P_D$ for each $D\in L$.

In particular, if the ratios of edge lengths are all rational, then the tiling by $\P_D$ is 
periodic.
\end{thm}

\begin{proof} By the discussion preceding the theorem, for each divisor $D\in \stable_G(\cl_{\valgroup})$ and $D'\in L$, we have $D+D'\in \stable_G(\cl_{\valgroup})$. Also, $\ssG_{D+D'}=\ssG_D$ and $\ssq_{D+D'}=\ssq_D+\ssq_{D'}$.

By Theorem~\ref{thm:voronoi-semistability}, each $\P_D$ for $D\in \stable_G(\cl_{\valgroup})$ is the translation by $\ssq_D$ of the Voronoi cell of the graph $\ssG_D$ associated to $D$, if $\ssG_D$ is connected, or of a product of Voronoi cells if not; see Remark~\ref{rmk:voronoi-semistability}.  We infer that the tiling by $\P_D$ for $D\in \stable_G(\cl_{\valgroup})$, is invariant under addition by the $\ssq_D$ for $D\in L$. This proves the first claim.

If the ratios of edge lengths are all rational, then the lattice formed by the $q_D$ for $D\in L$ has full rank in $H_0$, and we conclude that the tiling is periodic. 
\end{proof}

There is much more to say about periodicity properties of the tilings produced by Theorem~\ref{thm:semistability-tiling}. For example, by keeping track as well of the frequency of appearance of each edge in the spanning subgraphs $\ssG_D$ of $G$ associated to the divisors $D$, we can recover the metric graph $\Gamma$ from the tiling. A  more thorough discussion of this topic, important for the sequel work, is beyond the scope of this paper and is postponed to a future work.

\subsection{Bricks and polymatroids}
Our Theorems~\ref{thm:semistability-polytope} and~\ref{thm:voronoi-semistability} give a description of the polytopes $\P_D$ associated to $G$-admissible divisors $D$ in a given linear equivalence class. 

Figure~\ref{fig:mixed-tiling} depicts an example of tiling issued from Theorem~\ref{thm:mixed-tiling}. It would be interesting to have a combinatorial description of the polytopes that show up in Theorems~\ref{thm:reduction-admissible-tiling} and~\ref{thm:mixed-tiling}.

\begin{question} Describe the combinatorics of the polytopes $\Q_D$, for $G$-admissible $\valgroup$-rational divisors $D$ linearly equivalent to $\tau(\varD)$.
\end{question}

Here is a partial result in this direction, which turns out to be important in our forthcoming work~\cite{AEG23}. For each element $\beta \in \Delta_{r+1}$, define the polytope $\B_\beta$ as the set of all points $q \in H_{r+1}$ that verify the inequalities 
\[\lfloor \beta(I)\rfloor \leq q(I) \leq \lceil \beta(I)\rceil\quad \forall \,\, I \subseteq V.\]
We refer to $\B_\beta$ as the \emph{brick} associated to $\beta$. We have the following theorem.

\begin{thm} The collection of bricks $\B_\beta$ for $\beta \in \Delta_{r+1}$ gives a tiling of $\Delta_{r+1}$. Furthermore, the base polytope $\Q$ of a polymatroid with integral vertices  in $\Delta_{r+1}$ is tiled by the set of bricks included in $\Q$. In particular, the tilings in Theorem~\ref{thm:reduction-admissible-tiling} are coarsenings of the tiling by bricks.
\end{thm}

\begin{proof} Each point $\beta \in \Delta_{r+1}$ obviously belongs to $\B_\beta$. Furthermore, since the supermodular function $\nu$ defining a base polytope $\Q$  is integer-valued, if $\beta\in \Q$, then $\B_\beta$ is entirely in $\Q$. This proves that the union of $\B_\beta$ for $\beta \in \Delta_{r+1}$ is the full simplex, and the union for those $\beta\in \Q$ is $\Q$. Finally, for $\beta, \beta' \in \Delta_{r+1}$, if nonempty, the intersection $\B_{\beta}\cap\B_{\beta'}$ is the set of $q\in\B_{\beta}$ satisfying the additional condition 
$\lfloor \beta'(I)\rfloor \leq q(I) \leq \lceil \beta'(I)\rceil$ for each $I\subseteq V$. Depending on $\beta'$, each such condition amounts to no additional condition, to $q(I)=\lceil \beta(I)\rceil$, or to $q(I)=\lfloor \beta(I)\rfloor$. In any case, $\B_{\beta}\cap\B_{\beta'}$ is a face of $\B_{\beta}$. The theorem follows.
\end{proof}

\bibliographystyle{alpha}
\bibliography{bibliography}

\end{document}